\documentclass[10pt,a4paper]{article}
 \usepackage{enumitem}
 \usepackage{color}
\newcommand{\Sta}{\pmb{\mathbb{S}}}

\newcommand{\ztn}{\vartheta}
\newcommand{\argmin}{{\rm argmin}}
 \newcommand{\Mk}{M^\ztn_{\ztn-k}}
\newcommand{\Mh}{M^{\ztn+h}_\ztn}
 \newcommand{\Rk}{R^{\ztn-k}_\ztn}
\newcommand{\Rh}{R^\ztn_{\ztn+h}}
 \newcommand{\Gzt}{\Gamma_\zt}
  \newcommand{\Lzt}{\Lambda_\zt}
  
  \usepackage{amssymb,amsmath,exscale}
  \usepackage{bm}

    \usepackage{bbm}
 \usepackage[utf8]{inputenc} 
\usepackage{makeidx}
\usepackage{color}
  \usepackage{bbm}
  \usepackage{bm}
\usepackage{amsmath} 
\usepackage{amsfonts}
\usepackage{amssymb}
   \usepackage{amssymb,amsmath,exscale}
 
\usepackage{epsfig}
\usepackage{graphics}

 \newcommand{\Dom}{{\rm Dom}}

\newtheorem{Theorem}{Theorem}
\newtheorem{Corollary}[Theorem]{Corollary}
\newtheorem{Lemma}[Theorem]{Lemma}

\newtheorem{Remark}[Theorem]{Remark}
\newtheorem{Definition}[Theorem]{Definition}

\newtheorem{Assumption}[Theorem]{Assumption}

\newcommand{\zaa}{\alpha}
 
\newcommand{\zg}{\gamma}

 \newcommand{\ZEP}{\epsilon}
\newcommand{\ZSI}{\sigma}

 \newcommand{\ZOMq}{\Omega} 
 \newcommand{\zthe}{\theta}
\newcommand{\zt}{\tau}

\newcommand{\zzr}{\mathbb{R}}

\newcommand{\intT}{\int_0^T}
\newcommand{\intt}{\int_0^t}

\newcommand{\ZCD}{(\cdot)} 
\newcommand{\zdiaform}{\mbox{~~\zdia}}

 \newcommand{\ZR}{\rangle}
\newcommand{\ZL}{\langle}
 \newcommand{\zdia}{~~\rule{1mm}{2mm}\par\medskip}
 
\newcommand{\ZIN}{\infty}
\newcommand{\zProof}{{\bf\underbar{Proof}.}\ }
 \newcommand{\ZD}{\;\mbox{\rm d}}
 
 \newcommand{\ZLA}{\label}

\author{
L. Pandolfi\thanks{Dipartimento di Scienze Matematiche ``Giuseppe Luigi Lagrange'', Politecnico di Torino, Corso Duca degli Abruzzi 24, 10129 Torino, Italy (retired) (luciano.pandolfi@formerfaculty.polito.it)}
}
\title{The regulator problem for the   wave equation with high internal damping controlled on the boundary:  a new look via systems with memory}

\begin{document}

\maketitle

\begin{abstract}
We study the quadratic regulator problem on a finite time horizon  for the  wave equation 
with high internal damping controlled on the boundary by square integrable controls. The approach in this paper transforms  the wave equation 
with high internal damping to an equation with persistent memory controlled on the boundary.

One of the results   of this paper is the introduction of a state space which is an extended Hilbert space, so a time dependent Hilbert   space.
We prove that the unique optimal control   can be represented as a feedback control  via a Riccati operator which solves a suitable version of the Riccati equation.   Both the feedback operator and the Riccati equation acts on such time dependent space.
The derivation of these main results requires a very precise analysis of the properties of the derivatives of the value function      and we find an explicit form for the derivative of the Riccati operator.
\end{abstract}

\medskip

\textbf{20020 Mathematical Subject Classification.} Primary 49N10. Secondary
45D05, 93A10.
 
\section{Introduction}

We study the following quadratic regulator problem on a   bounded interval of time $[0,T]$. Let $v(t;v_0,v_1,u)=v(t)=v(x,t)$ solve
\begin{subequations}
\begin{equation}
\ZLA{eq:IntroNOTAZ} 
v''=\Delta v+\Delta v' \quad t\in(0,T)\,,\ x\in\ZOMq\qquad\left\{\begin{array}
{l} v(x,0)=v_0\,,\ v'(x,0)=v_1\\
v(x,t)=u(x,t)\quad x\in \partial\ZOMq \,.
\end{array}\right.
\end{equation}
We intend that $\ZOMq$ is a bounded region with   smooth   boundary.
We associate to~(\ref{eq:IntroNOTAZ}) the following quadratic cost and optimization problem:
 
\begin{equation}
\ZLA{eq:costoWAVEdamp}
\left\{\begin{array}{l}
\displaystyle
J_0(v_0,v_1;u)=
\intT\left [
\|v(t;v_0,v_1,u)\|^2_H+\|u(t)\|^2_ U
\right ]\ZD t \\
\displaystyle
 \min_{u\in L^2(0,T;U)} J_0(v_0,v_1;u)
\end{array}\right.
\end{equation}
 \end{subequations}
where
\[
H=L^2(\ZOMq)\,,\quad U=L^2(\partial\ZOMq) \quad \mbox{(real Hilbert spaces)}\,.
\]
 
 A quadratic regulator problem for the system~(\ref{eq:IntroNOTAZ}) has been studied in~\cite{BucciTESI}. This paper assumes    that $u\in H^1(0,T;U)$ and that the cost functional penalizes also $\|u'(t)\|^2_U$. Problem~(\ref{eq:IntroNOTAZ})-(\ref{eq:costoWAVEdamp}) was then studied    with a   method different from the one we are going to present here. One of the reason of interest 
  of the present paper is  the method we propose, different from that used before and for this reason
 comparison  of the methods and comments on the literature are postponed to Sect.~\ref{Sect:COMPARISON}, when the rationale behind both the methods   can be appreciated. 
 Here we confine ourselves to roughly state that the main   idea previously  used    is to assume first that  $u\in H^1 $ and to decouple $u(0)$ from $u\ZCD$. At the expenses of introducing highly unbounded operators in the representation of the system, it is then possible to     consider $u(0)$ as an additional parameter of an optimization problem respect to $u\in L^2$. The way for doing that is explained in Sect.~\ref{Sect:COMPARISON}. Here instead we avoid the decoupling of $u(0)$ from $u\ZCD$ and the introduction of such highly unbounded operators at the expenses of reducing~(\ref{eq:IntroNOTAZ}) to a heath equation with persistent memory controlled on the boundary.
 
 We are not asserting superiority of one of the method over the other. We believe that both deserve to be known but we conjecture that the method we introduce here can be adapted to study the quadratic regulator problem for
 other classes of systems, in particular for distributed systems with persistent memory controlled on the boundary (Sect.~\ref{subsEXTENS} suggests possible extensions) and this conjecture is   a reason of   our interest in this method. 
   A second reason of interest is the description of a distributed system with memory as a system in a time dependent Hilbert space (as explained in~\cite{Pandolfi24EECTtrackingRd}, this time dependent Hilbert space is an extended Hilbert space according to the definition in~\cite{DolezalBOOKmonotoneExtendedHspaces}). 
 
 %%%%%%%%%%%%%%  
   Finally we mention
  that a quadratic regulator problem for a class of systems with memory in a Hilbert space and with \emph{distributed control} is studied in~\cite{acquistaBUCCI24}. This paper extends to Hilbert spaces the results in~\cite{PandolfiVOLTERRAieee18,Pandolfi24EECTtrackingRd} for systems with persistent memory in $\zzr^d$.

 \subsection{Summary of the results}
 
We outline the main steps used in this paper to study the quadratic regulator problem~(\ref{eq:IntroNOTAZ})-(\ref{eq:costoWAVEdamp}).

\subsubsection*{Step 1: we transform equation~(\ref{eq:IntroNOTAZ}) to an heat equation with persistent memory} To achieve this result we borrow an idea from~\cite{BucciPandolfiArxive,DelloroPATA17MGTeVISCO}. 

In Sect.~\ref{SECT:soluzione} we perform formal computations which transform Eq.~(\ref{eq:IntroNOTAZ})  to the  heat equation with persistent memory~(\ref{eq:laFORMintegDIv}) below. The solutions  of~(\ref{eq:laFORMintegDIv}) (hence of~(\ref{eq:IntroNOTAZ})) are  those of the Volterra integral equation~(\ref{eq:laFORMintegDIvCONuVOLTE}).
The goal of Sect.~\ref{SectCondiINIcondiEsoluzL2} is the justification of such formal computations and the definition of the solutions in a   specified weak sense  and in  suitable spaces which are identified by the following considerations.
The  solutions of Eq.~(\ref{eq:IntroNOTAZ}), i.e. of~(\ref{eq:laFORMintegDIv}), depend on the control $u$ and on the initial conditions $v(0)  $ and $v'(0) $. 
 The dependence on $v(0) $ involves also $\Delta v(0)$ but it turns out that   $\Delta v(0)$ and $v'(0)$ are collected in a certain vector $\hat y$ which is linearly independent of $v_0=v(0)$ since $v(0)$ and $v'(0)$ are independent. So,  we can study the solutions as functions of $v_0$, $\hat y$ and   $u$. 
It is proved in Sect.~\ref{SectCondiINIcondiEsoluzL2} that the formal computations can be justified when  
 \[
v_0=v(0)\in H^{1/2+\ZEP}(\ZOMq)\subseteq  \Dom(-A)^{1/4-\ZEP/2}\,,\quad \hat y\in (\Dom(-A)^{3/4-\ZEP/2})'
 \]
($A$ is the laplacian with homogeneous Dirichlet conditions in $L^2(\ZOMq)$). The control $u$ is square integrable and it turns out that the problem is well posed on these spaces. In particular, the solution $v$ of~(\ref{eq:laFORMintegDIv})  is an $H$-valued  square integrable function which   depends continuously on the data, so that the quadratic regulator problem can be studied.

The precise statements are in the theorems~\ref{teoREgoDiv} and~\ref{teoLasolUdebDeVolte}.

\subsubsection*{Step~2: Problem~(\ref{eq:IntroNOTAZ})-(\ref{eq:costoWAVEdamp}) is embedded in a larger family of quadratic regulator problems} 
First it is proved that Eq.~(\ref{eq:laFORMintegDIv}) is well posed also if $(v_0,\hat y)\in L^2(\ZOMq)\times (\Dom\,A)'$ and the control $u$ is square integrable  (statement~\ref{I1teoREgoDiv} of Theorem~\ref{teoREgoDiv}). 
Then problem~(\ref{eq:IntroNOTAZ})-(\ref{eq:costoWAVEdamp}) is embedded in a family of quadratic regulator problems, each one on $[\zt,T]$, any $\zt\in [0,T)$ (the problem is described by the equations~(\ref{eq:laEQsu(ztT)}) and~(\ref{eq:costoWAVEdampSU(ztT)})).   The initial conditions 
of~(\ref{eq:laEQsu(ztT)}) are two vectors, $\hat v_\zt\in H $ and  $\hat y_\zt\in (\Dom\,A)'$, and a function $ \xi_\zt\in L^2(0,\zt;H)$ which depends on the \emph{past memory,} on the interval $(0,\zt)$, of the system. The function $\xi_\zt\in L^2(0,\zt;H)$ needs not be a ``segment'' of trajectory: it can be arbitrarily be assigned.
 
Theorem~\ref{Teo:preREGOLARITAstato} states that the solutions of~(\ref{eq:laEQsu(ztT)})  is square integrable and depends continuously on  
\[
(\hat v_\zt, \xi_\zt,\hat  y_\zt,u)\in H\times L^2(0,\zt;H)\times (\Dom\,A)'\times L^2(\zt,T;U)\,.
\]

\subsubsection*{Step~3: the study of the quadratic regulator problem of Eq~(\ref{eq:laEQsu(ztT)}) with cost~(\ref{eq:costoWAVEdampSU(ztT)})}

Existence and unicity of the optimal control is proved by a standard variational method, which   represents the optimal control as the unique solution of a Fredholm integral equation (which is the weak form of a hamiltonian system). A consequence of this fact is the 
continuity of the optimal control which moreover depends continuously on the  data $(\hat v_\zt, \xi_\zt,\hat  y_\zt)$ (Theorem~\ref{TEOValue functionCONTI}).

\subsubsection*{Step~4: the \emph{state}
of the system and the feedback form of the optimal control} 
In this step we replace 
  the optimization problem in $L^2(\zt,T;U)$ (any $\zt\in [0,T)$)  by a family of optimization problems in $U $ parametrized by $t\in [\zt,T]$. This is done in Sect.~\ref{SecDissiIneq} by using the standard Bellman principle and  the dissipation inequality. 
The dissipation inequality is studied in Sect.~\ref{SectValueFunDERIV}.  The final results obtained by this study are (see Sect.~\ref{OPtiCONTeRICCATI}):
\begin{enumerate}
\item the optimal control is a feedback control (see~(\ref{feedbackFORMcontrol}));
\item the feedback form of the optimal control is expressed via a \emph{Riccati operator,}   an operator  which is differentiable in a suitable 
weak sense and which solves   a 
differential equation of Riccati type, Eq.~(\ref{eqRiccati}).
 
\end{enumerate}

The term ``feedback'' has to be  explained: the control is a feedback of the \emph{state} of the system. The definition of the state and a differential equation for the evolution of the state are in Sect.s~\ref{SectSTATEdefi} and~\ref{SubsEQUAstate}. We believe that the results in these sections are crucial for understanding the feedback structure of the optimal control.

\section{\ZLA{SECT:soluzione}The solutions of~(\ref{eq:IntroNOTAZ}) and its well posedness}
We perform a chain of formal computations and we find a formula which can be interpreted as a definition of the solutions of Eq.~(\ref{eq:IntroNOTAZ}). 
Then, in Sect.~\ref{SectCondiINIcondiEsoluzL2}, 
we study the space  in which the problem is well posed.
We follow an idea   introduced in~\cite{DelloroPATA17MGTeVISCO} to study the stability of the MGT equation and then used in~\cite{BucciPandolfiArxive} to study the regularity properties of its solution.

We rewrite the system as
\[
(v'-\Delta v)'=-(v'-\Delta v)+v'
\]
so that
 \begin{multline}\ZLA{eq:laFORMintegDIv}
v'(t)-\Delta v(t)=e^{-t}(v_1-\Delta v_0)+\intt e^{-(t-s)}v'(s)\ZD s\\
= e^{-t}(v_1-v_0-\Delta v_0)+v(t)-\intt e^{-(t-s)} v(s)\ZD s\,.
 \end{multline}
Provided we specify the space in which $(v_1-v_0-\Delta v_0)$ has to be chosen, Eq.~(\ref{eq:laFORMintegDIv}) is a special instance of the heat equation with memory studied in~\cite[Ch.~3]{PandolfiLIBRO21}. As in this reference, we take into account the boundary condition as follows. We let $A$ be the Dirichlet laplacian in $L^2(\ZOMq)$:
\begin{equation}\ZLA{eq:DIRIcLaplac}
\Dom\, A=H^2(\ZOMq)\cap H^1_0(\ZOMq)\,,\qquad A\phi=\Delta\phi
\end{equation}
and we use the operator $D$ which associates to $u\in U=L^2(\partial\ZOMq)$ its harmonic extension in $\ZOMq$ so that $\Delta Du=0$ while $ADu\in (\Dom\,A)'$
(as usual,     $A$ denotes both the operator in~(\ref{eq:DIRIcLaplac}) and its extension $(A^*)'$ to $(\Dom\, A)'$. 
See~\cite[Ch.~2]{PandolfiLIBRO21} and~\cite[Appendix]{PPZ} for a detailed discussion).

 Problem~(\ref{eq:IntroNOTAZ}) can be recast  as the following integrodifferential equation
\begin{equation}
\ZLA{eq:laFORMintegDIvCONuINTEGROdiff}
v'(t)=Av(t)+e^{-t}\hat y+v(t)-\intt e^{-(t-s)} v(s)\ZD s-AD u(t)\,,\quad \hat y= v_1-v_0-\Delta v_0 \,.
\end{equation}
This equation makes sense in $(\Dom\,A)'$ (and also in a smaller space, as we discuss below).

We use the fact that $A$ generates a holomorphic exponentially stable semigroup and we rewrite~(\ref{eq:laFORMintegDIvCONuINTEGROdiff}) as the following integral equation of Volterra type in $(\Dom\,A)'$ :
\begin{multline}
\ZLA{eq:laFORMintegDIvCONuVOLTE}
v(t)= e^{At} v_0+\intt e^{A(t-s)} e^{-s}(v_1-v_0-\Delta v_0)\ZD s -\intt e^{A(t-s)} AD u(s)\ZD s 
\\
+\intt N(t-s) v(s)\ZD s
\end{multline}
where
\begin{equation}\ZLA{eq:defiNUcleoN}
N(t)= e^{At}-\intt e^{-(t-s) } e^{As} \ZD s\,.
\end{equation}

\begin{Remark}\ZLA{REMA:VolteINvarIspazi}{\rm
Note that:
\begin{itemize}
\item the operators $A$, $e^{At}$ and $N(t)$ are selfadjoint
and  $N(t)$ commutes with $A$ and with its fractional powers. So,~(\ref{eq:laFORMintegDIvCONuVOLTE}) is a Volterra integrodifferential equation  in any space $ \Dom(-A)^\zaa$ and in the dual   $(\Dom(-A)^\zaa)'$, provided that the affine term belongs to this space.
\item if the affine term of~(\ref{eq:laFORMintegDIvCONuVOLTE}) is of class $C^1$ or $H^1$ then the solution $v$ is of class $C^1$ or respectively $H^1$ too. 
\end{itemize}
A consequence is  that when studying (first order) differentability of the solution or its regularity in terms of the fractional powers of $-A$, we can confine ourselves to study the corresponding property of the affine term.\zdia
}
\end{Remark}
If $v$ is a $C^1$ solution of~(\ref{eq:laFORMintegDIvCONuINTEGROdiff}) then it solves~(\ref{eq:laFORMintegDIvCONuVOLTE}). Conversely, 
a solution of~(\ref{eq:laFORMintegDIvCONuVOLTE}) solves~(\ref{eq:laFORMintegDIvCONuINTEGROdiff})
in the weak  sense specified in
  Theorem~\ref{teoLasolUdebDeVolte} at the end of   section~\ref{SectCondiINIcondiEsoluzL2}  which follows from  the properties of the solutions of the Volterra integral equation~(\ref{eq:laFORMintegDIvCONuVOLTE}) which we are going to derive.
Here we note  that two of the tree addenda of the affine term of~(\ref{eq:laFORMintegDIvCONuVOLTE}) are standard terms encountered in the study of the quadratic regulator problem for the heat equation with Dirichlet  boundary control:
  \begin{enumerate}
  \item the transformation $v_0\mapsto e^{At } v_0$  is continuous from $v_0\in H=L^2(\ZOMq)$ to $C([0,T];H)$. Even more: we recall the following fact which holds for exponentially stable holomorfic semigroups. Let $\zaa\in[0,1]$. For every $t>0$, $e^{At}(-A)^\zaa$ admits a continuous extension to $H$ and
  \begin{equation}
  \ZLA{Eq:IneqPERsemiOLOM}
  \| e^{At}(-A)^\zaa\|_{\mathcal{L}(H)}\leq \dfrac{M}{t^\zaa}\,.
  \end{equation}
  So, if $\zaa\in [0,1]$ and if $v_0=(-A)^{\zaa}\tilde v_0$, i.e. $v_0\in(\Dom\,(-A)^\zaa)'$, then
  $v_0\mapsto e^{A\cdot} v_0\in \mathcal{L}\left ((\Dom\,(-A)^\zaa)',L^p(0,T;H)\right )$, any $p\in [1,1/\zaa)$.
  \item the transformation
  \begin{equation}\ZLA{Eq:MappaCONTROtoStaPARZ}
u\mapsto   \int_0^{\cdot} e^{A(\cdot-s)} AD u(s)\ZD s 
  \end{equation}
  is linear and continuous from $L^2(0,T;U)$ to $L^2(0,T;H)$ (and to a smaller space, see statement~\ref{I2Teo:DiseFonda} of Theorem~\ref{Teo:DiseFonda}).
  \end{enumerate}

The novelty of   the study of the  wave equation with strong damping is the memory and the term
\begin{equation}\ZLA{eq:NoveltyDAMPING}
y(t)=
\intt e^{A(t-s)} e^{-s} \hat y\ZD s\,,\qquad \hat y=v_1-v_0-\Delta v_0 
\end{equation}
which is discussed in the next section. This discussion leads to the definition of the   solutions of Eq.~(\ref{eq:IntroNOTAZ}) in the weak  sense specified
 in Theorem~\ref{teoLasolUdebDeVolte}.

\subsection{\ZLA{SectCondiINIcondiEsoluzL2}The definition of the solutions  of~(\ref{eq:IntroNOTAZ})
}

The following fact are well known (see~\cite{Fujiwara,LasiePoteFRAZ,LionsMagenesBOOKvol1}):
\begin{itemize}
\item we have:
\begin{equation}\ZLA{Eq:IncluDomiPoteFraz}
\left\{
\begin{array}{l}
\Dom\,(-A)^{\zaa} = \left\{\begin{array}{l}
 \mbox{$H^{2\zaa}(\ZOMq)$   if  $ \zaa\in [0, 1/4)$}\ \mbox{(with $H^0(\ZOMq)=L^2(\ZOMq)=H$)}\\[2mm]
 \mbox{$H^{2\zaa}_0(\ZOMq)$    if   $\zaa\in (1/4,3/4)$}
 %\\[2mm]
% \mbox{$H^{2\zaa}(\ZOMq)$   if  $ \zaa\in [0, 1/4)$}
\end{array}\right. 
\\[4mm]
\mbox{it follows:}\\[1mm]
\zaa\in (1/4,3/4)\ \implies\ H^{-2\zaa}(\ZOMq)=\left (H^{2\zaa}_0(\ZOMq)\right )'=\left (\Dom\,(-A)^{\zaa}\right )'\\
H^{1/2+\ZEP}(\ZOMq)\subseteq H^{1/2}(\ZOMq)\subseteq H^{1/2-\ZEP}(\ZOMq)= \Dom(-A)^{1/4-\ZEP/2}\quad\mbox{($\ZEP>0$)}
\end{array}     
\right.
\end{equation}
%(we intend   $H^0(\ZOMq)=L^2(\ZOMq)=H$). 

\item   if $u\in U=  L^2(\partial\ZOMq)$ we have
\[
Du\in H^{1/2}(\ZOMq)\subseteq \Dom\,(-A)^{ \frac{1}{4}-\ZEP} \quad\mbox{(any $\ZEP>0$, strict inequality)}
\] 
  and the map $u\mapsto Du$ belongs to $\mathcal{L}\left (U,
\Dom\,(-A)^{ \frac{1}{4}-\ZEP}\right )$.

\item
we denote $\zg_0$ the trace operator $\zg_0 v= v_{|_{\partial\ZOMq}}$ when $v$ is defined on $\ZOMq$. We have:
 $\zg_0\in\mathcal{L}(H^{1/2+\ZEP}(\ZOMq) , H^{\ZEP}(\partial\ZOMq)) $ for any $\ZEP>0$ (note the strict inequality) while $D\in\mathcal{L}(H^s(\ZOMq),H^{s+\frac{1}{2}}(\ZOMq))$ for every $s\geq 0$. So, for $\ZEP>0$ (strict inequality)
\begin{equation}\ZLA{eqPROPEtraccia}
 v\in H^{1/2+\ZEP}(\ZOMq)\ \implies\  \mbox{$v-D\zg_0  v\in H^{1/2 +\ZEP}_0(\ZOMq)=\Dom(-A)^{1/4+\ZEP/2}$}\,.\\
%& \mbox{provided that $  v\in H^{1/2+\ZEP}(\ZOMq)$ (any $\ZEP>0$)}\,.
\end{equation}

\end{itemize}

Inequality~(\ref{Eq:IneqPERsemiOLOM}) and the previous properties imply:

\begin{Theorem}\ZLA{Teo:DiseFonda}
The following properties hold:
\begin{enumerate}
\item\ZLA{I1Teo:DiseFonda} there exists $\ZSI\in (0,1)$ such that ${\rm im}\,D\subseteq \Dom\,(-A)^{1-\ZSI} $ so that
\begin{align*}
& AD\in\left (\Dom\,(-A)^{\ZSI}\right )'\\
&\| e^{At}AD\|_H=
\left \|e^{At}(-A)^{\ZSI}\left ((-A)^{1-\ZSI}D\right )\right \|_H
 \leq M\dfrac{1}{t^\ZSI}\quad t\in [0,T]\,,\quad \ZSI\in (0,1)
\end{align*}
 (more precisely, $\ZSI\in(3/4,1)$).
\item\ZLA{I2Teo:DiseFonda} there exists $p_0\in (1,2)$ (more precisely $p_0\in (1,1/\ZSI)\subseteq (1,4/3)$) such that  $e^{At}AD\in L^{p_0}(0,T,\mathcal{L}(U,H)$.
So, Young's inequalities (recalled in the Appendix~\ref{AppePROpreHzt}) imply
\begin{enumerate}
\item\ZLA{I21Teo:DiseFonda}
The map~(\ref{Eq:MappaCONTROtoStaPARZ}) is linear and continuous from $L^2(0,T;U)$ to $L^r(0,T;H)$ where
\begin{equation}\ZLA{defiEXPOr}
\dfrac{1}{r}
=\dfrac{1}{p_0}-\dfrac{1}{2}\quad {\rm i.e.}\quad r=\dfrac{2p_0}{2-p_0}.
\end{equation}
Note that $r>2$ since $p_0>1$.
\item\ZLA{I22Teo:DiseFonda}
Let $\ZEP>0$. The map~(\ref{Eq:MappaCONTROtoStaPARZ}) belongs to  
 \begin{multline*}
\mathcal{L}\left (
L^{\ZIN}([0,T];U),C([0,T];(\Dom(-A)^{1-\ZSI-\ZEP})\right )
 \\
\subseteq \mathcal{L}\left (
L^{\ZIN}([0,T];U),C([0,T];H)\right ) \,.
\end{multline*}
\end{enumerate}
\end{enumerate}
\end{Theorem}

  To help the reader, we collect   here the sharp values   of $\ZSI$, $p_0$ and $r$ encountered in the theorem.
We let $\ZEP\sim  0$ be a ``small'' \emph{strictly positive} number, not the same at every occurrence. Then we have
\begin{equation}\ZLA{eqDefiEXPLIdi]ZSIeP}
\left\{\begin{array}{ll}
\ZSI\in (3/4,1)\,:&\ZSI=3/4+\ZEP\,, 
\\
p_0\in (1,4/3)\,:&p_0=1+\ZEP\,,
\\
r=2p_0/(2-p_0)\,:
&r=2+\ZEP 
\,.
\end{array}\right.
%\mbox{$\ZSI\in (3/4,1)$ and $p_0\in (1,1/\ZSI)$, in particular  $p_0\in (1,4/3)$}
\end{equation}

  Now we examine $y(t)$, i.e. the integral in~(\ref{eq:NoveltyDAMPING}).

Let $v_0\in  H^{1/2+\ZEP}(\ZOMq)$ ($\ZEP>0$ small). From~(\ref{eqPROPEtraccia}) we have
 \begin{multline*}
\Delta v_0=\Delta(v_0-D\zg_0 v_0)=A(v_0-D\zg_0 v_0)\\
=-(-A)^{3/4-\ZEP/2}\left [ (-A)^{1/4+\ZEP/2}(v_0-D\zg_0 v_0)   \right ] 
\in\left (\Dom\,(-A)^{3/4-\ZEP/2}\right )' 
\,.
 \end{multline*}
 From the inequality~(\ref{Eq:IneqPERsemiOLOM}),
  \begin{equation}\ZLA{DefiZSI1p1}
  \left\{\begin{array}{l}
 \|e^{At}\Delta v_0\|_{H} 
\leq \frac{M}{t^{\ZSI_1}}\| v_0-D\zg_0v_0\|_{\Dom\,(-A)^{1/4+\ZEP/2}} 
  \leq \frac{M}{t^{\ZSI_1}} \|v_0\|_{H^{1/2+\ZEP}(\ZOMq)} \,, \\[2mm]
   e^{At}\Delta v_0\in L^{p_1}(0,T;H)\,, \\
\ZSI_1={3/4-\ZEP/2}\in (0,3/4)\,,\ p_1=\dfrac{4}{3-2\ZEP}>\dfrac{4}{3}\,.
\end{array}\right.
  \end{equation}

  This observation suggests to choose the following spaces (we recall that $\ZEP$ denotes an arbitrary ``small'' positive number)
  \begin{multline}\ZLA{eq:DefiHATy}
  \left\{\begin{array}{l}
  v_0\in H^{1/2+\ZEP}(\ZOMq)\subseteq \Dom(-A)^{1/4-\ZEP/2}\ \\
v_1\in(\Dom\,(-A)^{3/4-\ZEP/2} )'
\end{array}
\right.
%\ \mbox{so that}
\\
  \mbox{so that}\quad  % \\
 \left\{\begin{array}{l}
\hat y=v_1-v_0-\Delta v_0\in (\Dom\,(-A)^{3/4-\ZEP/2} )'\\
\| e^{At}\hat y\|_{H}\leq \dfrac{M}{t^{3/4-\ZEP/2}}\|\hat y\| _{(\Dom\,(-A)^{3/4+\ZEP} )'}\\
e^{At}\hat y\in L^{p_1}(0,T;H)\subseteq L^{p_0}(0,T;H)\,,\ p_1>4/3>1
\,.
\end{array}\right.
  \end{multline}
   
   We recall from~(\ref{eq:NoveltyDAMPING}) and the previous considerations:
 \begin{equation}
  \ZLA{Eq:DefiRefeSigy}
  \left\{\begin{array}{l}
%  \begin{array}{ll}
%  \displaystyle
%  y(t) &=   \intt e^{-(t-s)}e^{A s }  \hat y\ZD s=
% \\[1mm]
%  \displaystyle &=    e^{At}A^{-1}\hat y -e^{-t}A^{-1}\hat y -\intt e^{-(t-s)}e^{As}A^{-1}\hat y\ZD s   
% \end{array}\\
 %
%\\ 
  \displaystyle\left\{\begin{array}{l}
  \displaystyle v_0\in H^{1/2+\ZEP}(\ZOMq)\\
  \displaystyle v_1\in \left (\Dom\,(-A)^{3/4-\ZEP/2}\right )'
\end{array}\right.\ \implies\ \left\{\begin{array}{l}
  \displaystyle\hat y\in (\Dom\,(-A)^{3/4-\ZEP/2} )'\\
  \displaystyle A^{-1}\hat y\in \Dom\,(-A)^{1/4-\ZEP/2} 
\end{array}\right.
\\[2mm]
\begin{array}{ll}
  \displaystyle
  y(t) &=   \intt e^{-(t-s)}e^{A s }  \hat y\ZD s=
 \\[1mm]
  \displaystyle &=    e^{At}A^{-1}\hat y -e^{-t}A^{-1}\hat y -\intt e^{-(t-s)}e^{As}A^{-1}\hat y\ZD s   \,.
 \end{array} 
\end{array}\right.
  \end{equation}
 
  In order to understand the following theorem we note 
  \[
(\Dom\,(-A)^{3/4-\ZEP/2} )'\subseteq (\Dom\,A)'  
  \]
  and that the integration by parts in~(\ref{Eq:DefiRefeSigy}) shows that $y\ZCD$ takes values in $H$ even if $\hat y\in (\Dom\,A)'$.
  
  \begin{Theorem}\ZLA{teo:regolarDELriferimento} 
 Let $y\ZCD$ be the function in~(\ref{eq:NoveltyDAMPING}).
The following properties hold:
  \begin{enumerate}
    \item \ZLA{I0teo:regolarDELriferimento} if $\hat y\in(\Dom\, A)'$ then $y\in C([0,T];H)$ and the transformation $\hat y\mapsto y\ZCD$ is continuous between the indicated spaces.
\item\ZLA{I1teo:regolarDELriferimento} 
  the map $\hat y\mapsto y\ZCD$ belongs to
 $
    \mathcal{L}\left ((\Dom\,(-A)^{3/4-\ZEP} )',C([0,T];\Dom\,(-A)^{1/4-\ZEP/2})\right)$ and so also to
    $  
   \mathcal{L}\left ((\Dom\,(-A)^{3/4-\ZEP} )',C([0,T];H)\right ) 
  $;
     \item\ZLA{I2teo:regolarDELriferimento}  Let $\hat y\in (\Dom\,(-A)^{3/4-\ZEP} )'$. We have
    \[
y'\ZCD=    e^{A\cdot}\hat y-y\ZCD\in L^{p_1}(0,T;H)\subseteq L^{p_0}(0,T;H)
    \]
($p_1 $ is in~(\ref{DefiZSI1p1}), $p_1>4/3$)
  and    $y\in C^1\left( [0,T];(\Dom\,(-A)^{3/4-\ZEP/2} )'\right ) $.
 \item \ZLA{I3teo:regolarDELriferimento} If   $\hat y\in H$ then  $y\ZCD$ takes values in $\Dom\,A$     and the map $\hat y\mapsto y\ZCD$ belongs to $\mathcal{L}( H , C([0,T];\Dom\,A))$. 
 \end{enumerate}
  \end{Theorem}

  Theorems~\ref{Teo:DiseFonda} and~\ref{teo:regolarDELriferimento} imply that the Volterra integral equation~(\ref{eq:laFORMintegDIvCONuVOLTE}) admits a unique solution in $L^2(0,T;H)$ which depends continuously on $v_0\in H^{1/2+\ZEP}(\ZOMq) $,    $ v_1 \in (\Dom\,(-A)^{3/4+\ZEP} )' $ and $u\in L^2(0,T;L^2(\ZOMq))$.   
 
  \begin{Definition}\ZLA{Defi:defiSOLU}
  Let  $u\in L^2(0,T;U)$, $v_0\in H^{1/2+\ZEP}(\ZOMq)$, $v_1\in (\Dom\,(-A)^{3/4-\ZEP} )'$.   
  The unique (weak) solution of~(\ref{eq:IntroNOTAZ}) is  (by definition!)
   the unique solution of~(\ref{eq:laFORMintegDIvCONuVOLTE}) in $L^2(0,T;H)$.  
  \end{Definition}
  
  \begin{Remark} \ZLA{RemaSUiniv0inH}{\rm
We note:
 \begin{enumerate}
 \item\ZLA{I0RemaSUiniv0inH} $\hat y$ depends on $v_0$ and $v_1$ but not on $u(0)$ which is not even defined since $u\in L^2(0,T;U)$.
 \item\ZLA{I1RemaSUiniv0inH}
   the fact that both $v_0$ and $v_1$ are arbitrary in the respective spaces imply that $v_0$ and $\hat y$ are independent.\zdia
 
\end{enumerate}
  }\end{Remark}
 
 The following observation has a particular interest and we state it in a separate remark:
 \begin{Remark}[Important observation]\ZLA{I2RemaSUiniv0inH}{\rm
   Once that the contribution of $\Delta v_0$ has been absorbed   in $\hat y$,   the Volterra integral equation~(\ref{eq:laFORMintegDIvCONuVOLTE})
can be considered as depending $u\in L^2(0,T;U)$ and  on the independent parameters $v_0\in H^{1/2+\ZEP}(\ZOMq)$ and $\hat y \in (\Dom\, (-A)^{3/4-\ZEP/2})' $. The important observation is that we can relax these assumptions on $v_0$ and $\hat y$ and we can study Eq.~(\ref{eq:laFORMintegDIvCONuVOLTE}) under the weaker conditions
\begin{equation}\ZLA{Eq:I2RemaSUiniv0inH}
v_0\in H\,,\qquad \hat y\in (\Dom\,A)'\,.\zdiaform
\end{equation}

 }\end{Remark}

\emph{From now on, we relay on the important Remark~\ref{I2RemaSUiniv0inH}
and in particular we use the conditions~(\ref{Eq:I2RemaSUiniv0inH}).}

 Theorem~\ref{teo:regolarDELriferimento}, the    Remarks~\ref{REMA:VolteINvarIspazi} and~\ref{I2RemaSUiniv0inH} and the properties of the solutions of the heath equation with Dirichlet boundary control imply the following result concerning the regularity of the solutions of the Volterra integral equation~(\ref{eq:laFORMintegDIvCONuVOLTE}):
 \begin{Theorem}\ZLA{teoREgoDiv}
We use the  specific values of $\ZSI$, $p_0$ and $r$
 in~(\ref{eqDefiEXPLIdi]ZSIeP}) and    we denote $\ZEP$ a \emph{strictly positive} number, not the same at every occurrence.
 
 The solution $v$ of the Volterra integral equation~(\ref{eq:laFORMintegDIvCONuVOLTE}) has the following properties:
 \begin{description}
 \item[\bf $\bullet$ In terms of $v_0$, $\hat y$ and $u$.]
 
 {~}
 \begin{enumerate}
 \item\ZLA{I2teoREgoDivbis} if $v_0\in (\Dom\, (-A)^{1/r-\ZEP})'$ , $\hat y\in(\Dom\,A)'$, $u\in L^2(0,T;U)$
 then $v\in L^r(0,T;H)\subseteq L^2(0,T;H)$.

\item\ZLA{I1teoREgoDiv} we have:

 $v\in C([0,T];H)$ provided that 
 \[\left\{\begin{array}{l}
v_0\in H\,,\quad \hat y\in (\Dom\, A )'\,,\\
 \mbox{$ u\in L^q([0,T];U)$    any $q>p_0/(p_0-1)$}\,.
 \end{array}\right.
 \]
In particular, $v$ is continuous when
$v_0\in H$, $ \hat y\in (\Dom\, A )'$ and
   $u\in L^\ZIN(0,T;U)$.

  \item\ZLA{I0teoREgoDiv}   If $v_0\in  \Dom\,(-A)^{1-\ZSI-\ZEP}$, $\hat y\in(\Dom\,(-A)^{\ZSI+\ZEP})'$ and $u\in L^{\ZIN}( 0,T ;U)$ then $v\in C([ 0,T];\Dom\,(-A)^{1-\ZSI-\ZEP})$, $v'\in L^{\ZIN}( 0,T;(\Dom\,(-A)^{\ZSI+\ZEP})')$. If $u$ is continuous then $v$ and $v'$ (with values in the indicated spaces) are continuous too.

 \item\ZLA{I5teoREgoDiv} We have $v\in C^1([0,T];H)$ provided that $v_0\in H$,
 $\hat y\in (\Dom\,A)'$,  $u\in C^1([0,T];U)$ 
 and   the following consistency condition hold:
 \[
A(v_0- Du(0))+\hat y \in H\,. 
 \]

 \end{enumerate}
  \item[\bf $\bullet $ In terms of $v_0$, $v_1$     and $u$.] 
  {~}

\begin{enumerate}[resume]
 
%  \begin{enumerate}\setcounter{enumi}{5}
  \item 
Let $v_0\in H^{1/2+\ZEP}(\ZOMq)\subseteq \Dom(-A)^{1/4-\ZEP/2}$, $v_1\in (\Dom\,(-A)^{3/4-\ZEP/2})'$ and $u\in L^2(0,T;H)$. Then    $v\ZCD\in L^r(0,T;H) $ and the map $(v_0,v_1,u)\mapsto v\ZCD$ is linear ad continuous from  
the specified spaces  to $L^r(0,T;H)$.

  \item If  $u\in C^1([0,T];U)$, $v_0\in H^{1/2+\ZEP}(\ZOMq)$, $v_1\in (\Dom\,(-A)^{3/4-\ZEP/2})'$ and if
$
A(v_0-Du(0))+v_1-\Delta v_0\in H 
$
then $v\in C^1([0,T];H)$.

  \end{enumerate}
 \end{description}
 \end{Theorem}
 
 Further properties of $v$ in terms of  of $v_0$, $v_1$     and $u$ are easily derived from those in terms of $v_0$, $\hat y$ and $u$ but we don't need such information.
 
The properties in Theorem~\ref{teoREgoDiv} (in particular the property   in the statement~\ref{I0teoREgoDiv}) and   $\ZSI>3/4$ imply that $v(t)$ solves~(\ref{eq:laFORMintegDIvCONuINTEGROdiff}) in the weak  form~(\ref{eq:laFORMintegDIvCONuINTEGROdiffDEB}) below:

\begin{Theorem}\ZLA{teoLasolUdebDeVolte}

Let
$v_0\in H$, $ \hat y\in (\Dom\,A)'$, $u\in L^2(0,T;U)$. Let $\xi\in\Dom\,(-A)^{1+\ZSI}$. 
The function $\ZL v(t),\xi\ZR$ belongs to $H^1(0,T;H)$ and we have
\begin{multline}\ZLA{eq:laFORMintegDIvCONuINTEGROdiffDEB}
\dfrac{\ZD}{\ZD t} \ZL v(t),\xi\ZR=\ZL v(t), (A+I)\xi\ZR +\ZL e^{-t}\hat y,\xi\ZR\\
 -\ZL Du(t),A\xi\ZR-\intt\ZL e^{-(t-s)} v(s),\xi\ZR\ZD s\,. 
\end{multline}
 
\end{Theorem}

  \section{\ZLA{SECTandFAMILYquad}A family of quadratic regulator problems and the properties of the optimal control}

Let $\zt\in (0,T)$.  The restriction to $[\zt,T)$ of the weak  solution  of~(\ref{eq:laFORMintegDIvCONuINTEGROdiff})  is the weak  solution of
\begin{multline}\ZLA{eq:laEQsu(ztT)ANTE}
v'(t)=Av(t) +v(t)-\int_\zt^t e^{-(t-s)} v(s)\ZD s-\int_0^\zt e^{-(t-(\zt-s))}v(\zt-s)\ZD s\\
-AD u(t) +e^{-(t-\zt)}\hat y_{\zt} \,,
 \qquad  \hat y_{\zt}=   e^{-\zt}\hat y \quad v(\zt^+)=v(\zt^-)
\,.
\end{multline} 
The condition $v(\zt^+)=v(\zt^-)$ makes sense in a space in which the solution is continuous. So, in general the equality makes sense in the space    $(\Dom\,A)'$.

  Eq.~(\ref{eq:laEQsu(ztT)ANTE}) makes sense and defines a well posed dynamical system on $(\zt,T)$ even if the triple $(v(\zt^-),v_{|_{(0,\zt)}},e^{-\zt} \hat y)$  is replaced by an arbitrary element  $(\hat v_\zt, \xi ,\hat   y_\zt)$ in the space $H\times L^2(0,\zt;H)\times (\Dom\,A)'$. So, we consider the following family of  equations parametrized by $\zt\in[0,T)$:
\begin{subequations}
\begin{multline}
\ZLA{eq:laEQsu(ztT)}
v'(t)=Av(t) +v(t)-\int_\zt^t e^{-(t-s)} v(s)\ZD s-AD u(t)
\\
-e^{(t-\zt)}\int_0^\zt e^{-s}  \xi_\zt(\zt- s)\ZD s +e^{-(t-\zt)}\hat  y_{\zt}\,,\\ 
v(\zt )=\hat v_\zt\,,\qquad  \xi_\zt\ZCD  \in L^2(0,\zt;H)\,,\qquad\hat  y_{\zt}   \in (\Dom\,A)' \,.
\end{multline} 
\emph{We repeate that Eq.~(\ref{eq:laEQsu(ztT)}) can be studied when $\hat v_\zt$ is \emph{any}  elements of $H$,  $\xi_\zt\ZCD $ is \emph{any} element of $ L^2(0,\zt;H)$ and $\hat  y_{\zt}$ is \emph{any} element of $(\Dom\, A)'$, not related to a solution of the system on $[0,\zt]$.}

  We introduce   the notations 
\begin{equation}\ZLA{EQpreDEFIstate}
\begin{array}{l}
\displaystyle M^2_\zt=
 H\times L^2(0,\zt;H)   \\[1mm]
\displaystyle 
 \Xi_\zt=
(\hat v_\zt, \xi_\zt(\zt-\cdot))\in M^2_\zt 
\\[1mm]
\displaystyle 
    \Sta_\zt=
(\hat v_\zt, \xi_\zt(\zt-\cdot) ,\hat y_{\zt})=(\Xi_\zt,\hat y_ {\zt})  \in M^2_\zt\times (\Dom\,A)'
\end{array}
\end{equation}
 (if $\zt=0$ then   $M^2_0=H\times\{0\}$ is isomorphic to $H$. In this case we can put $M^2_0=H$ and $\Xi_0=\hat v_0$). 
 
 The considerations in Sect.~\ref{SECT:soluzione} and in its subsection~\ref{SectCondiINIcondiEsoluzL2} concerning the problem~(\ref{eq:laFORMintegDIvCONuINTEGROdiff}),
  in particular Theorem~\ref{teoREgoDiv} and Theorem~\ref{teoLasolUdebDeVolte},  are easily extended to the problem~(\ref{eq:laEQsu(ztT)}). So we can state in particular  that
 Eq.~(\ref{eq:laEQsu(ztT)}) is uniquely solvable and also that the solution $v\ZCD\in L^2([\zt,T];H)$ depends linearly and continuously on    $(\Sta_\zt, u\ZCD) 
 \in (M^2_\zt\times( \Dom\, A) ')\times  L^2(\zt,T;H)$. 
 
We associate to~(\ref{eq:laEQsu(ztT)}) the optimization problem
\begin{multline}
\ZLA{eq:costoWAVEdampSU(ztT)}
\min_{u\in L^2(\zt,T;U)} J_\zt(\Sta_\zt;u)=
\min_{u\in L^2(\zt,T;U)} J_\zt(( \Xi_\zt ,\hat  y_{\zt});u) \\
J_\zt(\Sta_\zt;u)= \int_\zt^T\left [\|v(t;\zt,\Sta_\zt ,u )\|^2_H+\|u(t)\|^2_u\right ]\ZD t 
\end{multline}
where  $v(t;\zt,\Sta_\zt ,u ) $ is the solution of~(\ref{eq:laEQsu(ztT)}).
\end{subequations}

Continuity of the linear
 transformation 
\[%\begin{multline*}
(\Sta_\zt,u) \mapsto v(t;\zt,\Sta_\zt,u )\,:\quad  M^2_\zt\times (\Dom\, A )'\times L^2(\zt,T;U)\mapsto L^2(\zt,T;H)
\]%\end{multline*}
    implies that    the cost, as a function of $u\in L^2(\zt,T;U)$, is coercive 
    for every fixed triple $ \Sta_\zt $ and we can state:

\begin{Theorem}\ZLA{Teo:WELLposedINM2zt} 
Let $\Sta_\zt\in M^2_\zt\times(\Dom\,A)'$ be fixed.
Problem~(\ref{eq:laEQsu(ztT)})-(\ref{eq:costoWAVEdampSU(ztT)})  admits a unique optimal control  
\[
u^+_\zt(\Sta_\zt) =u^+_\zt(\cdot;,\Sta_\zt)=u^+_\zt(\cdot;( \Xi_\zt ,\hat y_{\zt}) )\,.
\]

\end{Theorem}

When $u=u^+_\zt(\Sta_\zt)$ in~(\ref{eq:laEQsu(ztT)}), the solution is denoted
$v^+_\zt(\cdot; \Sta_\zt )$:
\[
v^+_\zt( \Sta_\zt )= 
v^+_\zt(\cdot; \Sta_\zt )=v(\cdot ;\zt,\Sta_\zt,u^+_\zt(\Sta_\zt))\,.
\]

 Finally we observe that the differential equation~(\ref{eq:laEQsu(ztT)})   in $(\Dom\,(A))'$ is the weak form of the following Volterra integral equation in $H$:
 \begin{multline}
 \ZLA{eq:COMEvolteEQsu(ztT)}
 v(t)=e^{A(t-\zt)}\hat v_\zt  
+\int_\zt^t e^{A(t-s)}e^{-(s-\zt)}\left [ \hat  y_\zt-  \int_0^\zt e^{- \nu } \xi_\zt (\zt-\nu)\ZD\nu\right ] \ZD s\\
 -\int_\zt^t e^{A(t-s)} AD u(s)\ZD s 
 +\int_\zt^t N(t-s)v(s)\ZD s 
 %\,,\quad y(t)=e^{-(t-\zt)}\hat  y_\zt
 \end{multline}
 where $N(t)$ is the function in~(\ref{eq:defiNUcleoN}).

We see from~(\ref{eq:laEQsu(ztT)})  or~(\ref{eq:COMEvolteEQsu(ztT)})   that when studying the properties of the solutions
 the contribution of $ \xi_\zt\ZCD$    can be treated similarly to that of $\hat y_\zt$ and we can even assume that $  \xi_\zt\ZCD$ is a $(\Dom\,A)'$-valued function.  Theorem~\ref{teoREgoDiv} gives in particular:
 
\begin{Theorem}
\ZLA{Teo:preREGOLARITAstato}The solution $v$ of~(\ref{eq:COMEvolteEQsu(ztT)}) has the following properties:
\begin{enumerate}
\item\ZLA{I1Teo:preREGOLARITAstato} 
Let $u=0$. We have:
\begin{enumerate}
\item\ZLA{I1ATeo:preREGOLARITAstato} if $\hat v_\zt\in H$, $\hat y_\zt\in (\Dom\,A)'$, $ \xi_\zt\in L^2(0,\zt;(\Dom\,A)')$ then $v\in C([\zt,T];H)$.
\item\ZLA{I1BTeo:preREGOLARITAstato} if $\hat v_\zt\in\Dom\, A$, $\hat y_\zt\in H$, $ \xi_\zt\in L^2(0,\zt;H)$ then $v\in C([\zt,T];\Dom\,A)\cap C^1([\zt,T];H)$.
\item\ZLA{I1CTeo:preREGOLARITAstato} let $\hat v_\zt\in\Dom\, A$, $\hat y_\zt\in H$, $ \xi_\zt\in H^1(0,\zt;H)$ and let the compatibility condition $ \xi_\zt(0)=\hat v_\zt$ hold. The following function  is of class $H^1(0,T;H)$:
 
\begin{equation}\ZLA{Eq:primaDefiEXTE}
\left\{\begin{array}{lll}
  \xi_\zt(t-\zt)&{\rm if}&0\leq t\leq \zt\\
v(t)& {\rm if}&\zt\leq t\leq T\,. 
\end{array}\right.
\end{equation}

\end{enumerate}

\item\ZLA{I2Teo:preREGOLARITAstato}
Let $\hat v_\zt\in H$, $ \xi_\zt\in L^2(0,\zt;(\Dom\,A)')$, $\hat y_\zt\in  (\Dom\,A)' $.
\begin{enumerate}
 
\item\ZLA{I2ATeo:preREGOLARITAstato} if $u\in L^2(\zt,T;U)$ then $v\in L^r(\zt,T;H)\subseteq L^2(0,T;H)$ ($r$ is the exponent in~(\ref{eqDefiEXPLIdi]ZSIeP})). 
\item\ZLA{I2BTeo:preREGOLARITAstato}
if $u\in L^\ZIN (\zt,T;H)$  then $v\in C([\zt,T];H)$. If furthermore $  \xi_\zt\in C([0,\zt];H)$ and   the compatibility condition $  \xi_\zt(0)=\hat v_\zt$ holds then the function~(\ref{Eq:primaDefiEXTE}) is continuous on $[0,T]$.
\end{enumerate}
 
\item\ZLA{I3Teo:preREGOLARITAstato}  Let $\hat v_{\zt}\in H$, $\hat y_{\zt}\in (\Dom\,A)'$, $ \xi_\zt \in L^2(0,\zt;(\Dom\, A)')$ and $u\in C^1([0,\zt];U)$ and let the following consistency condition hold:
\[
A(\hat v_\zt
-D u(\zt))+\left (\hat y_\zt+\int_0^\zt e^{-\nu}  \xi_\zt( \nu)\ZD \nu\right )\in H\,.
\]
Under these conditions, $v\in C^1([\zt,T];H)$.

\end{enumerate}
\end{Theorem} 
 
 \begin{Remark}{\rm
For the sake of completeness, in the previous theorem we considered the general case that    $ \xi_\zt$ takes values in $ (\Dom\,A)'$.  We don't need such generality and in the rest of the paper and we confine ourselves to  the case   $ \xi_\zt\in L^2(0,\zt,H)$ i.e. we assume condition~(\ref{EQpreDEFIstate}).\zdia
 }\end{Remark}
  
\subsection{A variation of constants formula for~(\ref{eq:laEQsu(ztT)})} 
Let $Z\ZCD\in C([0,T]; \mathcal{L}(H)) $ be the unique strongly continuous solution of the following Volterra integral equation:

\begin{multline}\ZLA{Eq:VoltePerZ}
Z(t)v=e^{At}v +\intt N(t-s)Z(s)v\ZD s\\
=
e^{At}v+\intt e^{A(t-s)}Z(s)v\ZD s-\intt e^{-(t-s)} \int_0^s e^{A(s-r)} Z(r)v\ZD r\,\ZD s
\\
=e^{At}v+\intt e^{A(t-s)}Z(s)v\ZD s-\intt e^{A(t-s)}\int_0^s e^{-(s-r)}Z(r)v\ZD r\,\ZD s
\,.
\end{multline}
  
 We state: 
\begin{Theorem}\ZLA{TeoPropriZt} We recall $\ZEP>0$, not the same at every occurrence.
The function $Z(t)$ has the following properties:
\begin{enumerate}
\item\ZLA{I1TeoPropriZt} $Z(t)=Z^*(t)$ for every $t\geq 0$.
\item\ZLA{I2TeoPropriZt} for every $\xi\in\Dom\, A$ we have $AZ(t)\xi=Z(t)A\xi$. In this sense we say that $Z(t)$ \emph{commutes with $A$.}
\item\ZLA{I5TeoPropriZt}
the function
\begin{equation}\ZLA{equDiPropI5TeoPropriZt}
\int_0^{\cdot}Z(\cdot-s)y(s)\ZD s 
\end{equation}
has the following properties:
\begin{enumerate}
\item\ZLA{I5aTeoPropriZt} if $y\ZCD\in C^1([0,T];H)$ then ~(\ref{equDiPropI5TeoPropriZt}) belongs to  $C^1([0,T];H)$.
\item\ZLA{I5bTeoPropriZt}
if $y\ZCD\in C^1([0,T];(\Dom\,A)')$ then~(\ref{equDiPropI5TeoPropriZt}) belongs to $ C([0,T];H)$.
 \end{enumerate}

\item\ZLA{I3TeoPropriZt} 
  $Z\ZCD ADu\in L^{p_0}(0,T;H)$ for every $u\in U$ (the number $p_0$ is the one in~(\ref{eqDefiEXPLIdi]ZSIeP})).
 
\item\ZLA{I4TeoPropriZt} the map 
\begin{equation}\ZLA{eqMAPPAdiTeoTeoPropriZt}
u\ZCD\mapsto (-A)^{\zaa} \int_0^{\cdot} Z(\cdot-s)ADu(s)\ZD s
\end{equation}
has the following properties:
\begin{enumerate}
\item\ZLA{I4ATeoPropriZt} for every $\zaa\in [0,1-\ZSI)$ there exists $q>2$ such that  the map~(\ref{eqMAPPAdiTeoTeoPropriZt}) belongs to $\mathcal{L}(L^2(0,T;U),L^q(0,T;H))$. In particular when $\zaa=0$ it
belongs to $\mathcal{L}(L^2(0,T;U),L^r(0,T;H))$
(the number $r$ is the one in~(\ref{eqDefiEXPLIdi]ZSIeP}));
\item \ZLA{I4BbisTeoPropriZt} if $\zaa+\ZSI<1$ then the map~(\ref{eqMAPPAdiTeoTeoPropriZt}) is linear and continuous from the space 
$L^{\ZIN}([0,T];U)$ to $C([0,T]; \Dom(-A)^{1-(\zaa+\ZSI)-\ZEP}  )\subseteq C([0,T];H)$.

\item\ZLA{I4BterTeoPropriZt} if $\zaa+\ZSI\geq 1$ then the map  the map~(\ref{eqMAPPAdiTeoTeoPropriZt}) is linear and continuous from the space 
$L^{\ZIN}([0,T];U)$ to $C\left([0,T];  (\Dom(-A)^{\zaa+\ZSI-1+\ZEP})'\right )$.
 
\end{enumerate}

\end{enumerate}
\end{Theorem}
\zProof 
We introduce the notations $\chi$ and $E$ to denote the functions $\chi(t)=e^{-t}$,
$E(t)=e^{At}$ so that $N(t)\xi=E\xi-\chi*E\xi$ ($*$   denotes the convolution). Then:
\begin{equation}\ZLA{serieDiZ}
Z(t)\xi=\sum _{k=0}^{+\ZIN}\left [ (N^{(*k)}*E)(t)\right ]\xi  
\end{equation}
where $ ^{(*k)}$ denotes iterated convolution and we intend
$N^{(*1)}=N$,
 $N^{(*0)}*E=E$.

The series is uniformly convergent both respect to $t$ on compact intervals and to $\xi$ in bounded balls.

\emph{Property~\ref{I1TeoPropriZt})} follows from~(\ref{serieDiZ}) since $N(t)$ is selfadjoint and $N(s)$ and   $E(t)$ commute for every $t $ and $s$.

We prove the \emph{statement~\ref{I2TeoPropriZt}).} If $\xi=A^{-1}\eta\in\Dom\,A$ then
\[
Z(t)\xi=Z(t)A^{-1}\eta=\sum _{k=0}^{+\ZIN}A^{-1} N^{(*k)}*E\eta=A^{-1}
\sum _{k=0}^{+\ZIN}  N^{(*k)}*E\eta=A^{-1}Z(t)\eta
\]
(uniform convergence of the series is used here). This proves $Z(t)\xi\in \Dom\,A$ and
\[
AZ(t)\xi=Z(t)\eta=Z(t)A\xi\,.
\]

We prove \emph{statement~\ref{I5TeoPropriZt}.} 

\emph{Statement~\ref{I5aTeoPropriZt}} follows since
\[
\dfrac{\ZD}{\ZD t}\int_0^t Z(t-s)y(s)\ZD s=
\dfrac{\ZD}{\ZD t}\int_0^t Z( s)y(t-s)\ZD s=Z(t)y(0)+\intt Z(s)y'(t-s)\ZD s\,.
\]
We prove the \emph{statement~\ref{I5bTeoPropriZt}.}
 
By assumption, $y(t)=A \tilde y(t)$ and $\tilde y\ZCD\in C^1([0,T];H)$.
First we consider the special case $Z(t)=e^{At}$. In this case the statement follows from an integration by parts:
\begin{equation}\ZLA{eq:I5bTeoPropriZt}
\intt e^{A(t-s)}A \tilde y(s)\ZD s=-\tilde y(t)+e^{At}\tilde y(0)+\intt e^{A(t-s)}\tilde y'(s)\ZD s\,.
\end{equation} 
The general case follows since
\begin{multline*}
\intt Z(t-s)A\tilde y(s)\ZD s=\intt e^{A(t-s)}A\tilde y(s)\ZD s
%+\intt e^{A(t-s)}A\int_0^s Z(s-r)\tilde y(r)\ZD r\,\ZD s
\\
+\intt e^{A(t-s)}A\left [\int_0^s Z(s-r)\tilde y(r)\ZD r-
\int_0^s e^{-(s-\nu)}\int_0^\nu Z(\nu-r) \tilde y(r)\ZD r\,\ZD\nu 
\right ] \ZD s\,. 
\end{multline*}
Statement~\ref{I5bTeoPropriZt} follows from~(\ref{eq:I5bTeoPropriZt}) since the square bracket in the last integral is of class $C^1$ thanks to the statement~\ref{I5aTeoPropriZt}. 

To prove the \emph{statement~\ref{I3TeoPropriZt}} we  recall
the definition of $p_0$: We know that if $\|\xi\|\leq 1$  then

\[
\|e^{At} (-A)^\ZSI \xi\|<\dfrac{M}{t^{\ZSI}} \quad \mbox{and}\quad u\mapsto e^{A\ZCD}ADu\in\mathcal{L}(U,L^{p_0}(0,T;H) )\,.
\]
Then,
  $Z(t)ADu$ solves the following Volterra integral equation in $L^{p_0}(0,T;H)$:
\[
Z(t)ADu=e^{At}ADu+\intt N(t-s)\left [ Z(s)AD u\right ]\ZD s\,. 
\]
its solution belongs to $L^{p_0}(0,T;H)$ and depends continuously on $u\in U$.

The proof of \emph{statement~\ref{I4TeoPropriZt}} is similar to that 
  of statement~\ref{I5bTeoPropriZt}: it is obtained by combining   property~\ref{I3TeoPropriZt} of this theorem,  the representation~(\ref{Eq:VoltePerZ}) of $Z(t)$,
statement~\ref{I22Teo:DiseFonda} of Theorem~\ref{Teo:DiseFonda} 
  and Young inequalities.\zdia

The following  result and its consequence     Theorem~\ref{Teo:PROpreidelRappoINCREvU} are used in the derivation of the Riccati equation.

\begin{Theorem}\ZLA{Teo:limRappINCREz}
Let $\ZEP>0$ and let $v \in\Dom\,(-A)^{\ZEP}$. Let $T>0$. There exists $p>1$ such that 
\begin{subequations}
\begin{equation}\ZLA{equaDELTeo:limRappINCREzVALUElim}
\lim _{h\to 0 } \dfrac{Z(t+h)v-Z(t)v}{h} =Z'(t)v=(A+I)Z(t)v-\intt e^{-(t-s)} Z(s)v\ZD s 
\end{equation}
and the limit exists for every $t\geq 0$ (if $t=0$ then we consider the right limit). The convergence is  in $L^p(0,T;H)$  (any $T>0$). 
 Furthermore, let $h\in (0,1)$.
 there exists $M=M_{T,\ZEP } $ (which does not depend on $h $) such that
 \begin{equation}\ZLA{equaDELTeo:limRappINCREz}
\left \| \dfrac{Z(\cdot+h)v-Z(\cdot)v}{h}\right \|_{L^p(0,T;H)} \leq M\|(-A)^{\ZEP}v\|_H\qquad \forall t\in [0,T]\,.
 \end{equation}
\end{subequations}
If $\ZEP\geq 1$ then $Z(t)v$ is of class $C^1$.
\end{Theorem}
\zProof First we prove the result in the special case  $Z(t)= e^{At} $. By assumption, $v=(-A)^{-\ZEP }\tilde v$ and $e^{A\cdot}A(-A)^{-\ZEP }\tilde v\in L^p(0,T;H)$ (any $p\in [1,1/(1-\ZEP))$\,).

 We consider the right limit at  $t\geq 0$. We have
\begin{subequations}
\begin{multline}\ZLA{Eq1diTeo:limRappINCREz}
\lim _{h\to 0^+}\dfrac{e^{A(t+h)}(-A)^{-\ZEP }\tilde v-e^{At}(-A)^{-\ZEP }\tilde v}{h} \\=
\lim _{h\to 0^+}e^{At}A(-A)^{-\ZEP }\dfrac{1}{h}\int_0^h e^{As }\tilde v\ZD s 
= 
e^{A t}A(-A)^{-\ZEP }\tilde v=e^{At}Av
\end{multline}
and there exists $M>0$ such that
\begin{equation}\ZLA{Eq1diTeo:limRappINCREzPOI}
\left \| e^{A t}A(-A)^{-\ZEP }\dfrac{1}{h}\int_0^h e^{As }\tilde v\ZD s\right \| \leq \dfrac{M }{t^{1-\ZEP}}\|\tilde v\|\,.
\end{equation}
The  equality~(\ref{Eq1diTeo:limRappINCREz}) is~(\ref{equaDELTeo:limRappINCREzVALUElim}) while~(\ref{Eq1diTeo:limRappINCREzPOI}) implies~(\ref{equaDELTeo:limRappINCREz}) when $Z(t)=e^{At}$.
Pointwise convergence, inequality~(\ref{Eq1diTeo:limRappINCREzPOI})    and dominated convergence show  that  the incremental quotient converges in the space $L^p(0,T;H)$.

The left limit at $t>0$ is studied in a similar way: the operator $e^{At}$ in~(\ref{Eq1diTeo:limRappINCREz}) is replaced by $e^{A(t-h)}$ which is strongly continuous and bounded on compact intervals.

Now we use the representation of $Z(t)$ in the third line of~(\ref{Eq:VoltePerZ}). The incremental quotient~(\ref{equaDELTeo:limRappINCREzVALUElim}) is the sum of three incremental quotients. The  limit of the first incremental quotient 
 is~(\ref{Eq1diTeo:limRappINCREz}) and this equality is used to prove inequality~(\ref{equaDELTeo:limRappINCREz}) in the general case and to compute the limits (in $L^p$) of the second and third incremental quotient which    are    respectively
\begin{align}\ZLA{Eq2diTeo:limRappINCREz}
&Z(t)v +A\intt e^{A(t-s)}Z(s)v\ZD s\,,
\\
&
\ZLA{Eq3diTeo:limRappINCREz}
-\intt e^{-(t-s)} Z(s)v \ZD s-A\intt e^{A(t-s)}\int_0^s e^{-(s-r)}Z(r)  v\ZD r\,\ZD s\,.
\end{align}
 The result is obtained by summing~(\ref{Eq1diTeo:limRappINCREz})  (\ref{Eq2diTeo:limRappINCREz}) and~(\ref{Eq3diTeo:limRappINCREz}). 
 
 The statement concerning the case $\ZEP\geq 1$ is obvious.\zdia
\end{subequations}

The previous theorem can be reformulated as follows: 

\begin{Corollary}\ZLA{coroConveIncreQuoZ}
  Let $\ZEP>0$. There exists $p>1$ such that the incremental quotient of $Z(t)v$ 
 converges in $L^p(0,T;(\Dom\,(-A)^{\ZEP})')$ for every $v\in H$ and there exists $M$ such that
 \begin{equation}\ZLA{Eq1diTeo:limRappINCREzPOIinSpazioduale}
\left \| \dfrac{Z(\cdot+h)v-Z(\cdot)v }{h} \right \|_{L^p(0,T;(\Dom\,(-A)^{\ZEP})')} \leq M\| v\|\,.
\end{equation}
 \end{Corollary}

We use $Z(t)$ to obtain a variation of constants formula for the solutions of~(\ref{eq:laEQsu(ztT)}):
\begin{Theorem}\ZLA{TheoVARIAConstFORM}Let\footnote{The special case of our interest is the case $y(t)=e^{-(t-\zt)}\hat y_\zt$, $\hat y_\zt\in (\Dom\,A)'$ but here we consider the general case of~(\ref{eq:laEQsu(ztT)}) with $e^{-(t-\zt)}\hat y_\zt$ replaced by any $ C^1([\zt,T],(\Dom\,A)')$-function $y$.} $y\in C^1([\zt,T],(\Dom\,A)')$,  $\Xi_\zt$ be as in~(\ref{EQpreDEFIstate}) and let $u$ belong to $ L^2(0,T;U)$.
The unique solution of problem~(\ref{eq:laEQsu(ztT)}) is  
\begin{multline}
\ZLA{VariatCONSTANTA}
v(t)=Z(t-\zt)\hat v_\zt  
 -\int_\zt^t Z(t-s)e^{-(s-\zt)}\int_0^{\zt} e^{-\nu } \xi_\zt(\zt-\nu)\ZD\nu\,\ZD s
 +\int_\zt ^t Z(t-s)y(s)\ZD s\\
-\int_\zt ^t Z(t-s)AD u(s)\ZD s 
\end{multline}
 
\end{Theorem}
\zProof To prove~(\ref{VariatCONSTANTA}), we put
\[
f(t)=
- e^{-(t-\zt)}\int_0^\zt e^{- s }  \xi_\zt(\zt-s)\ZD s +y(t)-AD u(t)\,.
\]
The (weak) solution of~(\ref{eq:laEQsu(ztT)}) is  the solution of the following Volterra integral equation:
\begin{equation}\ZLA{(eqUaPerVconF}
v(t)=e^{A(t-\zt)} \hat v_\zt+\int_\zt ^t e^{A(t-s)} f(s)\ZD s+\int_\zt^t N(t-s) v(s)\ZD s\,.
\end{equation}
The second addendum belongs to $L^r(\zt,T;H)$.
We use~(\ref{Eq:VoltePerZ}) to express $e^{At}\xi$ in terms of $Z(t)\xi$:
\[
e^{At}\xi=
Z(t)\xi-\intt N(t-s)Z(s)\xi\ZD s\,.
\]
We replace in the first two addenda of the right side of~(\ref{(eqUaPerVconF}) and we find the following equality:

\begin{multline*}
\left \{
v(t)-Z(t-\zt)\hat v_\zt-\int_\zt ^t Z(t-s) f(s)\ZD s
\right \}\\
=\int_\zt^t N(t-s)\left \{
v(s)-Z(s-\zt)\hat v_\zt-\int_\zt ^s Z(s-\nu)f(\nu)\ZD \nu
\right \}\ZD s\,.
\end{multline*}
This is a Volterra integral equation on the interval $[\zt,T]$ whose ``unknown'' is the brace. The known term is zero. So, the solution is zero, i.e the brace is zero and this is the required 
equality~(\ref{VariatCONSTANTA}).\zdia

The regularity properties of the function $v(t) $ in~(\ref{VariatCONSTANTA}) are stated in Theorem~\ref{Teo:preREGOLARITAstato}.

We present a 
consequence of Theorem~\ref{Teo:limRappINCREz} (partially stated in   Theorem~\ref{teoLasolUdebDeVolte}) which has a crucial role in the derivation of the Riccati equation. 
\begin{Theorem}\ZLA{Teo:PROpreidelRappoINCREvU}
Let    $\hat v_\zt\in\Dom\,A $, 
$y(t)=e^{-(t-\zt)}\hat y_\zt $ with
$\hat y_\zt\in H$,   $ \xi_\zt\in L^2(0,\zt;H)$
and let $u\in L^\ZIN(\zt,T;U)$. Let  $\ZEP>0$ and let $v$ be the function in~(\ref{VariatCONSTANTA}).
 
  There exists a constant $M=M_{\ZEP,T}$  which does not depend on  $\zt$ and $h$) such that
\begin{equation}\ZLA{eq:LimiQuozIndFivU}
\begin{array}{l}
\displaystyle \left \| \dfrac{v (t+h)-v (t)}{h}\right \|_{C([0,T],(\Dom\,(-A)^{\ZSI+\ZEP})')}\\[5mm]
\qquad<M\left\{
\hat v_\zt\|_{H}+\left \|  \hat y_\zt +\int_0^\zt e^{-(\zt-r)} \xi_\zt(r)\ZD r\right \|_H+
\|u\|_{C([\zt,T];U)}\right\}\,,
\\[3mm]
\displaystyle
 \lim _{h\to 0^+}\dfrac{v (t+h)-v (t)}{h} =v' (t)\quad \mbox{in $C([0,T],(\Dom\,(-A)^{\ZSI+\ZEP})')$}\,.
 \end{array}
\end{equation}
\end{Theorem}
\zProof The  integration by parts in~(\ref{Eq:DefiRefeSigy}) on the interval $[\zt,T]$ and with $\hat y$ replaced by
\[
 \hat y_\zt+\int_0^\zt e^{-(\zt-s)}  \xi_\zt(\zt-s)\ZD s
\]
shows that it is sufficient that we prove the statement for  the addendum 
\[
v_u(t)=-\int_\zt^t Z(t-s)ADu(s)\ZD s \,.
\]
  
We must study the incremental quotient
\begin{multline*}
\dfrac{(-A)^{-\ZSI-\ZEP}\left( v_u(t+h)-v_u(t)\right )}{h}=
\dfrac{1}{h}\int_t^{t+h}Z(t+h-s)(-A)^{-\ZEP} (-A)^{1-\ZSI}Du(s)\ZD s\\
+\int_\zt^t\dfrac{[Z(t+h-s)-Z(t-s) }{h}(-A)^{-\ZEP} (-A)^{1-\ZSI}Du(s)\ZD s\ZD s\,.
\end{multline*}
Boundedness and convergence of the first addendum in the right hand side follows since the integrand is continuous.  Boundedness and convergence of the second addendum follows from Young inequalities since $(-A)^{1-\ZSI}Du\in C([0,T];H)$ and  the incremental quotient belongs to a space $L^p$ (with $p>1$) with norm   bounded by a constant $M$ which does not depend on $\zt$, $t\in [\zt,T]$ and $h$ 
(see Theorem~\ref{Teo:limRappINCREz}).

The limit is   given by the equality~(\ref{equaDELTeo:limRappINCREzVALUElim}).\zdia

\subsection{\ZLA{sec:ContinuityOPTIMALcontrol}The optimal control}
We introduce  
\begin{equation}\ZLA{eq:DefiOPERfonDam}
\left\{\begin{array}{l}
\displaystyle K(t)=Z(t)AD\,,\ K\ZCD\in  \mathcal{L}(U,L^{p_0}(0,T; \mathcal{L}(U,H)))\,,\\[2mm]
\displaystyle (L_\zt y)(t)=\int_\zt^t Z(t-s) y(s)\ZD s \\[1mm]
\displaystyle \qquad L_\zt\in \mathcal{L}(C([\zt,T];(\Dom\,A)'))\cap \mathcal{L}(C^1([\zt,T];(\Dom\,A)'),C([\zt,T];H))\,,
\\[2mm]
\displaystyle     \left (\Gzt \Xi_\zt\right )(t)
=
Z(t-\zt)\hat v_\zt 
  -\int_\zt^t Z(t-s)\int_0^{\zt} e^{-(s-r)} \xi _\zt(\zt-r)\ZD r\,\ZD s \\[1mm]
  %%%%%%%%%%
%\displaystyle  =Z(t-\zt)\hat v_\zt 
%  -\int_\zt^t Z(t-s)e^{-(s-\zt)}\int_0^\zt e^{-r} \xi_\zt( r)\ZD r\,\ZD s  \\[1mm]
%\displaystyle =
% Z(t-\zt)\hat v_\zt 
% -\int_\zt^t Z(t-s)e^{-(s-\zt)}\int_0^\zt e^{-r} \xi (\zt- r)\ZD r\,\ZD s \\[1mm]
\displaystyle \qquad   \Gzt\in\mathcal{L}( M^2_\zt, C([\zt,T];H))\,,\\[2mm]
  %%%%%%%%%%%%
\displaystyle    \left (\Lzt u\right )(t)=-\int_\zt^t K(t-s)u(s)\ZD s  \\[1mm]
\displaystyle  \qquad  \Lzt\in\mathcal{L}( L^2(\zt,T;U),L^r(\zt,T;H))\subseteq 
   \mathcal{L}( L^2(\zt,T;U), L^2(\zt,T;H))\,,\\[2mm]
\displaystyle \hat Y_\zt(t)=\int_\zt ^t Z(t-s) e^{-(s-\zt)}\hat y_{\zt} \ZD s=\left [L_\zt e^{ -(\cdot-\zt)}\hat  y_\zt \right ](t)   \\[1mm] 
\displaystyle \qquad \hat Y\ZCD\in C([\zt,T];H)  \quad \forall \hat  y_{\zt}\in (\Dom\,A)'\,, \\[2mm] 
\displaystyle h_{\zt}(t)=\left (\Gzt\Xi_\zt\right )(t) +\hat Y_\zt(t)\\
\displaystyle =Z(t-\zt)\hat v_\zt -\int_\zt ^t Z(t-s) e^{-(s-\zt)}\left [\int_0^\zt e^{-(\zt-r)}\tilde\xi (r)\ZD r-\hat y_\zt\right ]\ZD s  \\
\displaystyle \qquad h_\zt\ZCD\in C([\zt,T];H)\,.
\end{array}\right.
\end{equation}
With these notations, the solution $v$ of~(\ref{eq:laEQsu(ztT)}) takes the form
\begin{multline}
\ZLA{Eq:soluCONoperatpri}
v(t;\zt,\Sta_\zt,u)=v(t;\zt,(\Xi_\zt,\hat y_{\zt}),u)
\\=\left (\Gzt \Xi_\zt\right )(t)+\left (\Lzt u\right )(t)+\hat Y_\zt (t)=h_{\zt}(t)+\left (\Lzt u\right )(t)\,.
\end{multline}

We recall the notation 
\begin{align*}
&\mbox{
$u^+_\zt= u^+_\zt(\cdot;\Sta_\zt)= u^+_\zt(\cdot;(\Xi_\zt,\hat y_{\zt}))$ for the optimal control}\\
&\mbox{and $v^+_\zt=v^+_\zt(\cdot;\Sta_\zt)=v^+_\zt(\cdot;(\Xi_\zt,\hat y_{\zt}))$ for the corresponding solution of~(\ref{eq:laEQsu(ztT)}).} 
\end{align*}
A standard variational approach   gives the following relation between the optimal control and the corresponding solution $v_\zt^+$:
  \begin{subequations}
 \begin{align}
 \ZLA{eq:controOTTIMO}
 u_\zt^+&=-\Lzt ^*v_\zt^+\quad{\rm} i.e.\quad  u_\zt^+(t)=-\int_t^T K^*(s-t)v_\zt^+(s)\ZD s\\
\ZLA{eq:evolOTTIMA}
v_\zt^+&=h_\zt-\Lzt \Lzt^*u_\zt^+\\
\nonumber &\mbox{and so}\\
 \ZLA{eq:evolOTTIMAOPEN}
 v_\zt^+&=(I+\Lzt\Lzt^*)^{-1}h_\zt\\
 \ZLA{eq:controOTTIMOOPEN}
 u_\zt^+&=-\Lzt^* (I+\Lzt\Lzt^*)^{-1}h_\zt=- (I+\Lzt^*\Lzt )^{-1}\Lzt^* h_\zt\,.
 \end{align}
 We replace~(\ref{eq:evolOTTIMAOPEN}) and~(\ref{eq:controOTTIMOOPEN}) in the cost~(\ref{eq:costoWAVEdampSU(ztT)}) and we find   
\begin{multline}
\ZLA{ValueFunction}
W_\zt(\Sta_\zt)%=W_\zt(\Xi_\zt,\hat y_{\zt})
=\min_{u\in L^2(\zt,T;U)} J_\zt(\Sta_{\zt};u) 
=\int_\zt^T [\|v^+_\zt(s)\|_H^2+\|u_\zt^+(s)\|_U^2 ]\ZD s \\
=\int_\zt^T\ZL [(I+\Lzt\Lzt^*)^{-1} h_\zt](s),h_\zt(s)\ZR_H \ZD s\,.
\end{multline}
 \end{subequations}

The function $\Sta_\zt\mapsto W_\zt(\Sta_\zt)$ is the  \emph{value function of the problem~(\ref{eq:laEQsu(ztT)})-(\ref{eq:costoWAVEdampSU(ztT)}) at the time $\zt$.}

 Now we observe that $v=v_{\zt}^+$ and $  u_\zt^+$ solve the following system of Fredholm integral equations on $[\zt,T]$: 
 
   \begin{equation}
    \ZLA{eq:controOTTIMOOPENalterna} \begin{array}{l}
    \displaystyle
    v(t)=-\int_\zt^t K(t-s)\left[\int_s^T K^*(\nu-s)v(\nu)\ZD\nu\right]\,\ZD s+h_\zt(t)\,,\\
 \displaystyle   u(t)=-\int_t^T K^*(s-t) \left [\int_\zt ^s K(r-s)  u(r)\ZD r\right] \ZD s\\
    \displaystyle \qquad -\int_t^T K^*(s-t) h_\zt
   (s)\ZD s\,.
   \end{array}
 \end{equation}

The function $ h_\zt $ is continuous. We apply  Theorem~\ref{teoAPPEregolarPRPsoluFred} in Appendix~\ref{AppePROpreHzt} (in particular
Statement~\ref{I1teoAPPEregolarPRPsoluFred})  to the Fredholm integral equation of $v^+_\zt$. We obtain continuity of $t\mapsto v^+_\zt(t)$  on $[\zt,T]$. Then we use~(\ref{eq:controOTTIMO}). We obtain:
\begin{Theorem}\ZLA{TEOValue functionCONTI}
Let \ $\Sta_\zt \in  M^2_\zt\times (\Dom\,A)'$. The following properties hold:
\begin{enumerate}
\item  $v_\zt^+(\cdot;\Sta_\zt)$ is continuous on $[\zt,T]$
and 
$
\Sta_\zt\mapsto v_\zt^+(\cdot;\Sta_\zt)
$ 
belongs to $\mathcal{L}(M^2_\zt\times ( \Dom\, A)', C([\zt,T];H))$. 
\item $u_\zt^+(\cdot;\Sta_\zt)$ is continuous on $[\zt,T]$
and 
$
\Sta_\zt\mapsto u_\zt^+(\cdot;\Sta_\zt)
$  
 belongs to $\mathcal{L}(M^2_\zt\times ( \Dom\, A)', C([\zt,T];U))$. 
\end{enumerate}

So we have
\begin{equation}
\ZLA{Value functionCONTI}
W_\zt(\Sta_\zt)
=\min_{u\in L^2(\zt,T;U)} J_\zt(\Sta_\zt;u)=
 \min_{u\in C([\zt,T];U)} J_\zt(\Sta_\zt;u)\,.
 \end{equation}
\end{Theorem}

The strict inequality $p_0>1$ is  used in the proof of Theorem~\ref{TEOValue functionCONTI}.

 \subsection{\ZLA{SectSTATEdefi}The state of the dynamical system~(\ref{eq:laEQsu(ztT)})}
 
The definition of the state of a dynamical system is deeply analyzed in~\cite{KalmanFalbArbib1969}. Without entering into details, we can roughly state: a dynamical system is a system which evolves in time  and such that unicity of solution (in the future) holds; i.e. at every time $t$ (say of a  right half line) it is possible to specify a certain object which uniquely determines the solution in the future. This object is the state at time~$t$.

We note that when defining the state of a control system, we assume that the control is    fixed and known. 
\begin{Remark}
{\rm
Very roughly speaking, we can think at the state at time $t$ as the set of the ``initial conditions which, when assigned at the time $t$, uniquely identify the evolution of the system in the future''. The fact that this intuition is rough is seen by considering the fact that  \emph{the state at the time $t$ is never uniquely identified.} So, an important problem has been the    identification of the  \emph{minimal state.} The search of the minimal state of lumped systems was initiated in the  30'-thies of the last century in the context of the synthesis of electrical circuits. The definitive results are in~\cite{KalmanFalbArbib1969}  then extended to certain general classes of distributed systems in the '70-this (see~\cite{FuhrmannBOOK81}). This quite general and abstract approach was not prosecuted. Instead, important and interesting results have been obtained for the class of systems with memory 
(see~\cite{DelfourManitOPERf1,FabrizioARCH2010}). An application of the 
results in~\cite{DelfourManitOPERf1} to the quadratic regulator problem of systems with delays   in $\zzr^n$ is in~\cite{DelfLEEmAnFOpeRICCA2}.\zdia
}
\end{Remark}
   
Here we study (one possible definition of) the state of   system~(\ref{eq:laFORMintegDIvCONuVOLTE}) \emph{with $L^\ZIN$-controls.} Bounded controls are sufficient to study the quadratic regulator problem   since the optimal control is even continuous. 
A comment on the state when the controls are square integrable   is in Remark~\ref{REMAsuL2CONTROLestato}.

 We recall that when defining the state, the control is assumed to be known. So the solution  for $t>\zt$  is uniquely identified by  $\hat v_\zt\in H$, $  \xi_\zt\in L^2(0,\zt;H)$ and $\hat y_\zt\in (\Dom\,A)'$ and so by
\begin{equation}\ZLA{EQpreDEFIstateREDEF}
 \Sta_\zt=(\Xi_\zt,\hat y_\zt)\in M^2_\zt\times (\Dom\,A)'\,. 
\end{equation}
\emph{We choose $\Sta_\zt=(\Xi_\zt,\hat y_\zt)\in M^2_\zt\times (\Dom\,A)'$   as  the \emph{state} at time $\zt$ and $M^2_\zt\times (\Dom\,A)'$ as the \emph{state space.}}

 This observation can be repeated for every $t_1\in (\zt,T]$. Eq.~(\ref{eq:laEQsu(ztT)}) on $[t_1,T]$ takes the form
 \begin{multline}
\ZLA{eq:laEQsu(ztT)maDAzt1}
v'(t)=Av(t) +v(t)-\int_{ t_1}^t e^{-(t-s)} v(s)\ZD s-AD u(t)
\\
-e^{-(t-t_1)}\int_0^{ t_1} e^{-s}  \xi_{ t_1}(t_1-s) \ZD s +e^{-(t-t_1) }\hat y _{t_1}\,,\qquad 
v(t_1)=\hat v_{ t_1} \,.
\end{multline}
\begin{Remark} {\rm In general, $\Sta_{t_1}=(\Xi_{t_1},\hat y _{t_1})$ is an \emph{arbitrary} element of $M^2_{t_1}\times (\Dom\, A)'$. In this case the 
solution of~(\ref{eq:laEQsu(ztT)maDAzt1}) is not related to that of~(\ref{eq:laEQsu(ztT)}). The solution of~(\ref{eq:laEQsu(ztT)maDAzt1}) is the restriction to $[t_1,T]$ of that of~(\ref{eq:laEQsu(ztT)}) when
 
\begin{multline} 
\ZLA{eq:DEFIvtANTE}
\left\{\begin{array}{lll}
\hat v_{t_1}&=&v(t_1^-)\\
\hat y _{t_1}&=& y(t_1)= e^{-(t_1-\zt)}\hat y_{\zt} 
\end{array}\right.  
\\
   \xi_{\zt_1}(s) =\left\{\begin{array}{ll}
\mbox{the function $ \xi_\zt(s)$ in~(\ref{EQpreDEFIstateREDEF})} &\mbox{when $s\in (0,\zt)$}
\\
%}\\
\mbox{$v(s)$, the solution of~(\ref{eq:laEQsu(ztT)}),}
&\mbox{when $s\in(\zt,t_1 )$}\,.\zdiaform
\end{array}\right.
\end{multline}
 }\end{Remark}

The following definition recapitulates the previous considerations.  
The time $t_1\in [0,T]$ being   \emph{arbitrary,}  in the next definition it is denoted $t$.
\begin{Definition}
{\rm
Let $0\leq\zt\leq t\leq T$. Let   $\Xi(t)=(v(t),\xi_t(t-\cdot))\in M^2_t$   defined in~(\ref{eq:DEFIvtANTE}) , $y(t)=e^{-(t-\zt) }\hat y_{\zt}$, $\hat y_{\zt}\in (\Dom\, A)'$.

The \emph{state space at time $t$} of system~(\ref{eq:laEQsu(ztT)}) is $M^2_t\times (\Dom\,A)'$ and the \emph{state at the time $t$} is $\Sta(t)=(\Xi(t), y(t))$.

When convenient, the state at time $t$ is also called ``the initial condition at time $t$'' (``initial'' respect to the future evolution of the system).\zdia

}
\end{Definition}
\begin{Remark}\ZLA{Rema:StatoEbponCONt1}
 {\rm
 We observe
 \begin{itemize}
 \item
 the reason why   \emph{we assumed boundedness of the control} is now clear: thanks to this assumption the solution $v$ is an $H$-valued continuous function for $t\geq \zt$ and $v(t)$ can be computed in $H$.
 \item we note a discrepancy in the notation: the state $\Sta(\zt)$ at the initial time $\zt$ is denoted also $\Sta_\zt$  and similar for the components $ v$, $\xi $, $y$ and   $\Xi= (v,\xi)$. 
  \item
 the component $y$ of the state does not depend on $\Xi_\zt$ and $u$. 
 \item We stress a fact which is typical of systems with memory on a finite interval of time: the state space does depend on time. When instead the system is studied on a half line, as in the Dafermo's semigroup approach (and with $y=0$)   the state space
is $H\times L^2_h(0,+\ZIN;H)$  and it does not depend on time. The index $_{h}$ denotes the fact that a suitable weight (which depends on the memory term) is used to define the $L^2$-norm.\zdia
 \end{itemize}
 }
 \end{Remark}

The results in Sect.~\ref{SECT:soluzione} are easily adapted to prove:

\begin{Theorem}
Let $0\leq \zt\le t\leq T$. The transformation  

\[ 
(\Sta_\zt ,u)\mapsto\Sta(t;\zt,\Sta_\zt, u)=  (\Xi(t;\zt,\Sta_\zt , u),y(t;\zt,\hat y_\zt))
\]
is (linear and) continuous from $M^2_\zt\times(\Dom\, A )' \times L^\ZIN(\zt,T;U)$
to $M^2_t \times H\subseteq M^2_t \times (\Dom\,A)'$.
\end{Theorem}
 
The last component $y$ of the state\footnote{Not to be confused with the 
function $y(t) $ in~(\ref{eq:NoveltyDAMPING}).}
is denoted $y(t;\zt,\hat y_\zt)$ since it does not depend either on $\Xi_\zt$ or on the control $u$.

Unicity of solutions can be reformulated as follows: let $0\leq \zt< t_1< t\leq T$. We have
\begin{multline}\ZLA{Eq.secDynaPro:DYnaProsoluANTE}
 \Sta(t;t_1,\Sta(t_1;\zt; \Sta_\zt,u),   u_1) 
 =\Sta(t;\zt,\Sta_\zt,
( u  \wedge_{ t_1 } u_1 )
)
\\
\mbox{where}\quad 
 \left (u  \wedge_{ t_1 } u_1\right )(s)
=\left\{\begin{array}{lll}
u(s) &{\rm if}& \zt< s< t_1\\
u_1(t) &{\rm if}& t_1<s<t\,.
\end{array}\right.
\end{multline}

 \begin{Remark}\ZLA{REMAsuL2CONTROLestato}
 {\rm
 In the definition of the state we assumed boundedness of the control in order to have continuity of the solution.  
 It is possible to define the state when $u$ is square integrable   at the expenses of replacing $H$ with $(\Dom(-A)^\ZSI)'$.\zdia
 } 
 \end{Remark}

\subsection{\ZLA{SubsEQUAstate}A differential  equation for the state}
In order to explain the synthesis of the optimal control, it is convenient to write an equation   which is satisfied (in a weak sense) by the state.  As in Sect.~\ref{SectSTATEdefi} we confine ourselves to the case that  $u$ is bounded, since this is the case we are interested in when studying the synthesis of the optimal control. We consider the evolution of the system on an interval $[\zt,T]\subseteq [0,T]$.

The function $y(t)$ solves the linear first order differential equation in $(\Dom\,A)'$
\begin{equation}
\ZLA{eqDINAMdiY}
y'=-y\,,\qquad y(\zt)=\hat y_{\zt}   \,.
\end{equation}
So, we write an equation for the component $\Xi(t)=(v(t),\xi_t)$ ($v$ 
solves~(\ref{eq:laEQsu(ztT)}) and $\xi _t$ is given by~(\ref{eq:DEFIvtANTE}))  
with an  arbitrary  $y\in C([\zt,T];(\Dom\,A)' )$. Then we  add a component which correspond to the equation~(\ref{eqDINAMdiY}).

We adapt the ideas first presented  in~\cite{PandolfiVOLTERRAieee18,Pandolfi24EECTtrackingRd}. 

We use~(\ref{eq:DEFIvtANTE}) 
and we denote either $\xi_t(s)$ or $\xi(s,t)$ thefollowing extension of the function $\xi_\zt$:
 
\begin{subequations}
\begin{equation} \ZLA{eq:DEFIXIt}
\xi_t(s)=\xi(s,t)=\left\{\begin{array}{lll}
v(t-s)
&{\rm if}& 0<s<t-\zt\\
 \xi_\zt(t-s)&{\rm if}& t-\zt<s<t\,.
\end{array}\right.
\end{equation}

Note that~(\ref{eq:DEFIXIt}) is equivalent to extending $v$ to $(0,\zt)$  as follows:
\begin{equation}
\ZLA{eq:DEFIvt}
v(s)=\left\{\begin{array}{l}
\mbox{$  \xi_\zt(s )$ if $s\in(0,\zt)$}\\
\mbox{$v(s;\zt,\Xi_\zt,  y,u)$ if $s>\zt$}\,.
\end{array}\right.
\end{equation}
 \end{subequations}

We use~(\ref{eq:DEFIXIt}) to represent the sum of the integral terms in~(\ref{eq:laEQsu(ztT)maDAzt1})   as
\begin{align*}
&
-e^{(-t-t_1)}\int_0^\zt e^{-s}  \xi_{t_1}(t_1-s)\ZD s-\int_\zt^t e^{-(t-s)} v(s)\ZD s 
=-\mathcal{E}_t \xi(\cdot,t)\\
&
 \mathcal{E}_t\,: L^2(0,t;H) \mapsto H\,,\qquad[ \mathcal{E}_t \xi(\cdot,t)](t )=
  \intt e^{-s}\xi(t-s)\ZD s\,.%= \intt e^{-s}\xi_t(s)\ZD s\,.
\end{align*}
The function $\xi(s,t)=v(t-s)$ ($v$ extended as in~(\ref{eq:DEFIvt})) is the (weak) solution of
 \begin{equation}\ZLA{eqDelloSHIFTprestato}
\partial_t \xi =-\partial_s\xi\,,\qquad \xi(0,t)=v(t)\,,\qquad t> 0\,,\ 0<s<t  \,.
 \end{equation}
On the interval of time $(\zt,T)$ we consider the system of the equations~(\ref{eq:laEQsu(ztT)maDAzt1}) and~(\ref{eqDelloSHIFTprestato}). With the notations just introduced, Eq.~(\ref{eq:laEQsu(ztT)maDAzt1}) (with $t_1=\zt$)  takes the form
\begin{equation}\ZLA{eq:laFORMintegDIvCONuINTEGROdiffDEBprestato}
v'(t)=(A+I)v(t) -\mathcal{E}_t\xi(\cdot,t)-AD u(t)+y(t)
  \,,\quad 
v(\zt )=\hat v_{\zt}\,.
\end{equation}
It is clear that any (weak) solution of~(\ref{eq:laFORMintegDIvCONuINTEGROdiffDEBprestato})
has the property that $\Xi(t)=(v(t), v(t-\cdot ))$
 solves the system~(\ref{eqDelloSHIFTprestato})-(\ref{eq:laFORMintegDIvCONuINTEGROdiffDEBprestato})
and conversely. In fact, a (weak) solution of~(\ref{eqDelloSHIFTprestato})
is of necessity $\xi(s,t)=v(t-s)$ and so $v$ solves~(\ref{eq:laFORMintegDIvCONuINTEGROdiffDEBprestato}).

We introduce the operators $ \widetilde{\mathcal{A}}_t$ and $ \widetilde{\mathcal{B}}_t$ as follows:
 
\begin{align*}
& \widetilde{\mathcal{A}}_t\,: M^2_t\mapsto M^2_t\,,\quad \widetilde{\mathcal{A}}_t\Xi= \widetilde{\mathcal{A}}_t\left[\begin{array}{c}
  v
\\
 \xi\ZCD
\end{array}
\right]=\left[
\begin{array}{c}
(A+I) \hat v-\mathcal{E}_t  \xi_t\\
- \partial_s\tilde\xi_t
\end{array}
\right]\\ 
&
\Dom\, \widetilde{\mathcal{A}}_t=\{
(\hat v, \xi_t)=(\hat v, \xi(t-\cdot))\,,\quad \hat v\in \Dom\,A\,, \quad  \xi\in H^1(0,t;H)\,,\quad \xi(0)=\hat v
\}\\
& \widetilde{\mathcal{B}}_t\,: H\mapsto M^2_t\,,\qquad 
 \widetilde{\mathcal{B}}_t u=\left[\begin{array}{c}
Bu\\
0
\end{array}\right]
\quad B=-AD\,.
\end{align*}
 System~(\ref{eqDelloSHIFTprestato})-(\ref{eq:laFORMintegDIvCONuINTEGROdiffDEBprestato}) takes the   form
\begin{equation}
\ZLA{Eq:sisteSOPAZIOstati}
\Xi'= \widetilde{\mathcal{A}}_t\Xi+ \widetilde{\mathcal{B}}_t u+Y(t)\,,\quad  Y(t)=\left[
\begin{array}{c}
y(t)
\\0 
\end{array}
\right]\,.
\end{equation}
This system has to be consider on $(\zt,T)$, $\zt\in [0,T]$ and the initial condition at the time $\zt$ is  
\begin{equation}\ZLA{INIdiEq:sisteSOPAZIOstati}
\Xi(\zt)=\Xi_\zt\,.
\end{equation}
The differential equation~(\ref{eq:laFORMintegDIvCONuINTEGROdiffDEBprestato})  has to be interpreted in the weak form~(\ref{eq:laFORMintegDIvCONuINTEGROdiffDEB}) and~(\ref{eqDelloSHIFTprestato}) has to be interpreted in the weak form
\[
\dfrac{\ZD}{\ZD t}\ZL\xi(\cdot,t),\eta\ZCD\ZR=\ZL\xi(\cdot,t),\partial_s \eta\ZCD\ZR-\ZL v(t),\eta(0)\ZR
\]
for every $\eta\in H^1(0,t;H)$ such that $\eta(t)=0$.
 
 Statement~\ref{I1CTeo:preREGOLARITAstato} of Theorem~\ref{Teo:preREGOLARITAstato} can be reformulated as follows:
\begin{Theorem}
\ZLA{Teo:REGOLARdelloSTATOparzia}
If $u=0$, $\Xi_\zt\in\Dom\,\widetilde{\mathcal{A}}_\zt$, $y(t)=e^{-(t-\zt)}\hat y_\zt$ with $\hat y_\zt\in H$ then $\Xi(t)\in  C([0,T];\Dom\,\widetilde{\mathcal{A}}_t).$
 
\end{Theorem}

 Now we combine the equations~(\ref{Eq:sisteSOPAZIOstati}) of $\Xi $ and the equation~(\ref{eqDINAMdiY}) of $y $. We get the following \emph{state space representation} of system~(\ref{eq:laEQsu(ztT)maDAzt1}):
  \begin{equation}\ZLA{StateEQUAcompleta}
 \left. \begin{array}{lll}
\displaystyle \Sta'(t)&=&
\left[\begin{array}{c}
\Xi(t)\\y(t)
\end{array}\right]'=\mathcal{A}_{ t}\left[\begin{array}{c}
\Xi(t)\\y(t)
\end{array}\right]+ \mathcal{B}_{ t}u(t)\\[4mm]
&=&\mathcal{A}_{ t}\Sta(t)+\mathcal{B}_{ t}u(t) \qquad t\geq \zt
 \\[2mm]
\displaystyle \Sta(\zt)&=&\Sta_\zt\in M^2_{\zt}\times (\Dom\,A)'
\end{array}\right.
  \end{equation}
 where
 \begin{align*}
%\ZLA{eq:OperAESTES}
& \mathcal{A}_{ t}=\left[\begin{array}{ccc}
A+I&-\mathcal{E}_t&I\\
0&-\partial_s &  0\\
0&0 & -I
\end{array}\right] 
\quad\Dom\,\mathcal{A}_{ t}=( \Dom\,\widetilde{\mathcal{A}}_{ t})\times H\,,
\quad 
 \mathcal{B}_{ t}= \left[\begin{array}{c} \widetilde {\mathcal{B}}_{ t}\\
 0\end{array}\right]
%\\
%&  
%\quad\Dom\,\mathcal{A}_{ t}=( \Dom\,\widetilde{\mathcal{A}}_{ t})\times H\,,\\
%& \mathcal{B}_{ t}={\rm col}\left[\begin{array}{cc} \mathcal{B}_{ t}&0\end{array}\right]  \qquad \mbox{($ \mathcal{B}_t$ is  the \emph{input operator})}\,.
\end{align*}
($ \mathcal{B}_t$ is  the \emph{input operator}). 
%  $ \mathcal{B}_t$ is  the \emph{input operator.}
 
Theorem~\ref{Teo:REGOLARdelloSTATOparzia} gives:
\begin{Theorem}\ZLA{Teo:REGOLARdelloSTATO}
If $u=0$ and $\Sta_\zt\in \Dom\mathcal{A}_\zt$ then $\Sta_t\in\Dom\,\mathcal{A}_t$ for every $t\in[\zt,T]$.
\end{Theorem}

Finally, we introduce the following operators: 
      \begin{subequations}
\begin{equation}\ZLA{DefiPROIEXdaSTATO}
\begin{array}{l}
\mathcal{C}^1_t\Sta(t)=v(t)\,,\\
   \mathcal{C }^{2 }_t\Sta(t)=v(t-\cdot)=\xi(\cdot,t)\,,\\
\mathcal{C}^{(1,2)}_t=\left[\begin{array}{cc}\mathcal{C}^1_t& \mathcal{C}^2_t\end{array}\right]\quad \mbox{(so that $\mathcal{C}^{(1,2)}_t\Sta(t)=\Xi(t)
 $)}\,,
\\
\mathcal{C }^3_t\Sta(t)=y(t) \,.
\end{array} 
\end{equation} 
We call the projection $\mathcal{C}^1_t$ the \emph{output operator} since
the cost   is
\[
\int_\zt^T \left [ \|\mathcal{C }^1_t \Sta(t)\|^2+\|u(t)\|^2 \right ]\ZD t\,.
\]
Moreover we introduce  
  
 \begin{equation}
  \ZLA{eqDefiHzt}
H_\zt\in \mathcal{L}(H)\,,\quad H_\zt =(I+\Lambda_\zt \Lambda^*_\zt )^{-1} 
  \end{equation}
      so that
  
      \begin{align}\ZLA{Eq:DEFIhCOMEopeSUsato}
     & h_\zt\ZCD=\left(\Gamma_\zt  \mathcal{C}^{(1,2)}_\zt+L_\zt \mathcal{C}^{(3)}_\zt\right )\Sta_\zt\,,\\
    &  \ZLA{COSTOottimoCOMEfarmaQUADRATICA}
      W_\zt(\Sta_\zt)=\Bigl\ZL H_\zt  
      (\Gamma_\zt  \mathcal{C}^{(1,2)}_\zt+L_\zt \mathcal{C}^{(3)}_\zt  )\Sta_\zt, (\Gamma_\zt  \mathcal{C}^{(1,2)}_\zt+L_\zt \mathcal{C}^{(3)}_\zt )\Sta_\zt\Bigr\ZR _{L^2(\zt,T;H)}\,.
      \end{align}
 \end{subequations}
 \begin{Remark}\ZLA{RemaEXTEprojeC3}{\rm
If  $\hat y\in H$ then $y\in C^1([0,T];H)$ and   $\mathcal{C}_t^{(3)} \Sta'(t)$ belongs to $H$.\zdia
 }\end{Remark}
 
\section{\ZLA{SecDissiIneq}The dynamical properties of the optimal control and the dissipation inequality}
 
 Also in this section we assume that the control is bounded since this is sufficient to study the quadratic regulator problem. 

\begin{Remark}{\rm For clarity, we recall explicitly the following facts:
 
 \begin{enumerate}
 \item let $\zt\in [0,T)$ be the initial time. The state $\Sta(\zt)$ at the initial time $\zt$, i.e. the initial condition, is also denoted  $\Sta_\zt$.  
\item  the equality~(\ref{eq:DEFIvtANTE}) for the the $L^2$-component of the state at time $t\geq\zt$.

\item  The   optimal control is denoted $u^+_\zt(\Sta_\zt)= u^+_\zt(\cdot;\Sta_\zt)$. The corresponding solution is denoted $v^+_\zt(\cdot;\Sta_\zt)$
 i.e.
$
 v^+_\zt(\cdot;\Sta_\zt)=
v(t;\zt,\Sta_\zt, u^+_\zt(\Xi_\zt) )  
$.
 We use also the corresponding notations $\Sta^+_\zt(t\cdot;\Sta_\zt)$,   $\Xi^+_\zt(t\cdot;\Sta_\zt)$. 

\end{enumerate}
 }
 \end{Remark}
 A standard consequence of the unicity of the solution is:
 \begin{Theorem}\ZLA{Teo:dynaPROPoprimalCONTRP}
 let $0\leq \zt<t_ 0<t$. Let $\Xi_\zt$ be the initial condition at the time $\zt$. 
 For $t\in[t_0,T]$ we have:
 \begin{subequations}
 \begin{align}
 \ZLA{Eq:dynaPROPoprimalCONTRP}
 &u^+_{t_0}(t;\Sta^+_{\zt}(t_0;\Sta_\zt))=u^+_\zt(t;\Sta_\zt)  \,,\\
 %%%%
  \ZLA{Eq:dynaPROPozeroEVOL}
  &\Sta^+_{t_0}(t;\Sta^+_{\zt}(t_0;\Sta_\zt))=\Sta^+_\zt(t;\Sta_\zt)  \ {\rm so:} 
 & \left\{
   \begin{array}{l} 
%% \\
%  \ZLA{Eq:dynaPROPoprimalEVOL}
% %&
  v^+_{t_0}(t;\Sta^+_{\zt}(t_0;\Sta_\zt))=v^+_\zt(t;\Sta_\zt)  \,,\\
%%\ZLA{Eq:dynaPROPoprimalSTATE}
%%&
 \Xi^+_{t_0}(t;\Sta^+_{\zt}(t_0;\Sta_\zt))=\Xi^+_\zt(t;\Sta_\zt)  \,.
 \end{array}\right.
 \end{align}
 \end{subequations}
 \end{Theorem}
The standard proof   follows from the \emph{Bellman principle,} i.e. it follows by applying the inequality 
\begin{multline*}
W_\zt(\Sta_\zt)=W_\zt(\Xi_\zt,  \hat y_\zt)\\
=\int_\zt^{t_0} [\| v^+_\zt(s;\Sta_\zt )\|^2+\|u^+_\zt(s;\Sta_\zt)\|^2]\ZD s +
\int _{t_0}^T[\| v^+_\zt(s;\Sta_\zt )\|^2+\|u^+_\zt(s;\Sta_\zt)\|^2]\ZD s\\ 
\leq \int_\zt^{t_0} [\| v^+_\zt(s;\Sta_\zt )\|^2+\|u^+_\zt(s;\Sta_\zt)\|^2]\ZD s 
+
\int _{t_0}^T[\| v(s;t_0;\Sta^+_{\zt}(t_0;\Sta_\zt),u)\|^2+\|u(s)\|^2]\ZD s\,.
\end{multline*}
Unicity of the optimal control implies that
the inequality holds with the equal sign     if and only if $u $ is the restriction of 
$u^+_\zt(\cdot;\Sta_\zt)$ to $[t_0,T)$ so that~(\ref{Eq:dynaPROPoprimalCONTRP}) holds.
The equality~(\ref{Eq:dynaPROPozeroEVOL}) is obtained from~(\ref{Eq.secDynaPro:DYnaProsoluANTE}) (i.e. from the     unicity of   solutions of Eq.~(\ref{eq:laEQsu(ztT)maDAzt1})) and~(\ref{Eq:dynaPROPoprimalCONTRP}). 

Similarly,  Bellman principle   leads to the disssipation inequality. Let us fix a \emph{continuous} control $u$ on $[\zt,t]$ and any square integrable  control $u_1$ on $ (t,T)$. We have
\begin{multline*}
W_\zt(\Sta_\zt) \leq \int_\zt ^t [ \|v(s;\zt,\Sta_\zt,u)\|^2+\|u(s)\| ^2 ]\ZD s\\+
\int_t^T[ \|v(s;t,\Sta(t;\zt,\Sta_\zt,u),u_1)\|^2+\| u_1(s)\|^2]\ZD s\,.
\end{multline*}
We compute the minimum respect to  $u_1\in L^2(t,T;U)$. We get the following inequality, which is the \emph{dissipation inequality in integral form:}
\[
-\int_\zt ^t [ \|v(s;\zt;\Sta_\zt,u)\|^2+\|u(s)\| ^2 ]\ZD s\leq W_t( \Sta(t;\zt;\Sta_\zt ,u)  )-W_\zt(\Sta_\zt)\,.
\]
This inequality holds for every control $u$ which is continuous on $[\zt,t]$. Unicity of the optimal control implies  that the inequality is strict, unless $u$ is the restriction to $[\zt,t]$ of the optimal control.
 
Let $\ztn\in [\zt,t)$.
The same argument on $[\ztn,t]$ and with the initial condition $\Sta(\ztn;\zt, \Sta_\zt ,u)$ gives
 
\begin{multline}\ZLA{eq:PreDIdiffe}
 0\leq \int_{\ztn} ^t [ \|v(s;\ztn,\Sta(\ztn;\zt, \Sta_\zt ,u),u)\|^2+\|u(s)\| ^2 ]\ZD s\\
+W_t( \Sta(t;\zt;\Sta_\zt ,u) )-W_{\ztn}(\Sta(\ztn;\zt, \Sta_\zt ,u))\,.
\end{multline}
We divide both the sides with $t-\ztn>0$ and we pass to the limit for $t\to \ztn^+$.  We  get:
\begin{equation}
\ZLA{LOIdiffeFORM}
0\leq \left [\, \|v(\ztn;\zt,\Sta_\zt,u)\|^2+\|u(\ztn)\|^2\, \right ]+ D_+ W_{\ztn}(\Sta(\ztn;\zt, \Sta_\zt ,u))
\end{equation}
where $D_+$ denotes the first right Dini derivative computed along $\Sta(t)$:
 
 \begin{multline}\ZLA{Eq:defiDINIderiv}
D_+ W_{\ztn}(\Xi(\ztn;\zt;\Xi_\zt,\hat y,u),\hat y)
 \\
=\liminf _{t\to \ztn^+}\dfrac{W_t( \Sta(t;\zt;\Sta_\zt ,u) )-W_{\ztn}(\Sta(\ztn;\zt, \Sta_\zt ,u))}{t-\ztn}\,.
 \end{multline}

Inequality~(\ref{LOIdiffeFORM}) is the \emph{dissipation inequality in differential form.} Unicity of the optimal control implies that the inequality is strict, unless $u$ is the optimal control. So, for every $\ztn\in [\zt,T]$
(in particular also when $\ztn=\zt$) we have
\begin{equation}
\ZLA{eq:conseLOIdiffeFORM}
u^+_{\zt }(\ztn;\Sta_\zt)=\argmin\left \{
\|v(\ztn;\zt,\Sta_\zt,u)\|^2+\|u(\ztn)\|^2 + D_+ W_{\ztn}(\Sta(\ztn;\zt, \Sta_\zt ,u))
\right \}\,.
\end{equation}
 \emph{This equality shows that the dissipation inequality reduce the optimization problem in the space $L^2(\zt,T;U)$ to a family
 parametrized by $\ztn\in[\zt,T]$ 
 of optimization problems in $U$.}
 
 Alternatively, when $t>\zt$, we can divide both the sides of~(\ref{eq:PreDIdiffe}) with $\zthe-t<0$ and pass to the limit for $\zthe\to t^-$. We get the following alternative form of~(\ref{LOIdiffeFORM})
 \begin{align}
 \ZLA{LOIdiffeFORMvariant}
 &D_- W_t(\Sta(t;\zt,S_\zt,u)+\|v(t;\zt,\Sta_\zt,u)\|^2+\|u(t)\|^2\geq 0\\
\ZLA{Eq:defiDINIderivLEFT} &D_- W_t(\Sta(t;\zt,S_\zt,u)= \liminf _{s\to t^-}\dfrac{W_s(\Sta(s;\zt,\Sta_\zt,u)-W_t(\Sta(t;\zt,\Sta_\zt,u)}{s-t}
 \end{align} 
($D_- W_t(\Sta(t;\zt,S_\zt,u)$ is the first left Dini derivative computed along $\Sta(t)$).

The number $\ztn $ in~(\ref{LOIdiffeFORM}) and the number 
$t$ in~(\ref{LOIdiffeFORMvariant}) are arbitrary numbers in $[\zt,T]$. So, we can use the same notation for both of them, for example~(\ref{LOIdiffeFORMvariant}) holds as well with $t$ replaced by $\ztn $.

A result proved in the Appendix~\ref{AppeDeriValue} is that
when $\Sta_\zt\in\Dom\,\mathcal{A}_\zt$ and $u$ is continuous
 both the right and left derivatives in~(\ref{LOIdiffeFORM}) and~(\ref{LOIdiffeFORMvariant}) are the ordinary directional derivatives (i.e. the $\liminf$ is actually the limit) and furthermore they have the same value so that the value function is differentiable at any $\ztn >\zt$ and the  right derivative at $\zt$ exists. More precisely we prove:
 \begin{Theorem}\ZLA{TeoDERIVveraVALUEfunction}
Let $\zt\in [0,T)$,   $\Sta_\zt\in\Dom\,\mathcal{A}_\zt$ and let $u$ be continuous. Then  
 \begin{equation}\ZLA{Eq:ValueFUNcI}
\ztn  \mapsto W_{\ztn}(\Sta(\ztn;\zt, \Sta_\zt ,u)) \quad\mbox{belongs to $C^1([\zt,T])$}\,.
 \end{equation}
Both the inequalities~(\ref{LOIdiffeFORM}) and~(\ref{LOIdiffeFORMvariant}) take   the following forms    on $[\zt ,T]$:
 
\begin{equation}
 \ZLA{LOIdiffeFORMvera}
  0\leq \left [\, \|v(\ztn;\zt,\Sta_\zt,u)\|^2+\|u(\ztn)\|^2\, \right ]+ \dfrac{\ZD}{\ZD\ztn} W_{\ztn}(\Sta(\ztn;\zt, \Sta_\zt ,u)) \\
 \end{equation} 
 (directional derivatives when $\ztn=\zt$ or $\ztn=T$).
 \end{Theorem}
Also the inequality~(\ref{LOIdiffeFORMvera}) is called the  \emph{dissipation inequality in differential form.} 

The goal now  is the  proof of Theorem~\ref{TeoDERIVveraVALUEfunction} and the explicit computation of the derivative so that   the equality~(\ref{LOIdiffeFORMvera}) can be effectively used to compute the optimal control.
\begin{Remark}\ZLA{RemaDELTeoDERIVveraVALUEfunction}{\rm The function $\ztn\mapsto W_{\ztn}(\Sta(\ztn;\zt, \Sta_\zt ,u))$ is continuous (Theorem~\ref{teo:SullaVALUEfunctEsueCOMPON}). So, the proof that~(\ref{Eq:defiDINIderiv}) is the standard right derivative and that the right derivative is continuous implies the property in~(\ref{Eq:ValueFUNcI}).   The alternative form~(\ref{LOIdiffeFORMvariant}) of the dissipation inequality and the explicit study of the left derivative are not needed to prove Theorem~\ref{TeoDERIVveraVALUEfunction}. 
We need to consider them explicitly for a reason we shall see in Sect.~\ref{OPtiCONTeRICCATI}.\zdia

}\end{Remark}
 
 \subsection{\ZLA{SectValueFunDERIV}The value function and its derivative}
  We recall the definition
 in~(\ref{DefiPROIEXdaSTATO})-(\ref{COSTOottimoCOMEfarmaQUADRATICA}). 
 We fix  any initial time   $\zt_0\in[0,T)$
and any   $\ztn\in [\zt_0,T)$. Our goal in this section is the computation of the Dini derivatives in~(\ref{Eq:defiDINIderiv}) and in~(\ref{Eq:defiDINIderivLEFT})   (this one when $\ztn>\zt_0$) at $\ztn$ along $\Sta(t)=\Sta(t;\Sta_{\zt_0},u)$.

    Equalities~(\ref{ValueFunction}) and~(\ref{COSTOottimoCOMEfarmaQUADRATICA}) give 
 
 \begin{equation}\ZLA{Eq:defiValueCOMEquadFORM}
W_\ztn(\Sta(\ztn))=\int_{\ztn}^T\ZL  \left [ H_{\ztn} h_{\ztn}  \right ](s),  h_{\ztn}    (s)\ZR\ZD s=\ZL P(\ztn)\Sta(\ztn),\Sta(\ztn)\ZR 
 \end{equation}
 where $P(\ztn)\in \mathcal{L}(M^2_{\ztn}\times (\Dom\,A)')$
 is defined by
%   \begin{subequations}
 \begin{multline}
 \ZLA{EQdefiRiccaOpe}
 \ZL P(\ztn)\Sta_1(\ztn),\Sta_2(\ztn)\ZR_{M^2_{\ztn}\times (\Dom\,A)'}=
 \int_{\ztn}^T\ZL  \left [ H_{\ztn} h_{1;\ztn}  \right ](s),  h_{2;\ztn}    (s)\ZR_H\ZD s \\
 =
 \Bigl\ZL H_\ztn 
      (\Gamma_\ztn  \mathcal{C}^{(1,2)}_\ztn+L_\ztn \mathcal{C}^{(3)}_\ztn  )\Sta_{1,\ztn}, (\Gamma_\ztn  \mathcal{C}^{(1,2)}_\ztn+L_\ztn \mathcal{C}^{(3)}_\ztn )\Sta_{2,\ztn}\Bigr\ZR _{L^2(\ztn,T;H)} \,.
 \end{multline}
Here $\Sta_{1,\ztn}$ and $\Sta_{2,\ztn}$ are two arbitrary elements of the state space $M^2_\ztn\times (\Dom\,A)'$.
 
 The selfadjoint operator $P(\ztn)$ is the \emph{Riccati operator} at time $\ztn\in[\zt_0,T]$ of the quadratic regulator problem.

We prove:
\begin{Theorem}
 \ZLA{teo:SullaVALUEfunctEsueCOMPON} 
The  real valued  function
$\ztn\mapsto W_{\ztn}(\Sta(\ztn))$
 is continuous on $[\zt_0,T)$ and $W_T(\Sta(T))=\lim _{\ztn\to T^-}W_{t}(\Sta(t))=0$. 
\end{Theorem}
\zProof 
The existence of the limit and the  equality $W_T(\Sta(T))=0$ are obvious properties, consequences of the boundedness of the function
$
\Sta(t) 
$.
We prove right continuity at $\ztn\in[\zt_0,T)$.

We have
\begin{multline*}%\ZLA{Eq:PERlaCONTIdiWcomeFUNrho}
W_{\ztn+h}(\Sta(\ztn+h))- W_{\ztn }(\Sta(\ztn ))=
-\int_{\ztn}^{\ztn+h} \ZL  \left [ H_{\ztn} h_{\ztn}  \right ](s),  h_{\ztn}    (s)\ZR\ZD s\\
+
 \int_{\ztn+h}^T\left[\ZL  \left [ H_{\ztn+h} h_{\ztn+h}  \right ](s),  h_{\ztn+h}    (s)\ZR -
\ZL  \left [ H_{\ztn} h_{\ztn}  \right ](s),  h_{\ztn}    (s)\ZR\right ]\ZD s
\end{multline*}
 The definition of $h_{\ztn}\ZCD$ and Theorem~\ref{teoAPPEregolarPRPsoluFred} in Appendix~\ref{AppePROpreHzt} show that the   integrand of the first  integral on the right side is continuous and  bounded. So, the first integral on the right hand side converges to $0$ when $h\to 0^+$. In order to prove that the second integral tends to $0$ too, we use Theorem~\ref{teoAPPEregolarPRPDIpedaTAU} in Appendix~\ref{AppePROpreHzt}  with $\zt=\ztn$ and  $g(\cdot;\ztn)=h(\cdot;\ztn)$
so that
 \[
\lim _{h\to 0^+}g(s;\ztn+h)= \lim _{h\to 0^+}  h_{\ztn+h}    (s) = h_{\ztn }    (s)
 \] 
 poitwise and boundedly and also the assumption~\ref{I3teoAPPEregolarPRPDIpedaTAU} of Theorem~\ref{teoAPPEregolarPRPDIpedaTAU} is satisfied so that 
 \[
 \lim _{h\to 0^+}    \left [ H_{\ztn+h} h_{\ztn+h}  \right ](s) =
  \left [ H_{\ztn } h_{\ztn }  \right ](s) 
 \]
 pointwise and boundedly. So, the second integral in the right side converges to $0$ too.

Left continuity at $ \ztn>\zt_0 $ is proved analogously.\zdia

Now we study  the Dini derivatives in~(\ref{Eq:defiDINIderiv}) and in~(\ref{Eq:defiDINIderivLEFT} ) along $\Sta(t ;,\zt_0,\Sta_{\zt_0},u )$.
As it is well known for distributed systems, the derivative of the value functions exists when the initial condition is smooth. So:
\begin{Assumption}{\rm
In order to compute the derivative  we assume
\begin{equation}\ZLA{AssumptionPERderivat}
 \left\{\begin{array}{l}
u\in C([\zt_0,T];U)\\
\Sta_{\zt_0}\in \Dom\,\mathcal{A}_{\zt_0}\,. 
\end{array}\right.
\end{equation}
 
}
\end{Assumption}

First we consider  the right derivative in~(\ref{Eq:defiDINIderiv}). So we study the right limit for $h\to 0^+$ of the incremental quotient
\begin{multline}\ZLA{INCREquiotieDIvalueFUNC}
\dfrac{1}{h}\left [\ZL P(\ztn+h)\Sta(\ztn+h),\Sta(\ztn+h)\ZR-
\ZL P(\ztn )\Sta(\ztn ),\Sta(\ztn )\ZR\right ]\\
=\dfrac{1}{h}\left \{
\int _{\ztn+h}^T [ \ZL H_{\ztn+h}h_{\ztn+h} ](s), h_{\ztn+h}  (s)\ZR\ZD s
-
\int _{\ztn }^T \ZL [H_{\ztn }h_{\ztn } ](s), h_{\ztn }  (s)\ZR\ZD s 
\right \}\,.
\end{multline}

As usual, by adding and subtracting suitable terms, we reduce this incremental quotient to the sum of simpler incremental quotients but, due to to the presence of the memory, we need a bit of care. In fact:
\begin{itemize}
\item we should introduce the difference $\Gamma _{\ztn+h}\Xi(\ztn+h)-\Gamma_\ztn\Xi(\ztn+h)$ but the term $\Gamma_\ztn\Xi(\ztn+h)$ is meaningless because $\Gamma_\ztn$ acts on $M^2_{\ztn}$ while $ \Xi(\ztn+h)\in M^2_{\ztn+h}$, a different space. We introduce the operators $ \Mh\in\mathcal{L}(M^2_{\ztn+h}, M^2_\ztn) $ defined as follows:
 
\[
\Mh\Xi(\ztn+h)=
\Mh(\hat v,  v(\ztn+h-s)_{s\in (0,\ztn+h)})=(   \hat v,                v(\ztn+h-s)_{s\in (0,\ztn )}) \,.
\]
 Then we add and subtract $\Rh\Gamma_\ztn\Mh\Xi(\ztn+h)$ where
\[
 \Rh\in\mathcal{L}(L^2(\ztn,T;H),L^2(\ztn+h,T;H))\ \mbox{ is the restriction}\,.
\] 
This way we have the meaningfull difference
\[
\Gamma _{\ztn+h}\Xi(\ztn+h)-\Rh\Gamma_\ztn\Mh\Xi(\ztn+h)\in L^2(\ztn+h,T;H)\,.
\]
We can simply write $\Gamma _{\ztn+h}\Xi(\ztn+h)- \Gamma_\ztn\Mh\Xi(\ztn+h)$
when the restriction is forced by  an integration on $[\ztn+h,T]$.

\emph{Note that if $\Xi_{\ztn+h}\in \Dom\, \mathcal{A}_{\ztn+h}$ then $\Mh \Xi_{\ztn+h}\in\Dom\, \mathcal{A}_\ztn$.}
 
\item we need also the difference $H_{\ztn+h}h_{\ztn+h}-H_{\ztn } h_{\ztn}$
which is integrated on $[\ztn+h,T]$, so in fact the difference is
$H_{\ztn+h}h_{\ztn+h}-\Rh H_{\ztn } h_{\ztn}$.
In order to examine this difference we put 
 
\begin{equation}\ZLA{eqDIFEdIphistheta}
 \phi(s;\ztn)=[H_{\ztn} h_\ztn](s)\quad  \mbox{for $s\in[ \ztn, T]$}\,.
\end{equation} 
Then, for $s\in [\ztn+h,T]$ we have
\begin{multline*}
[H_{\ztn+h}h _{\ztn+h} ](s)-[H_{\ztn }h_ \ztn  ](s)\\
%%%%
=[H_{\ztn+h}h_{\ztn+h}](s)-[\Rh H_{\ztn } h_{\ztn}](s)
%%%
=
\left [ \, \phi(s,\ztn+h)-\phi(s,\ztn)\,\right ]\\
=\int _{\ztn}^{\ztn+h} K(s-r)\int_r^TK^*(\nu-r)\phi(\nu,\ztn)\ZD\nu\,\ZD r+\left\{   h_{\ztn+h}  (s)-\Rh h_{\ztn } (s)\right \}\\
-\int_{\ztn+h}^sK(s-r)\int_r^TK^*(\nu-r)\left [\,\phi(\nu,\ztn+h)-\phi(\nu,\ztn)\, \right ]\ZD\nu\,\ZD r
\,.
\end{multline*}
I.e. on $[\ztn+h,T]$ we have
\[%begin{multline*}
H_{\ztn+h}h _{\ztn+h}   -                      H_{\ztn }h_ \ztn 
 =H_{\ztn+h}[h_{\ztn+h} -\Rh h_{\ztn} ] 
+
H_{\ztn+h}\Phi(\cdot;\ztn,h) 
\]%end{multline*}
where 
\begin{equation}\ZLA{eqDefiPHImaiuscolo}
\Phi(s;\ztn,h)= 
 \int _{\ztn}^{\ztn+h}K(s-r)\int_r^TK^*(\nu-r)\phi(\nu,\ztn )\ZD\nu\,\ZD r \,.
\end{equation}
\end{itemize} 

We use the operators $\Mh$ and $R^{\ztn}_{\ztn+h}$ and the function  $\Phi$ just defined in order to elaborate the incremental quotient~(\ref{INCREquiotieDIvalueFUNC}).
We see that we must study the limit for $h\to0^+$ of the following sum of
eight addenda from~(\ref{eqRaIncre2}) to~(\ref{subequatioPERdopobmepsilon}). 
The addenda are grouped in such a way to stress the different roles they have in the computation of the derivative of the value function. 

For the sake of concision,
 we do not explicitly write the restriction operator $\Rh$ when the restriction is forced by the domain of integration   $[\ztn+h,T]$ and  we do not   indicate the variable of integration.

% \newpage

 \begin{subequations}
 \ZLA{eqINCRequoDOLOstato}
\begin{align}
\nonumber
&\dfrac{1}{h}\left \{\ZL P(\ztn+h)\Xi(\ztn+h),\Xi(\ztn+h)\ZR-
\ZL P(\ztn )\Xi(\ztn ),\Xi(\ztn )\ZR\right\}\\[3mm]
\ZLA{eqRaIncre2}
&=\int_{\ztn+h}^T\Bigl\ZL 
 H_\ztn h_{\ztn}, \Gamma_{\ztn}\dfrac{1}{h}\left [ \Mh\Xi(\ztn+h)- \Xi(\ztn)\right ]
\Bigr\ZR\ZD s
\\[2mm]
\ZLA{eqRaIncre3}
&
+\int_{\ztn+h}^T\Bigl\ZL 
 H_{\ztn+h}\Rh\Gamma_{\ztn} 
  \frac{1}{h}\left [   \Mh\Xi(\ztn+h)-  \Xi(\ztn)\right ]
,h_{\ztn+h}
 \Bigr\ZR\ZD s
 \end{align}
% \vskip .5cm 
\end{subequations}
 \vskip-.2cm
\begin{subequations}
\ZLA{eqINCRequoDOLOriccati}
\begin{align}
\ZLA{subequatioPERdopobmsigma}
\quad & -\dfrac{1}{h}\int_\ztn^{\ztn+h} \ZL 
H_\ztn h_\ztn ,h_\ztn
\ZR\ZD s \\[2mm]
 &  + 
\int_{\ztn+h}^T \Bigl \ZL  \dfrac{1}{h}
\Phi(\cdot;\ztn,h)
  , H_{\ztn+h}h_{\ztn+h}  \Bigr \ZR\ZD s \ZLA{subequatioPERdopobmzaa}
  \\[2mm]
\ZLA{subequatioPERdopobmbeta}
 &  +\int _{\ztn+h}^T \Bigl\ZL H_{\ztn}h_{\ztn} ,\dfrac{1}{h}\left [\Gamma_{\ztn+h}\Xi(\ztn+h) -\Gamma_{\ztn}\Mh\Xi(\ztn+h) \right ]\Bigr\ZR\ZD s 
 \\[2mm]
 .\ZLA{subequatioPERdopobmgamma}
&+\int _{\ztn+h}^T\Bigl\ZL H_{\ztn+h}  \dfrac{1}{h}\left [\Gamma_{\ztn+h}\Xi(\ztn+h)-\Rh\Gamma_\ztn \Mh\Xi(\ztn+h) \right ], h_{\ztn+h} \Bigr\ZR\ZD s 
\end{align}
\end{subequations}
   \vskip -.6cm
\begin{subequations}
\ZLA{DaSUBEQparteCOMUNE}
\begin{align}
\ZLA{subequatioPERdopobmdelta} 
& +\int_{\ztn+h}^T\Bigl\ZL 
  H_\ztn h_\ztn, \dfrac{1}{h}\left [ \hat Y_{\ztn+h}-\hat Y_{\ztn}  \right ]\Bigr\ZR\ZD s 
  \hskip 5cm
\\[2mm] 
\ZLA{subequatioPERdopobmepsilon}
 &
+\int _{\ztn+h}^T 
  \Bigl\ZL   
  H_{\ztn+h} \dfrac{1}{h}\left[ \hat Y_{\ztn+h} -\Rh \hat Y_{\ztn} \right ] ,h_{\ztn+h} 
  \Bigl\ZR \ZD s  \,.
\end{align}
\end{subequations}
\emph{ We must compute the limit of this expression for $h\to 0^+$.}

 The rationale behind the grouping of the addenda is that the increment $h$  in both the addenda~(\ref{eqINCRequoDOLOstato}) affects only the incremental quotient of $\Xi$, so of the state;  it affects the incremental quotient of the Riccati operator but not that of the   state
in the four addenda~(\ref{eqINCRequoDOLOriccati}). 
Instead,   the incremental quotient  of the state as well as that of the operator are affected in both the integrals~(\ref{DaSUBEQparteCOMUNE}).

When using the alternative form~(\ref{LOIdiffeFORMvariant}) of the dissipation inequality, the increment is negative, $h=-k<0$. In this case we have the following equality. In this equality the lines are indicated by   the same number as the corresponding line when $h>0$, followed by $'$. 
 
\[\begin{array}{l}
\hskip -3cm \dfrac{1}{-k}\left \{\ZL P(\ztn-k )\Xi(\ztn-k ),\Xi(\ztn-k )\ZR-\ZL P(\ztn )\Xi(\ztn ),\Xi(\ztn )\ZR 
\right\}\\[2mm]
\hskip -3cm\quad = \dfrac{1}{k}\left \{\ZL P(\ztn )\Xi(\ztn ),\Xi(\ztn )\ZR-
\ZL P(\ztn-k )\Xi(\ztn-k ),\Xi(\ztn-k )\ZR\right\}
\end{array}
\]
\[
\begin{array}{l}\qquad ={\displaystyle  \int_{\ztn }^T}\left [\Bigl\ZL 
 H_{\ztn-k} h_{\ztn-k}, \Gamma_{\ztn-k}\dfrac{1}{k}\left [ \Mk\Xi(\ztn )- \Xi(\ztn-k)\right ]
\Bigr\ZR \right.\\[1mm]
\qquad\quad  \left.+\Bigl\ZL 
 H_{\ztn }\Rk\Gamma_{\ztn-k} 
  \frac{1}{k}\left [   \Mk\Xi(\ztn )-  \Xi(\ztn-k)\right ]
,h_{\ztn }
 \Bigr\ZR\right ]\ZD s
 \end{array} \eqno{(\ref{eqRaIncre2}')+(\ref{eqRaIncre3}')}
 \]
% \vskip -2mm
\[
 \hskip .11cm\quad -\dfrac{1}{k}\int_{\ztn-k}^{\ztn }\Bigl \ZL 
H_{\ztn-k}  h_{\ztn-k} ,h_{\ztn-k}\Bigr\ZR\ZD s
  - 
\int_{\ztn }^T \Bigl \ZL  \dfrac{1}{k}
\Phi(\cdot;\ztn ,-k)
  , H_{\ztn }h_{\ztn }  \Bigr \ZR\ZD s \eqno{(\ref{subequatioPERdopobmsigma}')+(\ref{subequatioPERdopobmzaa}')}
  \]
  \vskip-.5cm
  \[ 
\begin{array}{l} +{\displaystyle \int _{\ztn }^T} \left [\Bigl\ZL H_{\ztn-k}h_{\ztn-k} ,\dfrac{1}{k}\left [\Gamma_{\ztn }\Xi(\ztn ) -\Rk\Gamma_{\ztn-k}\Mk\Xi(\ztn ) \right ]\Bigr\ZR \right.
 \\[1mm] 
\qquad \left. + \Bigl\ZL H_{\ztn }  \dfrac{1}{k}\left [\Gamma_{\ztn }\Xi(\ztn )-\Rk\Gamma_{\ztn-k} \Mk\Xi(\ztn) \right ], h_{\ztn} \Bigr\ZR\right]\ZD s\end{array} \eqno{(\ref{subequatioPERdopobmbeta}')+(\ref{subequatioPERdopobmgamma}')}
\]
 \vskip -.3cm
 \[
 +{\displaystyle\int_{\ztn }^T} \left [\Bigl\ZL 
  H_{\ztn-k} h_{\ztn-k}  , \dfrac{1}{k}\left [ \hat Y_{\ztn }- \hat Y_{\ztn-k}  \right ]\Bigr\ZR  
 + 
  \Bigl\ZL   
  H_{\ztn } \dfrac{1}{k}\left[ \hat Y_{\ztn } -\Rk \hat Y_{\ztn-k} \right ] ,h_{\ztn} 
  \Bigl\ZR \right]\ZD s  \,.
 \eqno{(\ref{subequatioPERdopobmdelta}')+(\ref{subequatioPERdopobmepsilon}')}\,.
 \]
\emph{ We must compute the limit of this expression for $k\to 0^+$.}

In conclusion, in both the cases we have the sum of eight addenda. We prove, in the Appendix~\ref{AppeDeriValue}, that 
under the assumption~(\ref{AssumptionPERderivat}) the limit of each one of them exists.
By collecting the results in the  Appendix~\ref{AppeDeriValue} we obtain
the proof of Theorem~\ref{TeoDERIVveraVALUEfunction} and
 the main conclusion of this paper presented in the next Sect.~\ref{OPtiCONTeRICCATI}: \emph{the optimal control can be represented as a feedback on the state and the feedback operator at the time $t$ is $-\mathcal{B}^*_tP(t)$. Furthermore, $P(t)$ solves a Riccati differential equation.} 
 
The results in the Appendix~\ref{AppeDeriValue}   give a detailed analysis of the differentiability of the Riccati operator.

   \section{\ZLA{OPtiCONTeRICCATI}The optimal control and the computation of the Riccati operator}
In this section we derive an explicit form of the optimal control as a linear  feedback of the state and, by suitably collecting the results in the  Appendix~(\ref{AppeDeriValue}), a Riccati type equation for the Riccati operator.

The right derivative of the value function is the limit of the eight addenda from~(\ref{eqRaIncre2}) to~(\ref{subequatioPERdopobmepsilon})
but it is convenient to split the limits of the two addenda in~(\ref{subequatioPERdopobmepsilon}) so to separate the contribution of the incremental quotient of the state and that of the operator $P(t)$. So we write
 
\begin{equation}\ZLA{formDERIvalFunctDIVIdueParti}
 \dfrac{\ZD}{\ZD\ztn}W_\ztn(\Sta(\ztn)) 
 =\left [
\mbox{(\ref{EqSommaDUEterminiINDaSUBEQparteCOMUNE})}+\mbox{(\ref{ValeLimiZTNeqRaIncre2})}+\mbox{(\ref{ValeLimiZTNeqRaIncre3})}\right]\\
 +\left[
\mbox{(\ref{EqLIMITEesplicitoDIsubequatioPERdopoDELTAsigma})}+
 \mbox{(\ref{EqLIMITEesplicitoDIsubequatioPERdopoALPHA})}+
 \mbox{(\ref{LasommaDIPARTEderiRiccaOper})}
\right]\,.
\end{equation}

This equality holds when $\Sta\ZCD$ is a solution of~(\ref{StateEQUAcompleta}) for $t\geq\zt_0$, provided that the conditions~(\ref{AssumptionPERderivat}) hold.

The right hand side is a continuous function of $\ztn\geq \zt_0$ and this observation is sufficient to complete the proof of Theorem~\ref{TeoDERIVveraVALUEfunction}, as we noted in Remark~\ref{RemaDELTeoDERIVveraVALUEfunction}.

 We use the definition~(\ref{EQdefiRiccaOpe}) of the Riccati operator and we see
\begin{multline}
\ZLA{CollegoINCREstato}
\left [
\mbox{(\ref{EqSommaDUEterminiINDaSUBEQparteCOMUNE})}+\mbox{(\ref{ValeLimiZTNeqRaIncre2})}+\mbox{(\ref{ValeLimiZTNeqRaIncre3})}\right]=\int_\ztn ^T\Bigl\ZL H_\ztn h_\ztn,[\Gamma_\ztn\mathcal{C}^{(1,2)}_\ztn+L\mathcal{C}^{(3)}_\ztn]\Sta'(\ztn)\Bigr\ZR \ZD s\\
+
\int_\ztn ^T\Bigl\ZL
[\Gamma_\ztn\mathcal{C}^{(1,2)}_\ztn+L\mathcal{C}^{(3)}_\ztn]\Sta'(\ztn),
  H_\ztn h_\zt\Bigr\ZR\ZD s\\
=\ZL P(\ztn)\Sta(\ztn),\Sta'(\ztn)\ZR +\ZL \Sta'(\ztn),P(\ztn)\Sta(\ztn)\ZR\,.
\end{multline} 
\begin{Remark}\ZLA{RemaREGOLARdiP}{\rm 
We note:
\begin{itemize}
\item the vector $\Sta'(\ztn)$ needs not belong to the state space  $M^2_{\ztn}\times (\Dom\,A)'$. In spite of this, the right hand side makes sense thanks to the smoothing properties of the Riccati operator proved in the sections~\ref{ParaDaSUBEQparteCOMUNE} and~\ref{paraeqINCRequoDOLOstato} (see Remark~\ref{remaREGULpropRIccati}).
\item
The    Hilbert spaces in this paper are real,   so that the two addenda have the same value, but we prefer to keep this more symmetric form.\zdia
\end{itemize}
}
\end{Remark}
The second addendum 
 in the right side of~(\ref{formDERIvalFunctDIVIdueParti}) 
depends on the incremental quotient of the operator $P(\ztn)$ but not on that of $\Sta(\ztn)$. So we put
 
\begin{equation}
\ZLA{eqSHORTperDERIpRICCA}
\mbox{(\ref{EqLIMITEesplicitoDIsubequatioPERdopoDELTAsigma})}+
 \mbox{(\ref{EqLIMITEesplicitoDIsubequatioPERdopoALPHA})}+
 \mbox{(\ref{LasommaDIPARTEderiRiccaOper})}
= 
  \ZL  P'(\ztn)\Sta_\ztn ,\Sta_\ztn  \ZR \,.
  \end{equation}
 
  The explicit definition, obtained by summing the addenda on the left side, is
  \begin{multline}\ZLA{eqSHORTperDERIpRICCAexplicitFORM}
   \ZL  P'(\ztn)\Sta_\ztn ,\Sta_\ztn  \ZR 
  =\|\mathcal{C}^{(1)}\Sta_\ztn\|^2_H
+\| [\Lambda_\ztn^*H_\ztn h_\ztn](\ztn)\|^2_H
 % \\
   -\ZL H_\ztn h_\ztn, \Gamma_\ztn \mathcal{C}^{(1,2)}\mathcal{A}_\ztn\Sta_\ztn\ZR_{L^2(\ztn,T;H)}\\
%%%   
  -\ZL \Gamma_\ztn \mathcal{C}^{(1,2)}\mathcal{A}_\ztn\Sta_\ztn,
  H_\ztn h_\ztn\ZR_{L^2(\ztn,T;H)}\qquad \forall\; \Sta_\ztn\in\Dom\,\mathcal{A}_\ztn\,.
\end{multline} 
  The left hand side of~(\ref{eqSHORTperDERIpRICCA}) is computed as a right limit, the limit for $h\to 0^+$ of the addenda in~(\ref{eqINCRequoDOLOriccati}) 
and of those in~(\ref{DaSUBEQparteCOMUNE}) which correspond to the incremental quotient of the Riccati operator.  
 \emph{The same expression  is obtained when computing the left limits, i.e. the sum of the limits for $k\to 0^+$ of the corresponding addenda   and this justifies the notation $P'(\ztn)$.}

% This way we define the derivative of the Riccati operator in a weak sense. \emph{The weak sense in which $P'$ is defined is the computation of the  (right or left) limits which concur to the equality~(\ref{eqSHORTperDERIpRICCAexplicitFORM}).
%  }

The computation of the limits in~(\ref{CollegoINCREstato}) use the conditions~(\ref{AssumptionPERderivat}). It is important to note that the full force of this assumption is not needed to compute $P'$. In  ordet to compute $P'$ it is sufficient to use  a function   $\Sta(t)$  constructed as described in the following assumption.
 \emph{We stress the fact that $\Sta(t)$ as defined in the next assumption needs not solve Eq.~(\ref{StateEQUAcompleta}) i.e. it needs not be the evolution of the system driven by a    control.}

\begin{Assumption}\ZLA{AssuSuallaConveeqPERdopoBETAeDopoGAMMA}
  We assume:
\begin{enumerate}

\item\ZLA{I0AssuSuallaConveeqPERdopoBETAeDopoGAMMA} $y\ZCD\in C^1([0,T];H)$. 

\item\ZLA{I10AssuSuallaConveeqPERdopoBETAeDopoGAMMA}  let  $\ztn\in [0,T)$. 
         \begin{enumerate}  
		\item We fix a function $v\ZCD\in C([0,T];H)$ defined on $ [0,T]$ with the following properties:
			\begin{enumerate}
% \item\ZLA{I1AssuSuallaConveeqPERdopoBETAeDopoGAMMA}   $v\ZCD\in %C([0,T];H)$;
%  \begin{enumerate}
			\item\ZLA{I2AssuSuallaConveeqPERdopoBETAeDopoGAMMA}  
$v_{|_{[0,\ztn]}}\in H^1(0,\ztn;H)$ 
			\item\ZLA{I2bisAssuSuallaConveeqPERdopoBETAeDopoGAMMA}  
 $v_{|_{[\ztn,T]}} \in H^1(\ztn,T ;(\Dom(-A)^{\ZSI+\ZEP})') $.
			 \end{enumerate}
%		\item\ZLA{I3AssuSuallaConveeqPERdopoBETAeDopoGAMMA} For every $t\in[0,T]$ and $\nu\in [0,t]$ we define $\xi_t(\nu )=v(t-\nu)$.
		\item\ZLA{I4AssuSuallaConveeqPERdopoBETAeDopoGAMMA} for   $t\in (\ztn,T]$ we define  
\[
\Xi_t=(  v(t),v(t-\cdot)\in M^2_t\,,\qquad 
\Sta(t)=\left (   \Xi_{t} , y(t )\right )  \,.\zdiaform
\] 
		\end{enumerate}

 \end{enumerate}
\end{Assumption}
\emph{We note that $\Sta(\ztn)\in \Dom\,  \mathcal{A}_{\ztn}$.}

A function $\Sta\ZCD$ with the property in Assumption~\ref{AssuSuallaConveeqPERdopoBETAeDopoGAMMA} needs not solve
the state equation~(\ref{StateEQUAcompleta}) and it does not depend on a control but the  properties
 in Assumption~\ref{AssuSuallaConveeqPERdopoBETAeDopoGAMMA} are satisfied in particular when $\Sta(t)$ solves the state equation~(\ref{StateEQUAcompleta})   provided that the 
 Assumptions~(\ref{AssumptionPERderivat}) hold (see the statement~\ref{I0teoREgoDiv} of Theorem~\ref{teoREgoDiv}). In this case $\Sta(t)$ is driven by a continuous control but   \emph{also in this case, the right hand side of~(\ref{eqSHORTperDERIpRICCA}) does not depend on $u(t)$, $t\geq\ztn$.}  
The independence on the control can be seen by looking at the explicit expression~(\ref{eqSHORTperDERIpRICCAexplicitFORM})\footnote{\ZLA{FOOTnoIndePprimoU}More directly    since the incremental quotient of the state is not used in the definition of the functions which give the limits on the right side of~(\ref{eqSHORTperDERIpRICCA}).}:
\emph{the right side of~(\ref{eqSHORTperDERIpRICCA}) does not depend on the control $u $    not even when $\Sta\ZCD $ is the state of the system driven by~$u$. } 
 
 This way we define the derivative of the Riccati operator 
along a function $\Sta(t)$  with the properties in the Assumption~\ref{AssuSuallaConveeqPERdopoBETAeDopoGAMMA}. The derivative is defined in a weak sense.
 \emph{The weak sense in which $P'$ is defined is the computation of the  (right or left) limits which concur to the equality~(\ref{eqSHORTperDERIpRICCAexplicitFORM}). 
  }

 \subsection{\ZLA{OPtiCONTeRICCATIthecontrol}The feedback form of the optimal control} 
We consider an \emph{arbitrary} initial time in $\zt_0\in  [0,T)$  and the corresponding 
  optimal control  
 $u^+_{\zt_0}=u^+_{\zt_0}(\cdot;\Sta_\zt)$.

We recall:
\begin{enumerate}
\item the optimal control is continuous
and depends continuously on $\Sta_{\zt_0}$ (see Theorem~\ref{TEOValue functionCONTI});  
\item for every $t\geq {\zt_0}$ (in particular also for $t={\zt_0}$!)
$u^+_{\zt_0}(t)$  minimizes the dissipation inequality in differential  form~(\ref{LOIdiffeFORM}).  The value of the minimum is zero, a fact   used in the derivation of the Riccati equation but not in the derivation of the feedback form of the optimal control.
\item we rewrite formula~(\ref{Eq:dynaPROPoprimalCONTRP}):
 \[
u_{\zt_0}^+( t;\Sta_{\zt_0}) 
=u^+_t(t;\Sta^+_{\zt_0}(t;\Sta_{\zt_0}))
\,.
\] 
This equality holds for every $t\geq  \zt_0$.  It shows that $u_{\zt_0}^+( t;\Sta_{\zt_0})$ depends solely on the state $\Sta^+_{\zt_0}(t;\Sta_{\zt_0})$  reached by the optimal evolution of the system at the time $t$, i.e. \emph{the optimal control is a fedback of the state.}  The goal here is the derivation of an \emph{explicit} expression for feedback form of the optimal control.
 
\end{enumerate}
 
We insert~(\ref{formDERIvalFunctDIVIdueParti}) in~(\ref{eq:conseLOIdiffeFORM})  and we use
 ~(\ref{CollegoINCREstato}) and~(\ref{eqSHORTperDERIpRICCA}).
 We find:
\begin{multline}\ZLA{eq:LaFAMIprobleOPTIM}
u_{\zt_0}^+(t)=   \argmin _{u\in U}\Bigl \{\|\mathcal{C}^1_t \Sta^+_{\zt_0} (t)\|^2+\|u \|^2+
\ZL P(t)\Sta^+_{\zt_0}(t), 
 \mathcal{A}_{t}\Sta^+_{\zt_0}(t)+ \mathcal{B}_{t}u  
 \ZR\\
\left. +
\ZL 
 \mathcal{A}_{t}\Sta^+_{\zt_0}(t)+ \mathcal{B}_{t}u 
 ,P(t)\Sta^+_{\zt_0}(t)\ZR   +
 \ZL P' (t)\Sta(t),\Sta(t)\ZR
 \right\}\,.
\end{multline}

The explicit expression~(\ref{eqSHORTperDERIpRICCAexplicitFORM}) of $ \ZL P' (t)\Sta(t),\Sta(t)\ZR$ is not used but we use   the fact that  it  does not depend on $u$. We use this fact to compute the minimum and we find the following equality for every $t\in [\zt_0,T]$: $u^+_{\zt_0}(t)   =- {\mathcal{B}} _{t}^*P(t)\Sta^+_{\zt_0}(t)$. 
     When $t=\zt_0$
  we get    \begin{equation}
   \ZLA{feedbackFORMcontrol}
   u^+_t(t)   =- {\mathcal{B}} _{t}^*P(t)\Sta ( t )  \,.
   \end{equation}
In this expression,  $\Sta ( t )=\Sta_t$ is the initial condition at time $t$, So, it is an \emph{arbitrary} element of $\Dom\,\mathcal{A}_t$.
In fact, the   computations which lead to~(\ref{feedbackFORMcontrol}) make sense when $\Sta_{t}\in\Dom\,\mathcal{A}_{t }$ but we 
proved in Theorem~\ref{TEOValue functionCONTI}   that the optimal control depends continuously on $\Sta_t$. So, the operator ${\mathcal{B}} _{t}^*P(t)$ is bounded and equality~(\ref{feedbackFORMcontrol}) can be extended by continuity to the entire state space $M^2_t\times (\Dom\,\mathcal{A}_t)'$.

In conclusion, the feedback form of the optimal control is~(\ref{feedbackFORMcontrol}), an equality which hold for \emph{ every} $\Sta_t\in M^2_t\times (\Dom\,A)'$.
\begin{Remark}{\rm
It may have an interest to note that the existence of the derivative of $P(t)$ is not really used in the computation of the feedback form of the optimal control. We used the fact that the \emph{right} derivative exists and that it does not depend on $u$, a fact that can be obtained without any explicit computation from the observation in the footnote~\ref{FOOTnoIndePprimoU}.\zdia
}\end{Remark}

   \subsection{\ZLA{OPtiCONTeRICCATItheRICCATI}The determination of the Riccati operator and the Riccati equation}
   
In this section we derive a differential equation of Riccati type  for the Riccati operator $P(t)$.    We use explicitly  that the value of the minimum in~(\ref{eq:LaFAMIprobleOPTIM}) is equal to zero and we consider the equality at $t=\zt_0$. So, $\Sta _{\zt_0}(\zt_0) =\Sta_t$ is an \emph{arbitrary} element of $\Dom\,\mathcal{A}_t$. 
    
   We replace~(\ref{feedbackFORMcontrol}) in~(\ref{eq:conseLOIdiffeFORM}).  
We get the following equality: 
  
   \begin{multline} \ZLA{eqRiccati}
  \ZL P'(t)\Sta_t,\Sta_t\ZR+\ZL P(t)\Sta_t,\mathcal{A}_t\Sta_t\ZR+\ZL
   \mathcal{A}_t\Sta_t,P(t)\Sta_t\ZR\\
   -\ZL \mathcal{B}_t^*P(t)\Sta_t, 
    \mathcal{B}_t^*P(t)\Sta_t\ZR+\ZL \mathcal{C}^1_t\Sta_t,\mathcal{C}^1_t\Sta_t\ZR=0\,,\qquad P(T)=0\,.
   \end{multline}
  The equality holds for every $\Sta_t\in\Dom\,\mathcal{A}_t$ and, as noted in the Remark~\ref{RemaREGOLARdiP}, the inner product is that of the state space $M^2_t\times (\Dom\,A)'$.

 The equation~(\ref{eqRiccati})   is the Riccati differential equation for the control system~(\ref{eq:laEQsu(ztT)}).
  
\section{\ZLA{Sect:COMPARISON}Comparison with previous results and final comments}

We recapitulate the main idea  used in this paper: the   wave equation  
 with high internal damping~(\ref{eq:IntroNOTAZ}) is a second order differential equation. It has been transformed 
 to a heat equation with memory (of the type studied in~\cite[Ch.~3]{PandolfiLIBRO21}) by a suitable integration. This equation is written as a state space system in a \emph{time variable} state space (which is an extended 
 Hilbert space as defined in~\cite{DolezalBOOKmonotoneExtendedHspaces}, see the discussion in~(\cite{Pandolfi24EECTtrackingRd})). After that we used Bellman principle and the dissipation inequality to express the unique optimal control in feedback form and to derive a suitable version of the Riccati equation.
 
 The   quadratic regulator problem~(\ref{eq:IntroNOTAZ})-(\ref{eq:costoWAVEdamp}) has been previously studied using a different idea in the papers~\cite{LasiePANtrigHighlyDAMPEDbis,LasiePANtrigHighlyDAMPED,TriggianiHIGHLYDIEdaSOLO}. The idea in these papers is to impose the boundary conditions
 \begin{equation}\ZLA{bouCONDIPLT}
\zg_0(v(t)-Du(t))=0  \,, \quad \zg_0(v'(t)-Du'(t))=0\,.
 \end{equation}
 of course, this is possible if $u$ is differentiable, a condition then removed as we are going to outline.
 Once the boundary conditions~(\ref{bouCONDIPLT}) has been imposed, the equation takes the form
 
\begin{align*}
&W'=\mathcal{A} W-\mathcal{A} \mathcal{D}(u(t)+u'(t))\,,\qquad 
W =\left[\begin{array}{c }
v\\w
\end{array}\right]\,,\\
&\mathcal{A}=\left[\begin{array}{cc}
0&I\\
A&-A
\end{array}\right]\,,\qquad 
\mathcal{D}=\left[\begin{array}{c}
0\\D
\end{array}\right]
\end{align*}
so that
\begin{equation}\ZLA{eq:diWcomeinLTP}
W(t)=e^{\mathcal{A}t}W(0)-\intt e^{\mathcal{A}(t-s)}\mathcal{D}u(s)\ZD s
%\\
-\intt e^{\mathcal{A}(t-s)}\mathcal{D}u'(s)\ZD s \,.
\end{equation}
This equality makes sense when $u\in H^1(0,T;U)$ and it was previously used in~\cite{BucciTESI}  to study the quadratic regulator problem when the cost functional penalizes not only the control $u$ but also its derivative. Instead, the authors of~\cite{LasiePANtrigHighlyDAMPEDbis,LasiePANtrigHighlyDAMPED,TriggianiHIGHLYDIEdaSOLO} follow  an idea   introduced in~\cite{PandAMODelaySINGULAR} in the study of the quadratic regulator problem for systems with delays in $\zzr^n$: the last integral in~(\ref{eq:diWcomeinLTP}) 
 is integrated by parts:
\begin{equation}\ZLA{eq:INTEHIHinteDOPOintegr}
 \intt e^{\mathcal{A}(t-s)}\mathcal{D}u'(s)\ZD s
 =\mathcal{D}u(t)-
 e^{\mathcal{A}t}\mathcal{D}u(0)-\intt e^{\mathcal{A}(t-s)}[\mathcal{A}\mathcal{D}]u (s)\ZD s\,. 
\end{equation}
The high internal damping in the equation~(\ref{eq:IntroNOTAZ}) is reflected in the high degree of unboundedness of the operator $\mathcal{A}\mathcal{D}$ in the right side of~(\ref{eq:INTEHIHinteDOPOintegr}).

Eq.~(\ref{eq:diWcomeinLTP}) is then rewritten as

%%%%%%%%%%%%%%%%%%%
\begin{multline}\ZLA{eq:diWcomeinLTPriscritta}
W(t)=e^{\mathcal{A}t}\left [W(0)+\mathcal{D}u_0\right ]-\intt e^{\mathcal{A}(t-s)}\mathcal{D}u(s)\ZD s
 \\
-
\left [
\mathcal{D}u(t) 
-\intt e^{\mathcal{A}(t-s)}[\mathcal{A}\mathcal{D}]u (s)\ZD s
\right ]
\end{multline}
where $u_0$ is now a free parameter, not equal to $u(0)$ since $u$ in~(\ref{eq:diWcomeinLTPriscritta})  is assumed solely square integrable.

The quadratic regulator problem is associated to~(\ref{eq:diWcomeinLTPriscritta}). This way we get a family
$J_{u_0} (W(0),u\ZCD ))$
 of quadratic regulator problems parametrized by $u_0$. A delicate analysis of the  system (in particular, of the last integral in~(\ref{eq:diWcomeinLTPriscritta})) shows that for every $u_0$ the unique optimal control $u_{u_0}^+(\cdot;W(0) )$ is continuous and it can be written in feedback form thanks to an operator which solves a suitable version of the Riccati differential equation. In general  $u_{u_0}^+(0;W(0) )\neq u_0$.
%%%%%%%%%%%%%%%%%%%

 The final step   is the study of the minimum of $u_0\mapsto J (W(0),u_{u_0}^+(\cdot;W(0) ))$. It is proved that the minimum    point $\tilde u_0$ exists and that the corresponding control satisfies
 $u_{\tilde u_0}^+(0;W(0) )=\tilde u_0$.
 
 In conclusion, we have now two alternative methods to study the quadratic regulator problem~(\ref{eq:IntroNOTAZ})-(\ref{eq:costoWAVEdamp}). Both of them deserve  to be known.
 
 \subsection{\ZLA{subsEXTENS}Extensions}

We mentioned that our approach has been motivated by that used in~\cite{BucciPandolfiArxive,DelloroPATA17MGTeVISCO} to study the properties of the solutions of the MGT equation. We conjecture that the approach we presented here can be adapted to study the quadratic regulator problem for the controlled  MGT equation (see for example~\cite{BucciLasieRegulatAcouPress,ClasonKalteOPTIM,LAStri2022OPTIMmgt} for previous results on this subject) and also that of different classes of damped equations, as the plate equation with internal damping which has been previously studied in~\cite{LasieLukesPAND2,LasieLukesPAND1,TriggianiTRENTO} with the approach in
 ~\cite{LasiePANtrigHighlyDAMPEDbis,LasiePANtrigHighlyDAMPED,TriggianiHIGHLYDIEdaSOLO}. 
 
 A more general version of the equation with memory~(\ref{eq:laFORMintegDIvCONuINTEGROdiff}) is the   system with persistent memory controlled on the boundary\footnote{If $K(t)$ is differentiable, MacCamy's trick can be used to rewrite the equation in an equivalent form but without the unbounded operators in the memory term.} described by
 
\[
v'(t)=\Delta v(t)+  \intt K(t-s) \Delta v(s)\ZD s +f(t) \quad\mbox{in $\ZOMq$}\,,\qquad v=u\ \mbox{on $\partial\ZOMq$}
\] 
The memory kernel $K(t)$ needs not be an exponential. We conjecture that similar ideas as those developed in this paper can be used to study this problem too.
 
%%%%%%%%%%%%%%%
%%%%%%%%%%%%%
\appendix
\section{Ancillary results and postponed  proofs}

We collect in this appendix few ancillary results whose proof is standard and the proofs we postponed.

\subsection{\ZLA{AppePROpreHzt}the properties of $(I+\Lzt \Lzt^*)^{-1}$}
 
We recall $K=e^{At}AD\in L^{p_0}(0,T;\mathcal{L}(U,H))$ with $p_0>1$. Let
\begin{equation}\ZLA{eq:AppePRIMAeqOperaForm}
\phi(\cdot;\zt)=(I+\Lzt \Lzt^*)^{-1} g(\cdot;\zt)\,.
\end{equation}
The functions  $\phi$ and $g$ are functions on $(\zt,T)$. The final time $T$ is fixed and we are interested in~(\ref{eq:AppePRIMAeqOperaForm}) for every $\zt\in[0,T]$  when either of the  following assumption hold: 
\begin{description}
\item[Assumption {\bf H1:}] $g=g(\cdot,\zt)\in L^2(\zt;T;H)$ and $\|g\|_{L^2(\zt;T;H)}<M$. The constant $M$ does not depend on $\zt$.
\item[Assumption {\bf H2:}] $g=g(\cdot,\zt)\in C([\zt;T];H)$ and $\|g\|_{C([\zt;T];H)}<M$. The constant $M$ does not depend on $\zt$.
\end{description} 
Note that assumption {\bf H2} implies {\bf H1} and that {\bf H1} implies  
\[
\phi(\cdot;\zt)\in L^2(\zt,T;H)\,,\qquad  \|\phi(\cdot;\zt)\|_{L^2(\zt;T;H)}<\|g\|_{L^2(\zt;T;H)}\,.
\]

We decouple Eq.~(\ref{eq:AppePRIMAeqOperaForm}) as follows:
\begin{subequations}
\begin{align}
\ZLA{Eq:appeDaPSIaPHI}
\phi(t;\zt)&=-\int_\zt^t K(t-\nu)\psi(\nu,\zt)\ZD\nu+g( t;\zt)
\\
\ZLA{Eq:appeDaPHIaPSI}
\psi(t;\zt)&=\int_t^T K^*(\nu-t)\phi(\nu;\zt)\ZD\nu\,.
\end{align}
\end{subequations}

The integrals in~(\ref{Eq:appeDaPSIaPHI}) and in~(\ref{Eq:appeDaPHIaPSI}) are both convolution integrals (of functions suitably extends to $\zzr$ with zero). 
So, Young inequalities  can be applied to both of them.  We recall the inequalities with reference to~(\ref{Eq:appeDaPHIaPSI}) but a similar statement holds for the convolution integral in~(\ref{Eq:appeDaPSIaPHI}): let $\phi\in L^{q}(\zt,T)$ and let
\[
\dfrac{1}{r}=\dfrac{1}{p_0}+\dfrac{1}{q}-1\,.
\]
Then we have:
\begin{enumerate}
\item if $r\leq 0$ then   $\psi(\cdot;\zt)\in C([\zt,T];U)$ (since $p_0>1$) and  
\[
\|\psi(\cdot;\zt)\|_{ C([\zt,T];U)}
\leq
\|K\|_{L^{p_0}(0,T;\mathcal{L}(U,H)) } \|\phi(\cdot;\zt)\|_{L^q(\zt,T;H)} 
\,.
\]
\item if $r>0$ then    $\psi(\cdot;\zt)\in L^r(\zt,T;U)$ and 
\[ 
\|\psi(\cdot;\zt)\|_{ L^r(\zt,T;U)}\leq
  \|K\|_{L^{p_0}(0,T;\mathcal{L}(U,H)) } \|\phi\|_{L^q(\zt,T;H)}
\,.
\]
\end{enumerate}

The important observation is that $\|K\|_{L^{p_0}(0,T;\mathcal{L}(U,H)) }$ does not depend on $\zt\in [0,T]$.

We apply this result to the system of the equations~(\ref{Eq:appeDaPSIaPHI}) and~(\ref{Eq:appeDaPHIaPSI}). We consider the case
\begin{equation}
\ZLA{eqAPPEipoFONDA}
\phi(\cdot;\zt)\in L^{q}(\zt,T;U)\,,\qquad \|\phi(\cdot;\zt)\|_{L^{q}(0,T;U)}<M\,.
\end{equation}
Condition~(\ref{eqAPPEipoFONDA}) holds with $q=2$ in the case of our interest, i.e. under either the assumption {\bf H1} or {\bf H2} since  $\phi$ is given by~(\ref{eq:AppePRIMAeqOperaForm}).

 We associate respectively   to~(\ref{Eq:appeDaPHIaPSI})  and to~(\ref{Eq:appeDaPSIaPHI}) 
 the numbers
 \begin{subequations}
\begin{align}
\ZLA{YoungPERpsi}
 \dfrac{1}{r_\psi} &=\dfrac{1 }{p_0}+\dfrac{1}{q}-1\\
 \ZLA{YoungPERphi}
\dfrac{1}{r_\phi}&=\dfrac{1}{p_0}+\dfrac{1}{r_\psi} -1=\dfrac{2}{p_0}+\dfrac{1}{q}-2
\end{align}
\end{subequations}
so that
\begin{equation}\ZLA{IneqCONSEyoungphi}
\dfrac{1}{r_\phi}-\dfrac{1}{q}=-2\left (1-\dfrac{1}{p_0}\right )=-\zaa<0\quad \mbox{ since $p_0 >1$}\,.
\end{equation}

We have
\begin{Lemma}\ZLA{LemmaAPPEperYOUNG}
  Assumption {\bf H1} implies:
  \begin{enumerate}
 \item\ZLA{I1LemmaAPPEperYOUNG} $\phi\in L^2(\zt,T;H)$ and $\|\phi\|_{L^2(\zt,T;H)}\leq M$ ($M$ does not depend on $\zt$);
  \item\ZLA{I2LemmaAPPEperYOUNG} $\psi\in L^{r_\psi}(0,T;U)$ with  $r_\psi>2$ and there exists $M$ which does not depend on $\zt$ for which $\|\psi\|_{L^{r_\psi}(0,T;U)}<M$.
  \end{enumerate}
\end{Lemma}
\zProof Statement~\ref{I1LemmaAPPEperYOUNG} follows from the definition~(\ref{eq:AppePRIMAeqOperaForm}) of $\phi$.  Statement~\ref{I2LemmaAPPEperYOUNG}  follows from~(\ref{YoungPERpsi}) and from the observation that 
\[
\dfrac{1}{r_\psi}=\dfrac{1}{p_0}-\dfrac{1}{2}<\dfrac{1}{2} \quad\mbox{because $p_0>1$}\,.\zdiaform
\]

We prove:
\begin{Theorem}\ZLA{teoAPPEregolarPRPsoluFred}
The function $t\mapsto \phi(t;\zt)$ has the following properties:
\begin{enumerate}
\item\ZLA{I1teoAPPEregolarPRPsoluFred} if $g\in C([\zt,T];H)$   then there exists $M$ such that

\[%\begin{equation}
%\ZLA{eq:PROPRIEPHIconti}
\mbox{$\phi\in  C([\zt,T];H)$\  and\quad $\|\phi\|_{C([\zt,T];H)}<M\|g\|_{C([\zt,T]}$}\,.
\]%\end{equation}
The number $M$ does not depend on $\zt$ and the transformation $g\mapsto \phi$ is linear and continuous in $C([\zt,T];H)$.
 
\item\ZLA{I2teoAPPEregolarPRPsoluFred} if $g\in L^{p_1}(\zt,T;H)$ with $p_1\geq 2$  then   
\[%\begin{equation}
%\ZLA{eq:PROPRIEPHIconti}
\mbox{$\phi\in   L^{p_1}(\zt,T;H)$ and  $\|\phi\| _{L^{p_1}(\zt,T;H)}<M\|g\|_{ L^{p_1}(\zt,T;H)}$}\,.
\]%\end{equation}
The number $M$ does not depend on $\zt$ and the transformation  $g\mapsto \phi$ is continuous in $L^{p_1}(\zt,T;H)$.
\end{enumerate}
 \end{Theorem}
\zProof We prove the statement~\ref{I1teoAPPEregolarPRPsoluFred}.
We know that~(\ref{eqAPPEipoFONDA}) holds with $q=2$ so that~(\ref{YoungPERpsi}) implies $\psi\in L^{r_\psi}(\zt,T;U)$ with $1/r_\psi=1/p_0-1/2$. Inequality~(\ref{YoungPERphi}) implies that the convolution integral in~(\ref{Eq:appeDaPSIaPHI}) belongs to $L^{r_\phi}(\zt,T;H)$ and
\[
\frac{1}{r_\phi}=\dfrac{2}{p_0}-\dfrac{3}{2}\quad\mbox{so that}\quad \dfrac{1}{r_\phi}=\dfrac{1}{2}-\zaa
\]
where $\zaa>0$ is the number in~(\ref{IneqCONSEyoungphi}). The assumed continuity of $g$ implies $\phi\in L^{r_\phi}(\zt,T;H)$ too. Young inequalities shows also  $\|\phi\|_{ L^{r_\phi}(\zt,T;H)}<M$, a number which does not depend on $\zt$.

The argument can be iterated with $r_\phi$ in the place of the exponent $2$: we find $\phi\in  L^{r_{\phi,2}}$ with
\[
\dfrac{1}{r_{\phi,2}}=\dfrac{1}{2}-2\zaa\,.
\]
$k$-iterations, with $k\geq 1/(2 \zaa )$, show that
\begin{equation}
\ZLA{eq:PROPRIEPHIconti}
\mbox{$\phi\in  C([\zt,T];H)$  and $\|\phi\|_{C([\zt,T];H)}<M$ ($M$ does not depend on $\zt)$}\,.
\end{equation}

The  proof of the statement~\ref{I2teoAPPEregolarPRPsoluFred} is similar. The iteration finishes when it is obtained that the convolution integral   is of class $L^{p_1}$ like $g$.\zdia

Now we compare $\phi(\cdot;\zt)$ and $\phi(\cdot;\zt+h)$ on the interval $[\max\{\zt,\zt+h\},T]$ (here $h$ needs not be positive).  
\begin{Theorem}\ZLA{teoAPPEregolarPRPDIpedaTAU}
Let $h\in\zzr$ and let us assume:
\begin{enumerate}
\item\ZLA{I1teoAPPEregolarPRPDIpedaTAU} the function $g(\cdot;\zt)$ is continuous on $[\zt,T]$ for every $\zt\in [0,T)$.
\item\ZLA{I2teoAPPEregolarPRPDIpedaTAU} there exists $M_0$ such that $\|g(t,\zt)\|_H<M_0$ on $0\leq \zt\leq t\leq T$. 
\item\ZLA{I3teoAPPEregolarPRPDIpedaTAU}
the following continuity property: 
\begin{align*}
&\mbox{for every $\ZEP>0$ there exists $h_{g,\ZEP}>0$ (independent of $\zt$) such that}\\
&\mbox{if\  $0\leq |h|<h_{g,\ZEP}$ then $\|g(\cdot,\zt+h)-g(\cdot;\zt)\|_{C([\max\{\zt,\zt+h\},T];H)}\leq \ZEP$}\,.
\end{align*}
\end{enumerate}
Then, the difference $\phi(\cdot,\zt+h)-\phi(\cdot;\zt)$ is continuous on $[\max\{\zt,\zt+h\},T]$
and there
 exist $h_{\phi,\ZEP}$ and a constant $M$ (both independent from $\zt$) such that 
\[
0\leq |h|<h_{\phi,\ZEP} \ \implies\ 
 \|\phi(\cdot,\zt+h)-\phi(\cdot;\zt)\|_{C([\max\{\zt,\zt+h\},T];H)}\leq M\ZEP \,.
\] 
\end{Theorem}
 \zProof 
We consider $h>0$ so that both the functions are defined on $[\zt+h,T]$. On this interval we have
\begin{multline*}
\phi(t;\zt+h)-\phi(t;\zt)=
-\int_{\zt+h}^tK(t-\nu)\int_\nu^T K^*(r-\nu)\left [
\phi(r;\zt+h)-\phi(r;\zt)
\right ]\ZD r\,\ZD\nu\\
+\int_{\zt}^{\zt+h}K(t-\nu)\int_\nu^T K^*(r-\nu)\phi(r;\zt)\ZD r\,\ZD\nu +\left [g(t,\zt+h)-g(t;\zt)\right ]\,.
\end{multline*} 
Theorem~\ref{teoAPPEregolarPRPsoluFred} can be applied to the difference $\phi(t;\zt+h)-\phi(t;\zt)$  in the interval $[\zt+h,T]$. We get that $\phi(\cdot;\zt+h)-\phi(\cdot;\zt)$ is continuous and we get the following inequality
 \begin{multline} \ZLA{EqDIteoAPPEregolarPRPDIpedaTAU}
 \|\phi(\cdot;\zt+h)-\phi(\cdot;\zt)\|_{C([\zt+h,T];H} 
\leq M\max_{t\in [\zt+h,T]}
\Biggl \{
\|g(t;\zt+h)-g(t;\zt)\|_H
\Biggr.\\
\left.+\left \| \int_\zt^{\zt+h}  K(t-s)\int_s^T K^*(\nu-s)\phi(\nu;\zt)\ZD\nu\,\ZD s\right \|_H
\right\}\,.
 \end{multline}
The required property of $\phi $ follows from the assumed property of $g$ and from the statement~\ref{I1teoAPPEregolarPRPsoluFred} of Theorem~\ref{teoAPPEregolarPRPsoluFred} applied to $\phi(\cdot;\zt)$.

The case $h<0$ is proved in a similar way. The difference to be noted is that
the maximum in~(\ref{EqDIteoAPPEregolarPRPDIpedaTAU}) is computed on the fixed interval $[\zt,T]$ and that the last integral in~(\ref{EqDIteoAPPEregolarPRPDIpedaTAU}) depends on $\phi(\nu;\zt-h)$. The statement~\ref{I1teoAPPEregolarPRPsoluFred} of Theorem~\ref{teoAPPEregolarPRPsoluFred} shows the existence of a bound for
$\|\phi(\cdot;\zt-h)\|$ which does not depend on $h$.\zdia

 \subsection{\ZLA{subseLIMITEincremQUOT}The computation of two limits}
We recall the assumptions
\begin{enumerate}
\item the control $u$ is continuous on $[\zt_0,T]$;
\item $\Sta_{\zt_0}\in \Dom\,\mathcal{A}_{\zt_0} $ i.e.  $\Xi_{\zt_0}\in\Dom\,\widetilde{\mathcal{A}}_{\zt_0}$ and $ \hat  y_{\zt_0}\in H$.
\end{enumerate} 

   We   compute  
 \begin{subequations}
  \begin{align}
  \ZLA{LimSECOeleProinteqRaIncre2}
&\lim _{h\to 0^+} \left [\Gamma_{\ztn}\dfrac{ \Mh\Xi(\ztn+h)- \Xi(\ztn)}{h}  \right ] \\
 \ZLA{LimSECOeleProinteqRaIncre2SINISTRA}
 & \lim _{k\to 0^+}
 \left [\Gamma_{\ztn-k}\dfrac{ \Mk\Xi(\ztn )- \Xi(\ztn-k)}{k}  \right ] 
  \end{align}
  \end{subequations}
where $\Xi_{r}=\mathcal{C}^{(1,2)}_r\Sta(r)$ is the $M^2_r$-component of the state.   
  Our goal  is not
   the sole computation of the limit for every $s$. \emph{We  prove   existence of the limit in an $L^p$-norm, $p>1$, so that   the limit 
can be exchanged with an integral.}    The operator $\Gamma_{\ztn}$ in front of the incremental quotient has its role in the proof.

The computation of the limits requires the computation of several intermediate limits and we obtain also   the important equality
\[
\lim _{h\to 0^+}\mbox{(\ref{LimSECOeleProinteqRaIncre2})}=\lim _{k\to 0^+}
\mbox{(\ref{LimSECOeleProinteqRaIncre2SINISTRA})}\,.
\]
This equality is proved explicitly for the first intermediate limit. Then we confine ourselves to study the limits for $h\to 0^+$ and we leave the limits for $k\to 0^+$ to the reader.

To proceed, we represent $\Gamma_\ztn$ as
  
\begin{multline}\ZLA{EqDEfiGamma12}
 \Gamma_\ztn\Xi_\ztn = \Gamma^1_{\ztn}\Xi_\ztn +
 \Gamma^2_{\ztn}\Xi_\ztn\quad  \mbox{($\Xi_{\ztn}=(v(\ztn),v(\ztn-\cdot))$ and $v\ZCD$     in~(\ref{eq:DEFIXIt})-(\ref{eq:DEFIvt}))}\,,\\[1mm]
 \left\{\begin{array}{lll}
\displaystyle  (\Gamma^1_{\ztn}\Xi_\ztn)\ZCD &=& Z(\cdot-\ztn)  v(\ztn)\\[1mm]
\displaystyle  (\Gamma^2_{\ztn}\Xi_\ztn)\ZCD&=& -\int_{\ztn}^\cdot Z(\cdot-r)
 e^{-(r-\ztn)}\left [\int_0^{\ztn} e^{-\nu}v(\ztn-\nu)\ZD\nu\right]\ZD r\,.
\end{array}\right.
\end{multline}

 In order to understand the following arguments, the reader should keep in mind the following facts:
 \begin{enumerate}
 \item $\Gamma^1_{\ztn} $ depends solely on the first component of $\Xi_\ztn$ while $\Gamma^2_{\ztn}$ depends solely on the second component. Both the components depend  on $\Sta_{\zt_0}$ and on $u$.
 
 \item
    both   $\Gamma_{\ztn}$ and $\Xi(\ztn)$ do depend on $\ztn$ but the
    incremental quotient which enters~(\ref{LimSECOeleProinteqRaIncre2}) and in~(\ref{LimSECOeleProinteqRaIncre2SINISTRA}) is solely that of    
    $ \Xi(\ztn)$.\zdia 
 \end{enumerate}

We use~(\ref{VariatCONSTANTA}) with $\zt=\zt_0$  to represent
\begin{align*}
&\Xi(t)=\Xi_u(t)+\Xi_{\Xi_{\zt_0}}+\Xi_{\hat y_{\zt_0}}\,,\\
&
v(t)=\left\{\begin{array}{lll}
v_u(t)+ v_{\Xi_{\zt_0}}(t)+ v_{  \hat y_{\zt_0}}(t)&{\rm if}&t\geq \zt_0\\
  \xi_{\zt_0}(t)&{\rm if}&0\leq t\leq \zt_0\,.
\end{array}\right.
\end{align*}
The notation is self explanatory.

\emph{Our assumptions  imply continuity of $v(t)$ for $t\geq 0$, see the statement~\ref{I2BTeo:preREGOLARITAstato} of  Theorem~\ref{Teo:preREGOLARITAstato}. }
  
    First we study the contribution to the limit of $\Gamma^1_\ztn  $ and then that of $\Gamma^2_\ztn$.
  
  \subsubsection*{The contribution of $\Gamma^1_\ztn  $}  
 
  First we consider to the limit of the contribution of $\Xi_u$  and we recall
  \begin{equation}\ZLA{eq:Controv1NellaGamma1ePOIgAMMA2}
v_u(t)=\left\{\begin{array}{lll}  \int_{\zt_0} ^t Z(t-s)(-A D) u(s)\ZD s&{\rm if}& t\geq \zt_0\\
0 &{\rm if}& 0\leq t\leq \zt_0 \,.
\end{array}\right.
  \end{equation}
We study   the limit~(\ref{LimSECOeleProinteqRaIncre2}) i.e. we study the limit for $h\to 0^+$ of
 
\begin{subequations}
\begin{align}
\nonumber
&\left (\Gamma^1_{\ztn}\dfrac{ \Mh\Xi_u(\ztn+h)- \Xi _u(\ztn)}{h}\right )(\cdot)\\
&= 
\ZLA{eqRaIncre1CONTR1diU}
  Z(\cdot-\ztn)(-A)^{\ZSI }\dfrac{1}{h}\int_\ztn ^{\ztn+h} Z(\ztn+h-r) (-A)^{1-\ZSI}D u(r)\ZD r\\
&+
\ZLA{eqRaIncre2CONTR1diU}
Z(\cdot-\ztn)(-A)^{ \ZSI+\ZEP}  \int_{\zt_0}^{\ztn}\dfrac{\left [   
Z(\ztn+h-r)-Z(\ztn-r) \right ](-A)^{-\ZEP}[(-A)^{1-\ZSI}Du(r)]}{h}   
\ZD r
\,.
\end{align}
\end{subequations}

The  integrand  in~(\ref{eqRaIncre1CONTR1diU}) is continuous   so that 
\begin{align*}
\lim_{h\to0^+}
&\dfrac{1}{h}\int_\ztn ^{\ztn+h} Z(\ztn+h-r) (-A)^{1-\ZSI}Du(r)\ZD r=(-A)^{1-\ZSI}Du(\ztn)\,,\\
&\left \|\dfrac{1}{h} \int_\ztn ^{\ztn+h} Z(\ztn+h-r) (-A)^{1-\ZSI}Du(r)\ZD r\right \|<M\,.
\end{align*}
It follows
  \begin{subequations}
 
   \begin{equation}\ZLA{eq:limite1DIeqRaIncre2uRHO}
\lim _{h\to0^+}\mbox{(\ref{eqRaIncre1CONTR1diU})}=
- Z(\cdot-\ztn)A D u(\ztn)\quad \mbox{  in the space $L^{p_0}(\zt,T;H)$}\,. 
   \end{equation}
 In order to study the limit of~(\ref{eqRaIncre2CONTR1diU}) we use Theorem~\ref{Teo:limRappINCREz}:  for every $v\in H$ we have
    \begin{align*}
& \lim_{h\to 0^+  }\dfrac{   
 Z(\cdot+h)-Z(\cdot)
}{h}(-A)^{-\ZEP}v= Z'(\cdot)(-A)^{-\ZEP}v\quad \mbox{ in an $L^p$ norm, $p>1$} \\
&\mbox{and inequality~(\ref{equaDELTeo:limRappINCREz}) holds}\,.
   \end{align*}
    
   The expression of $Z'(t)$ is given in~(\ref{equaDELTeo:limRappINCREzVALUElim}). The function $(-A)^{1-\ZSI}Du\ZCD$ is continuous. Young inequalities imply that the convolution in~(\ref{eqRaIncre2CONTR1diU}) is a continuous and bounded $H$-valued function of $\ztn$ and $h$ and that
  the limit of the integral (which exists in $H$) is equal to
   \[
\int _{\zt_0}^\ztn   Z'(\ztn-r)(-A)^{-\ZEP}[(-A)^{1-\ZSI}Du(r)]\ZD r\,.
   \]
Provided that $\ZEP\in (0,1-\ZSI)$, the function $Z(\cdot-\ztn)(-A)^{{ \ZSI+\ZEP}}\in L^p(\ztn,T;\mathcal{L}(H))$
  (for a suitable $p>0$). So, there exists $p>1$ such that in~$L^p(\zt_0,T;H)$ we have:
  \begin{equation}\ZLA{eq:limite2DIeqRaIncre2uRHO}
\lim_{h\to 0^+}  \mbox{(\ref{eqRaIncre2CONTR1diU})}=
Z(\cdot-\ztn)\left \{(A+I)v_u(\ztn)-\int_{\zt_0}^{\ztn}e^{-(\ztn-s)}v_u(s)\ZD s\right \}
\,.
  \end{equation}

   When studying the contribution to the  limit~(\ref{LimSECOeleProinteqRaIncre2SINISTRA})  we consider the limit for $k\to 0^+$ of
   \begin{multline*}
 \left (  \Gamma^{(1)}_{\ztn-k}\dfrac{\Mk \Xi(\ztn)-\Xi(\ztn-k)}{k}\right )\ZCD=
 Z(\cdot-(\ztn-k))\dfrac{v(\ztn)-v(\ztn-k)}{k}\\
 =
 Z(\cdot-(\ztn-k))\dfrac{1}{k}\int _{\ztn-k}^\ztn Z(\ztn-r)(-AD)u(r)\ZR r\\
  + Z(\cdot-(\ztn-k))\int_\zt^{\ztn-k}\dfrac{Z(\ztn-k-r)-Z(\ztn-r)}{-k}(-AD) u(r)\ZR r\,.
   \end{multline*}
  Computations similar to the ones above give
   \begin{multline} \ZLA{LimSECOeleProinteqRaIncre2destraEsinistra}
   \lim _{h\to 0^+}\left (\Gamma^1_{\ztn}\dfrac{ \Mh\Xi_u(\ztn+h)- \Xi _u(\ztn)}{h}\right )(\cdot)\\
   =\lim _{k \to 0^+}\left ( \Gamma^{(1)}_{\ztn-k}\dfrac{\Mk \Xi_u(\ztn)-\Xi_u(\ztn-k)}{k }\right )\ZCD =- Z(\cdot-\ztn)A D u(\ztn)\\
   +Z(\cdot-\ztn)\left \{(A+I)v_u(\ztn)-\int_{\zt_0}^{\ztn}e^{-(\ztn-s)}v_u(s)\ZD s\right \}\,.
   \end{multline}   
   \begin{Remark}\ZLA{Rema:REGlimCompTWOlimsGAMMA1}{\rm
The statements~\ref{I5TeoPropriZt} and~\ref{I4TeoPropriZt} of Theorem~\ref{TeoPropriZt}  and Theorem~\ref{Teo:PROpreidelRappoINCREvU}  show that the limit belongs to $C([\ztn,T];(\Dom (-A)^{ \ZSI+\ZEP})')$ (any $\ZEP>0$)   
and to $L^{p_0}(\ztn,T;H)$ and that  its $L^{p_0}$-norm  is bounded on $[\ztn,T]\subseteq [\zt_0,T]$, uniformly respect to $\zt_0\in [0,T]$ ($p_0$ is the exponent in~(\ref{eqDefiEXPLIdi]ZSIeP})).\zdia
   }
      \end{Remark}
    
   \emph{ From now on, we confine ourselves to study the limits respect to $h $ and we leave to the reader the corresponding limits respect to $k$.}
  
We consider the contribution of
 
 \begin{multline*} 
   v_{\Xi(\zt_0),  y _{\zt_0}}(t)=
  v_{\Xi_{\zt_0}}(t)+v_{ \hat  y_{\zt_0}}(t)\\
  =\left\{\begin{array}{l} 
  % \displaystyle
    Z(t-\zt_0) v(\zt_0)
  \\   
%   \displaystyle  
\quad   +\int_{\zt_0}^t Z(t-s)e^{-(s-\zt_0)}\left [ \hat   y_{\zt_0}
- \int _0^{\zt_0}e^{-r}  \xi_{\zt_0}(\zt_0-r)\ZD r\right ]
\ZD s\,,\qquad t\geq \zt_0 \\[2mm]
%\displaystyle
  \xi_{\zt_0} (t)\qquad t\in[0,\zt_0]\,.
\end{array}\right.
 \end{multline*}

  The first addendum $Z(t-\zt_0)v(\zt_0)$ in the right side is of class $C^1$ because the assumption  $\Xi_{\zt_0}\in\Dom\,\mathcal{A}_{\zt_0}$ implies  $v(\zt_0)= \hat v_{\zt_0}\in\Dom\,A$. The derivative of the second and third addenda can be computed 
  by using  Theorem~\ref{Teo:limRappINCREz}. The computations are similar to those used to compute the contribution of $v_u$ and it is omitted. We state the result: the limit   converges in an $L^p$-norm (with $p>1$) and the limit  is
  \begin{equation}\ZLA{eq:limiteDATIiniHATyIeqRaIncre2uRHO}
  Z(\cdot-\ztn)
   \left  \{
   (A+I)
    v_{\xi_{\zt_0},\hat y_{\zt_0}}
    -\int_0^\ztn e^{-s} v_{\xi_{\zt_0},\hat y_{\zt_0}}(\ztn-s)\ZD s+e^{-\ztn  }  y _{\zt_0}
\right\} \,.
  \end{equation}
 \end{subequations}

By summing~(\ref{eq:limite1DIeqRaIncre2uRHO})-(\ref{eq:limiteDATIiniHATyIeqRaIncre2uRHO}) (and similar limits respect to $k$) we get

\begin{multline}\ZLA{eqCONTRItotGamma1}
\lim _{h\to 0^+}
\left (\Gamma^1_{\ztn}\dfrac{ \Mh\Xi (\ztn+h)- \Xi  (\ztn)}{h}\right )(\cdot)\\
=
 \lim _{k \to 0^+}\left ( \Gamma^{(1)}_{\ztn-k}\dfrac{\Mk \Xi (\ztn)-\Xi (\ztn-k)}{k }\right )\ZCD \\
=Z(\cdot-\ztn)\left [(A+I)v(\ztn)-\int_0^\ztn e^{-s}v(\ztn-s)-ADu(\ztn)+e^{-\ztn}\hat y_{\zt_0}\right ]\\
=Z(\cdot-\ztn)\mathcal{C}^{1}_\ztn\left [ \mathcal{A}_\ztn\Sta(\ztn)+\mathcal{B}_\ztn u(\ztn)\right ] =\Gamma^1_{\ztn}\mathcal{C}^{1,2}_\ztn\left [ \mathcal{A}_\ztn\Sta(\ztn)+\mathcal{B}_\ztn u(\ztn)\right ]  
\end{multline}
since $\Gamma^1_{\ztn} $ depends on the first component of $\Xi(\ztn)$ hence of $\Sta(\ztn)$.

Convergence is in an $L^p $ norm, $p>1$.

  \subsubsection*{The contribution of $\Gamma^2_\ztn  $}  
Also in this case we confine ourselves to study the limit respect to $h$ and we leave to the reader the similar computations respect to $k$.  
  
 We   study convergence of
\begin{multline}\ZLA{Eq:GENEnelCAsoGamma^2}
-\int_\ztn^\cdot Z(\cdot-r) e^{-(r-\ztn)}\int_0^{\ztn}e^{-\nu}\dfrac{v(\ztn-\nu+h)-v(\ztn-\nu)}{h}\ZD\nu\,\ZD r
\\
=-\int_\ztn^\cdot Z(\cdot-r) e^{-(r-\ztn)}\int_0^{\ztn}e^{-\nu}\dfrac{\xi(\nu-h,\ztn)-\xi(\nu,\ztn)}{h}\ZD\nu\,\ZD r
\end{multline}
(the relation between $v$ and $\xi$ is in~(\ref{eq:DEFIXIt})-(\ref{eq:DEFIvt})).

We consider   the contribution of $ u$. The function $v_u$ is given in~(\ref{eq:Controv1NellaGamma1ePOIgAMMA2}): it is continuous on $[0,T]$. Theorem~\ref{Teo:PROpreidelRappoINCREvU} shows that the incremental quotient of $(-A)^{-\ZSI-\ZEP} v_u(t)$ converges in $C([\zt_0,T];H)$ for every $\ZEP>0$ and that inequality~(\ref{eq:LimiQuozIndFivU}) holds. It converges also in $C([0,\zt_0 ];H)$ since the function is zero on this interval. Hence this incremental quotient converges (in particular) in $H^1(0,T;H)$.

We use this observation and we write~(\ref{Eq:GENEnelCAsoGamma^2}) (with $v=v_u$) as
\[
-\int_\ztn^\cdot Z(\cdot-r)(-A)^{\ZSI+\ZEP} e^{-(r-\ztn)}\int_0^{\ztn}e^{-\nu}(-A)^{-(\ZSI+\ZEP)} \dfrac{v_u(\ztn-\nu+h)-v_u(\ztn-\nu)}{h}\ZD\nu\,\ZD r\,.
\]
 The limit for $h\to 0^+$ is
  \begin{subequations}
 \begin{multline} \ZLA{eqCONTRIvuAgAmmA2}
- \int_\ztn^\cdot Z(\cdot-r) e^{-(r-\ztn)}\int_0^{\ztn}e^{-\nu}(-\partial_\nu)  v_u(\ztn-\nu) \ZD\nu\,\ZD r\\=
- \int_\ztn^\cdot Z(\cdot-r) e^{-(r-\ztn)}\int_0^{\ztn}e^{-\nu}(-\partial_\nu)  \xi_u(\nu,\ztn) \ZD\nu\,\ZD r=\Gamma^2_\ztn \left (\widetilde{\mathcal{A}}_\ztn \Xi_u \right )
 \end{multline}
since $\Gamma^2_\ztn$ depends solely on the second component of $\Xi_\ztn$.

The contribution of $v_{\hat y_{\zt_0}} $ is treated in a similar way since by assumption $\hat y_{\zt_0}\in H$ so that $v_{\hat y_{\zt_0}}$ (wich is zero for $t<\zt_0$) is of class $H^1(0,T;H)$.
So, when $v=v_{ \hat  y_{\zt_0}}$, the limit~(\ref{Eq:GENEnelCAsoGamma^2}) is
\begin{equation}
\ZLA{eqCONTRIvhatyAgAmmA2}
\Gamma^2_\ztn  \left (\widetilde{\mathcal{A}}_\ztn   \Xi_{ \hat  y_{\zt_0}}\right )
=
\Gamma^2_\ztn \left [ \widetilde{\mathcal{A}}_\ztn   \Xi_{\hat  y_{\zt_0}} 
+Y(\ztn)
\right ]\,.
\end{equation}
The second equality is true since  the  second component of $Y $ is zero.

Finally we examine the contribution of  $\Xi_{\zt_0}$.
 We recall the assumption  $\Xi_{\zt_0}\in \Dom\,\widetilde{\mathcal{A}}_{\zt_0}$
i.e. 
$ \xi_{\zt_0}(\zt_0-\cdot)\in H^1(0,\zt_0;H)$ and $\xi(0) =v(\zt_0)\in\Dom\, A$. It follows that $t\mapsto v_{\Xi_{\zt_0}}(t)$ is of class $H^1(0,T;H)$ (see the statement~\ref{I1CTeo:preREGOLARITAstato} of Theorem~\ref{Teo:preREGOLARITAstato}).
Thanks to this observation it is easy to compute that the contribution of $\Xi_{\zt_0}$ to the limit is
 
\begin{equation}
\ZLA{eqCONTRIvXIAgAmmA2}
\Gamma_\ztn^2 \left (\widetilde{\mathcal{A}}_\ztn \Xi_{\Xi_{\zt_0}} \right ) \,.
\end{equation}
\end{subequations}

We sum~(\ref{eqCONTRItotGamma1})     with~(\ref{eqCONTRIvuAgAmmA2})-(\ref{eqCONTRIvXIAgAmmA2}). Then we perform similar computations for the 
limit~(\ref{LimSECOeleProinteqRaIncre2SINISTRA}). We get:
 
\begin{multline}\ZLA{sectLIMIINCREquotEquiFIN}
%\lim _{h\to 0^+}\mbox{(\ref{LimSECOeleProinteqRaIncre2})}=
\lim _{h\to 0^+}\left [\Gamma_{\ztn}\dfrac{ \Mh\Xi(\ztn+h)- \Xi(\ztn)}{h}  \right ] 
=
 \lim _{k\to 0^+}
 \left [\Gamma_{\ztn-k}\dfrac{ \Mk\Xi(\ztn )- \Xi(\ztn-k)}{k}  \right ] \\
=
\Gamma_\ztn\left (\,
\widetilde{\mathcal{A}}_\ztn  \Xi(\ztn)+\widetilde{\mathcal{B}}_\ztn u(\ztn)+Y(\ztn)
\,
\right )=
\Gamma_\ztn\mathcal{C}^{(1,2)}_\ztn\left (\,
 \mathcal{A}_\ztn  \Sta(\ztn)+\mathcal{B}_\ztn u(\ztn) 
\,
\right )
\end{multline}
 
 \begin{Remark}\ZLA{Rema:REGlimCompTWOlims}{\rm
The same observation as in Remark~\ref{Rema:REGlimCompTWOlimsGAMMA1} holds:  the limit belongs to $C([\ztn,T];(\Dom (-A)^{ \ZSI+\ZEP})')$ (any $\ZEP>0$)   
and to $L^{p_0}(\ztn,T;H)$ and   its $L^{p_0}$-norm  is bounded on $[\ztn,T]\subseteq [\zt_0,T]$, uniformly respect to $\zt_0\in [0,T]$ ($p_0$ is the exponent in~(\ref{eqDefiEXPLIdi]ZSIeP})).\zdia
   }
      \end{Remark}
%\emph{The limit converges in $L^p$ with $p>1$.}

 \subsection{\ZLA{AppeDeriValue}Computation of the derivative of the value function}
  
We state a lemma,  consequence of    Theorem~\ref{teoAPPEregolarPRPDIpedaTAU},  which is repeatedly used.  
\begin{Lemma}\ZLA{LemmaINAppeDeriValue} Let $0\leq\zt_0\leq \ztn<T$ and let $\Sta_\ztn\in\Dom\,\mathcal{A}_\ztn$. Let $h_{\ztn }$   be the function  in~(\ref{eq:DefiOPERfonDam}).
The function $s\mapsto  [H_{\ztn }h_{\ztn }](s)=[(I+L_\ztn L_\ztn^*)^{-1}h_{\ztn}](s)$  
 has the following properties:
\begin{enumerate}
\item it is   bounded on the set $\ztn\leq s\leq T$, uniformly respect to $\ztn\geq \zt_0$.
\item
Let $\ztn_0>0$. We have
\[
\lim _{h\to 0} [H_{\ztn+h }h_{\ztn+h }](s)= [H_{\ztn }h_{\ztn }](s)
\] 
uniformly for $ \ztn_0\leq s\leq T$.
\end{enumerate}
In particular, let us extend $[H_{\ztn+h }h_{\ztn+h }](s)=0$ when $\zt_0\leq s\leq \ztn+h$. 
We have:
\[
\lim _{h\to 0}  [H_{\ztn+h }h_{\ztn+h }]\ZCD= [H_{\ztn  }h_{\ztn  }]\ZCD 
\]
in $L^q(\zt_0,T;H)$ (hence also in $L^q(\ztn,T;H)$) for every $q\geq 1$.
\end{Lemma}

This lemma holds because $g(\cdot;\ztn)=h_\ztn\ZCD$ satisfies the assumptions of Theorem~\ref{teoAPPEregolarPRPDIpedaTAU}.

Now we compute the limits of the eight addenda in~(\ref{eqINCRequoDOLOstato})-(\ref{DaSUBEQparteCOMUNE}), respectively 
in~(\ref{eqINCRequoDOLOstato}')-(\ref{DaSUBEQparteCOMUNE}').
First we compute the limits of the two addenda in~(\ref{DaSUBEQparteCOMUNE}) and in~(\ref{DaSUBEQparteCOMUNE}')  because the increment in these terms involves both the operator and the state. So, the results we find have to be summed in part to the limits of the addenda in~(\ref{eqINCRequoDOLOstato})
(respectively in~(\ref{eqINCRequoDOLOstato}')) and in part to those of the addenda in~(\ref{eqINCRequoDOLOriccati}) (respectively in~(\ref{eqINCRequoDOLOriccati}')).

\subsubsection{\ZLA{ParaDaSUBEQparteCOMUNE}The limits of the addenda in~(\ref{DaSUBEQparteCOMUNE}) and  in~(\ref{DaSUBEQparteCOMUNE}')}

We recall that   $y(t)$ is the solution of
\begin{equation}\ZLA{EquaDIyNellappePERultiLIMI}
y'=-y\,,\quad y(\zt_0)=\hat y_{\zt_0}\in H 
\end{equation}
but we do the   computation with any $C^1$-function $y  \in C^1([0,T];H)$    
    as in the requirement~\ref{I0AssuSuallaConveeqPERdopoBETAeDopoGAMMA} of the Assumption~\ref{AssuSuallaConveeqPERdopoBETAeDopoGAMMA} and  for consistency with~(\ref{EquaDIyNellappePERultiLIMI}) we   assume $y(0)=y_{\ztn}$.
    In this generality, we use
    \[
\hat Y_\zt  (s)=\int_\ztn^s Z(s-r) y(r-\zt )\ZD r\,.    
    \]

 The incremental quotients in~(\ref{subequatioPERdopobmdelta})  and in~(\ref{subequatioPERdopobmdelta}') are the right, respectively left, incremental quotients of
\[
\ztn\mapsto \int_\ztn^\cdot Z(\cdot-r)y(r-\ztn)\ZD r\,.
\]
So:
 
\begin{align}
\nonumber
 & \lim _{h\to 0^+}\mbox{(\ref{subequatioPERdopobmdelta})} 
 = \lim _{h\to 0^+}\mbox{(\ref{subequatioPERdopobmdelta}')} \\
\ZLA{subequatioPERdopobmdeltaSOLOoperatore}
&\quad = -\int_\ztn^T \Bigl\ZL 
  [H_\ztn h_\ztn](s), Z(s-\ztn)\hat y_\ztn \Bigr\ZR\ZD s\\
 \ZLA{subequatioPERdopobmdeltaSOLOstato}
 &\quad +
\int_\ztn^T \Bigl\ZL 
  [H_\ztn h_\ztn](s),\left [ \int_\ztn^s Z(s-r)[-  y'(r-\ztn)]\ZD r\right ] \Bigr\ZR\ZD s\,.
\end{align}

Analogously:
\begin{align}
\nonumber
&  \lim _{h\to 0^+}\mbox{(\ref{subequatioPERdopobmepsilon})}
  =  \lim _{h\to 0^+}\mbox{(\ref{subequatioPERdopobmepsilon}')}  \\
\ZLA{subequatioPERdopobmepsilonSOLOoperatore}
=&-\int_\ztn^T \Bigl\ZL 
 Z(s-\ztn)\hat y_\ztn , [H_\ztn h_\ztn](s)\Bigr\ZR\ZD s\\
 \ZLA{subequatioPERdopobmepsilonSOLOstato}
&+
\int_\ztn^T \Bigl\ZL 
\left [ \int_\ztn^s  Z(s-r)[ - y'(r-\ztn)]\ZD r\right ], [H_\ztn h_\ztn](s) \Bigr\ZR\ZD s\,.
\end{align}

\emph{Note continuity respect to to $\ztn$ of these expressions.}

\begin{Remark} \ZLA{RemaINTElimiDaSUBEQparteCOMUNE}{\rm
In order to interpret these  equalities we note:
\begin{enumerate}
\item\ZLA{I1RemaINTElimiDaSUBEQparteCOMUNE} the derivation of the two addenda~(\ref{subequatioPERdopobmdeltaSOLOoperatore}) and~(\ref{subequatioPERdopobmepsilonSOLOoperatore}) uses the sole
incremental quotient of the operator and  continuity of $y$. Equation~(\ref{EquaDIyNellappePERultiLIMI}) is not used.   
Let $\Sta_\ztn\in \Dom\,\mathcal{A}_{\ztn}$ be such that $\mathcal{C}^{(3)}_\ztn \Sta_\ztn=\hat y_\ztn$. We have
 
\begin{multline}\ZLA{subequatioPERdopobmSOLOoperatore}
\mbox{(\ref{subequatioPERdopobmdeltaSOLOoperatore})}+
\mbox{(\ref{subequatioPERdopobmepsilonSOLOoperatore})}= 
 \ZL 
 [  H_\ztn h_\ztn]\ZCD  , 
Z(\cdot-\ztn)(-\hat y_\ztn)
  \ZR _{L^2(\ztn,T;H)}
  \\
  +
  %%%%%%%%%%
  \ZL 
 Z(\cdot-\ztn)(-\hat y_\ztn)
 , [ H_\ztn h_\ztn]\ZCD  \ZR _{L^2(\ztn,T;H)}
\,.
\end{multline}
\item\ZLA{I2RemaINTElimiDaSUBEQparteCOMUNE} to interpret  the two addenda~(\ref{subequatioPERdopobmdeltaSOLOstato}) and~(\ref{subequatioPERdopobmepsilonSOLOstato}) we note that    the computation uses the  incremental quotient of $y$ but not that of   the Riccati operator. So, now \emph{we use explicitly that  $y$ solves~(\ref{EquaDIyNellappePERultiLIMI})} and we see  
\begin{multline}
 \ZLA{EqSommaDUEterminiINDaSUBEQparteCOMUNE} 
\mbox{(\ref{subequatioPERdopobmdeltaSOLOstato})}+\mbox{(\ref{subequatioPERdopobmepsilonSOLOstato})} =\ZL H_\ztn h_\ztn, L \mathcal{C}^{(3)}_{\ztn} \Sta'(\ztn)\ZR _{L^2(\ztn,T;H)}+\ZL  L \mathcal{C}^{(3)}_{\ztn}  \Sta'(\ztn), H_\ztn h_\ztn\ZR _{L^2(\ztn,T;H)} 
\\
=
\ZL H_\ztn h_\ztn, L \mathcal{C}^{(3)}_{\ztn}\left [\mathcal{A}_{\ztn} \Sta (\ztn)
+\mathcal{B}_\ztn u(\ztn)\right ]\ZR _{L^2(\ztn,T;H)}+
\\
%%%%%%%%%%
 \ZL   L \mathcal{C}^{(3)}_{\ztn}\left [\mathcal{A}_{\ztn} \Sta (\ztn)
+\mathcal{B}_\ztn u(\ztn)\right ],H_{\ztn}h_\ztn \ZR _{L^2(\ztn,T;H)}\,.
\end{multline}
We use Remark~\ref{RemaEXTEprojeC3} to apply the projection $\mathcal{C}^{(3)}_{\ztn}$  to $\Sta'(\ztn)$. In fact, $y$ takes values in $H$ thanks to the assumption  $y_{\zt_0}\in H$.\zdia
\end{enumerate}
}
 \end{Remark}
  \subsubsection{\ZLA{paraeqINCRequoDOLOstato}The   limits in~(\ref{eqINCRequoDOLOstato}) and in~(\ref{eqINCRequoDOLOstato}')}  
  We consider the limit of~(\ref{eqRaIncre2}).
  
The limit of the function~(\ref{LimSECOeleProinteqRaIncre2}) is studied in Sect.~\ref{subseLIMITEincremQUOT}. 
We proved that the limit   converges in an $L^p$-norm and that the limit is~(\ref{sectLIMIINCREquotEquiFIN}).
 
 The first element of the inner product in~(\ref{eqRaIncre2}) is the \emph{continuous function}
 \[
 s\mapsto\left [
H_\ztn h_\ztn  \right ](s)\,.
 \]
 This function does not depend on $h$ so that
 \begin{multline}\ZLA{ValeLimiZTNeqRaIncre2}
\lim _{h\to 0^+} \mbox{(\ref{eqRaIncre2})}=\int_{\ztn}^T\Bigl\ZL
H_\ztn h_\ztn  ,\Gamma_\ztn\left (\,
\widetilde{ \mathcal{A}} _\ztn \Xi(\ztn)+\widetilde{\mathcal{B}}_\ztn u(\ztn)+Y(\ztn)\,
\right)\Bigr\ZR_H\ZD s\\
%%%%%%%%%%%%
= \int_{\ztn}^T
\Bigl \ZL
 [H_\ztn h_\ztn](s)  ,\left [\Gamma_\ztn\mathcal{C}^{(1,2)}_\ztn\left (\,
 \mathcal{A} _\ztn \Sta(\ztn)+\mathcal{B}_\ztn u(\ztn)  
\right)\right ](s)
\Bigr\ZR_H\ZD s\\
=\ZL H_\ztn h_\ztn,\Gamma_\ztn\mathcal{C}^{(1,2)}_\ztn\left (\,
 \mathcal{A} _\ztn \Sta(\ztn)+\mathcal{B}_\ztn u(\ztn)  
\right)\ZR _{L^2(\ztn,T;H)}\,.
 \end{multline}
 \begin{Remark}\ZLA{remaREGULpropRIccati}{\rm
 We note:
 \begin{enumerate}
 \item The fact that the limit is an $L^2(\ztn,T;H)$-inner product is a smoothing property of the Riccati operator.
 \item
 Continuity of $H_\ztn h_\ztn$ and $L^p$ integrability of the second factor of the inner product (stated in Remark~\ref{Rema:REGlimCompTWOlims}) imply continuity of the inner product respect to~$\ztn$.\zdia
\end{enumerate}
}\end{Remark}
    
In order to compute the limit of of~(\ref{eqRaIncre3}) we note
\begin{multline*}
\Bigl\ZL 
 H_{\ztn+h}\Rh\Gamma_{\ztn} \frac{  \Mh\Xi(\ztn+h)-  \Xi(\ztn)}{h},h_{\ztn+h} 
\Bigr\ZR\\
=
\Bigl\ZL{\Rh} \Gamma_{\ztn} \frac{  \Mh\Xi(\ztn+h)-  \Xi(\ztn)}{h},
 H_{\ztn+h}
h_{\ztn+h} 
\Bigr\ZR\,.
\end{multline*}  
  The computations done in the case of the limit~(\ref{eqRaIncre2}) and Lemma~\ref{LemmaINAppeDeriValue} imply
 \begin{multline}\ZLA{ValeLimiZTNeqRaIncre3}
\lim _{h\to 0^+} \mbox{(\ref{eqRaIncre3})}= 
 \int_\ztn^T \Bigl\ZL
\Gamma_\ztn\left (\,
\widetilde{ \mathcal{A}} _\ztn \Xi(\ztn)+\widetilde{\mathcal{B}}_\ztn u(\ztn)+Y(\ztn)\,
\right),H_{\ztn}h_{\ztn} 
\Bigr\ZR\ZD s\\ 
=
  \int_\ztn^T \Bigl\ZL H_\ztn\Gamma_\ztn\mathcal{C}^{(1,2)}_\ztn\left (\,
 \mathcal{A} _\ztn \Sta(\ztn)+\mathcal{B}_\ztn u(\ztn)  
\right),h_\ztn\Bigr\ZR\ZD s
\\
=\ZL  \Gamma_\ztn\mathcal{C}^{(1,2)}_\ztn\left (\,
 \mathcal{A} _\ztn \Sta(\ztn)+\mathcal{B}_\ztn u(\ztn)  
\right),H_\ztn h_\ztn\ZR _{L^2(\ztn,T;H)}\,.
\end{multline}

Similar  arguments show
\[
\lim _{k\to 0^+} \mbox{(\ref{eqRaIncre2}')}=\lim _{h\to 0^+} \mbox{(\ref{eqRaIncre2})}\,,\quad \lim _{k\to 0^+} \mbox{(\ref{eqRaIncre3}')}=
\lim _{h\to 0^+} \mbox{(\ref{eqRaIncre3})} 
\]
 and, we repeat, these limits are continuous functions of $\ztn$. 
   \subsubsection{The sum of the four limits in~(\ref{eqINCRequoDOLOriccati}) and  of those in~(\ref{eqINCRequoDOLOriccati}')}
We repeat  that in these limits the increment $h$ affects the incremental quotient of the Riccati operator but not that of the state.
 
\paragraph{The limit of~(\ref{subequatioPERdopobmsigma}) and of~(\ref{subequatioPERdopobmsigma}')} 
The integrand  is a continuous function of $s$. So the limit of the incremental quotient is
\[
\ZL
  [H_{\ztn}h_{\ztn}](\ztn),h_\ztn(\ztn)\ZR
 =
\ZL [H_{\ztn} \Gamma_{\ztn} \Xi({\ztn})](\ztn),[\Gamma_{\ztn}\Xi({\ztn})](\ztn)\ZR
=\ZL [H_{\ztn} \Gamma_{\ztn} \Xi({\ztn})](\ztn),v(\ztn)\ZR \,.
\]

From $\phi(s)=(H(\ztn)\Gamma_{\ztn}\Xi({\ztn}))(s)$ we see
  then $\phi(\ztn)=[\Gamma_{\ztn}\Xi({\ztn})](\ztn) =v(\ztn)$.  Hence, by taking the minus sign into account,
\[%  \begin{equation}
%\ZLA{EqLIMITEesplicitoDIsubequatioPERdopoDELTAsigma}
  \lim _{h\to 0^+}\mbox{(\ref{subequatioPERdopobmsigma})} =-\|\hat v_\ztn\|^2=-\|\mathcal{C}^{(1)}_\ztn \Sta_\ztn\|^2=
\]%  \end{equation} 
and
  \begin{equation}\ZLA{EqLIMITEesplicitoDIsubequatioPERdopoDELTAsigma}
   \lim _{h\to 0^+}{(\ref{subequatioPERdopobmsigma})}
   =  \lim _{k\to 0^+}{(\ref{subequatioPERdopobmsigma}')} =-\|\hat v_\ztn\|^2=-\|\mathcal{C}^{(1)}_\ztn \Sta_\ztn\|^2\,.
\end{equation}
To compute the limit of  (\ref{subequatioPERdopobmsigma}') we apply Theorem~\ref{teoAPPEregolarPRPDIpedaTAU} to $\phi(\cdot;\ztn-k)=H_{\ztn-k}h_{\ztn-k}$.

\paragraph{The limit of~(\ref{subequatioPERdopobmzaa}) and of~(\ref{subequatioPERdopobmzaa}')} 
  
  We recall the definition~(\ref{eqDefiPHImaiuscolo}) of $\Phi$ and we recall that  $s>\ztn+h$ so that the integrand is continuous and we have
  \[
\lim _{h\to 0^+} \dfrac{1}{h}\Phi(s;\ztn,h)= K(s-\ztn) \int _{\ztn}^T
K^*(\nu-\ztn)\phi(\nu,\ztn )\ZD \nu\qquad \forall s>\ztn 
  \]
($\phi(\nu,\ztn) $ is defined in~(\ref{eqDIFEdIphistheta}): $\phi(\cdot,\ztn)=H_{\ztn} h_\ztn$).

 We extend the functions in both the side of the previous equality with zero to $[\zt_0,\ztn]$ and, analogously, we extend with zero the function $s\mapsto [H_{\ztn+h}h_{\ztn+h} ](s)$   so that the integral 
 in~(\ref{subequatioPERdopobmzaa})  can be consider as an integral on the fixed interval $[\zt_0,T]$.

Pointwise convergence is preserved by the extension.
  
  We prove below the existence of $p>1$ and $M$ such that 
  \begin{equation}
  \ZLA{INEdaPROVAREperZAA}
  \left \|\dfrac{1}{h}\Phi(\cdot;\ztn,h)\right \|_{L^p(\ztn+h,T;H)} <M\,.
  \end{equation}

 We extend $\Phi(\cdot;\ztn,h)$ with zero to $(\zt_0,T)$ so that the previous inequality holds for the extension too.

  The \emph{strict inequality $p>1$} implies compactness of the set $\{\Phi(\cdot;\ztn,h)\}_{\ztn\geq\zt_0,h\in(0,1]}$ in the space $L^p(\zt_0,T;H)$. 
  
Mazur's theorem shows that   any weak limit point is equal to the pointwise limit of a subsequence. Unicity of the pointwise limit then implies that there is only one     weak  limit point, i.e. weak convergence:
  \[
{\rm w-}\lim _{h\to 0^+}  \Phi(\cdot;\ztn,h)=
 K(\cdot-\ztn) \int _{\ztn}^T
K^*(\nu-\ztn)\phi(\nu,\ztn )\ZD \nu\,.
  \]

 This observation and~Lemma~\ref{LemmaINAppeDeriValue} imply   the existence of the limit we are looking for: 
\begin{multline} \ZLA{EqLIMITEesplicitoDIsubequatioPERdopoALPHA}
\lim _{h\to 0^+}\mbox{(\ref{subequatioPERdopobmzaa})}=
 \int_ \ztn^T \Bigl\ZL
K(s-\ztn) \int _{\ztn}^T
K^*(\nu-\ztn)\phi(\nu,\ztn )\ZD \nu,
\left [H_{\ztn }h_{\ztn } \right ](s)
\Bigr\ZR_H\ZD s\\
=\Bigl \ZL \int_\ztn ^T K^*(\nu-\ztn)\phi(\nu,\ztn)\ZD\nu, \int_\ztn^T K^*(s-\ztn)[H_\ztn h_\ztn](s)\ZD s\Bigr\ZR_{H}\\
= \ZL [\Lambda_\ztn^*\phi(\cdot,\ztn)](\ztn), [\Lambda_\ztn^*H_\ztn h_\ztn](\ztn)\ZR_H=
 \| [\Lambda_\ztn^*H_\ztn h_\ztn](\ztn)\|^2_H\,.
\end{multline}
%%%%%%%%%%%%%
 
 An analogous computation gives the same value for the limit of~(\ref{subequatioPERdopobmzaa}'):
 \[
 \lim _{h\to 0^+}\mbox{(\ref{subequatioPERdopobmzaa})}=\lim _{h\to 0^+}\mbox{(\ref{subequatioPERdopobmzaa}')}\,.
 \]
 
In order to complete the previous arguments, we prove~(\ref{INEdaPROVAREperZAA}) in a suitable $p$-norm, $p>1$.

First we note
\begin{equation}
\ZLA{INEdaPROVAREperZAABIS}
\|
\phi(\nu,\ztn )\|<M\qquad \zt_0\leq\ztn\leq\nu\leq T\,.
\end{equation}
The number $M$ depends solely on $T$. This inequality follows from the definition~(\ref{eqDIFEdIphistheta}) of $\phi(\nu,\ztn)$.
 
Inequality~(\ref{INEdaPROVAREperZAA}) follows from~(\ref{INEdaPROVAREperZAABIS}) and from the following chain of inequalities:
 
\begin{multline*}
\left \|
  \dfrac{1}{h} \int _{\ztn} ^{\ztn+h} K(\cdot-r)\int _r^T K^{*} (\nu-r)h_{\ztn}(\nu )\ZD\nu\,\ZD r 
 \right \|^p_{L^p(\ztn+h,T;H)}\\
 \leq
  M\int _{\ztn+h}^T \left [
 \dfrac{1}{h} \int_{\ztn} ^{\ztn+h} \|K(s-r)\|\ZD r\right ]^p\ZD s 
\leq 
M\int _{\ztn+h}^T \left [
\dfrac{1}{h} \int_{\ztn} ^{\ztn+h} \dfrac{1}{(s-r)^{\ZSI}}\ZD r\right ]^p\ZD s \\
%=
%M\int _{\ztn+h}^T \dfrac{1}{(s-r(s;h))^{p\ZSI  }}\ZD s\qquad \mbox{(where $0<r(s,h)< \ztn+h<s$)}\\
\leq M\int _{\ztn+h}^T \dfrac{1}{(s-(\ztn+h))^{p\ZSI}}<+\ZIN\qquad 
\mbox{provided that $1\leq p<1/\ZSI$}\,. 
\end{multline*}  
 
\paragraph{The limits of~(\ref{subequatioPERdopobmbeta}) and~(\ref{subequatioPERdopobmgamma}) and those   of~(\ref{subequatioPERdopobmbeta}') and~(\ref{subequatioPERdopobmgamma}')}

 We study the limit of~(\ref{subequatioPERdopobmbeta}).
 The integrand is 
 \begin{equation}\ZLA{IntegrandDIsubequatioPERdopobmbeta}
\Bigl\ZL(-A)^{\ZEP}\phi(\cdot;\ztn), \dfrac{(-A)^{-\ZEP}\Gamma_{\ztn+h}\Xi(\ztn+h)-(-A)^{-\ZEP}\Gamma_{\ztn}
\Mh\Xi(\ztn+h) 
}{h} \Bigr\ZR_H
 \end{equation}
 where $[H_{\ztn} h_\ztn](s)=\phi(s;\ztn)$ is the function in~(\ref{eqDIFEdIphistheta}). It solves the Fredholm integral equation given by the system~(\ref{Eq:appeDaPSIaPHI})-(\ref{Eq:appeDaPHIaPSI}) with 
the continuous affine term $g=H_{\ztn} h_\ztn$ i.e.
 \[
 \phi (s;\ztn)=-\int_{\ztn} ^s K(s-r)\int_r^T K^*(\nu-r)\phi (\nu;\ztn)\ZD\nu\,\ZD r+[H_\ztn h_\ztn] (s)\,.
 \]
The affine term   $H_\ztn h_\ztn$ is continuous and, from the statements~\ref{I5TeoPropriZt} and~\ref{I4TeoPropriZt} of Theorem~\ref{TeoPropriZt}, there exists $\ZEP>0$ such that $(-A)^\ZEP h_\ztn$ is continuous too
and such that $(-A)^\ZEP K(t)$ is of class $L^p$ with an exponent $p>1$.
  It follows from Theorem~\ref{teoAPPEregolarPRPsoluFred} that
  
 \[
s\mapsto (-A)^{\ZEP }\phi (s;\ztn)\quad\mbox{is $H$-valued continuous}\,.
 \]
 
In order to study the incremental quotient in~(\ref{IntegrandDIsubequatioPERdopobmbeta})  we use the operators in~(\ref{EqDEfiGamma12}).
 The incremental quotient is the sum of the two incremental quotients \begin{subequations}
\begin{align}
\ZLA{eqDERIriccaCaso2A}
& \dfrac{(-A)^{-\ZEP}\Gamma^1_{\ztn+h}\Xi(\ztn+h)-(-A)^{-\ZEP}\Gamma^1_{\ztn}
\Mh \Xi(\ztn+h) 
}{h} &s\in (\ztn+h,T)\,,\\[2mm]
\ZLA{eqDERIriccaCaso2B}
& \dfrac{(-A)^{-\ZEP}\Gamma^2_{\ztn+h}\Xi(\ztn+h)-(-A)^{-\ZEP}\Gamma^2_{\ztn}
\Mh \Xi(\ztn+h) 
}{h}&s\in (\ztn+h,T)\,.
\end{align}
\end{subequations}

\begin{Remark}\ZLA{Rema1SuallaConveeqPERdopoBETAeDopoGAMMA}
{\rm
Now we recall: $\Xi(\ztn+h)=\mathcal{C}^{(1,2)}_{\ztn+h}\Sta(\ztn+h;\zt_0,\Sta_{\zt_0},u)$ with $\Sta_{\zt_0}\in \Dom\,\mathcal{A}_{\zt_0}$ and $u$ continuous. But, the fact that $\Xi(\ztn+h)$ is the $M^2$-component of the evolution of the state is not used in the following computations. In order to appreciate  the computations below and so the definition of $P'$, it is important to note that the computations are justified 
under the sole conditions concerning the $M^2$-component $\Xi$ in the Assumption~\ref{AssuSuallaConveeqPERdopoBETAeDopoGAMMA}.

}
\end{Remark}

We state explicitly:
\begin{Lemma}\ZLA{LemmaConseI2Rema1SuallaConveeqPERdopoBETAeDopoGAMMA}
Let the $H$-valued function $v$, defined on $[0,T]$, 
have the properties   in the Assumption~\ref{AssuSuallaConveeqPERdopoBETAeDopoGAMMA}. 
We have:
\begin{enumerate}
\item \ZLA{I1LemmaConseI2Rema1SuallaConveeqPERdopoBETAeDopoGAMMA}  we have
\[
\lim _{h\to 0^+} v(\ztn+h)=v(\ztn)\quad \mbox{in $H$}\,,
\qquad \lim _{h\to 0^+}\xi_{ \ztn+h} =\xi_{ \ztn }
\quad \mbox{in $C([0,T];H)$} 
  \,.
\] 
\item\ZLA{I2LemmaConseI2Rema1SuallaConveeqPERdopoBETAeDopoGAMMA}
We have
\begin{multline}\ZLA{EqI1LemmaConseI2Rema1SuallaConveeqPERdopoBETAeDopoGAMMA} 
\lim _{h\to 0}\int_0^{t}e^{- \nu}   \dfrac{   v(t-(\nu-h))- v(t-\nu)}{h}\ZD \nu\\
=
\int_0^{t}e^{- \nu} (-\partial_\nu)v(t-\nu)\ZD\nu \in  C([0,T];H) \,.
\end{multline}
 The convergence is  in $C([0,T];H)$.
\end{enumerate}
\end{Lemma}
\zProof
 Property~\ref{I1LemmaConseI2Rema1SuallaConveeqPERdopoBETAeDopoGAMMA} is obvious. Property~\ref{I2LemmaConseI2Rema1SuallaConveeqPERdopoBETAeDopoGAMMA} 
 is seen because the integral in the left side of~(\ref{EqI1LemmaConseI2Rema1SuallaConveeqPERdopoBETAeDopoGAMMA})
 is
 \begin{multline*}
-\dfrac{1}{h}\int_0^h e^{-(t-s)}v(s)\ZD s+\dfrac{1}{h} \int_t^{t+h}e^{-(t+h-s)}v(s)\ZD s\\+\int_h^t \dfrac{e^{-(t+h-s)}-e^{-(t-s) }}{h}v(s)\ZD s\,.\zdiaform
 \end{multline*}
 
 The incremental quotient~(\ref{eqDERIriccaCaso2A}) is
 \[
\dfrac{(-A)^{-\ZEP}Z(s-\ztn+h)v(\ztn+h)-
(-A)^{-\ZEP}Z(s-\ztn)v(\ztn+h)}{h}\qquad s\in[\ztn+h,T] \,.
 \]
As $v(\ztn+h)\to v(\ztn)$ for $h\to 0^+$,   uniformly in $ \ztn\in[\zt_0,T]$,  inequality~(\ref{equaDELTeo:limRappINCREz}) implies that
\[
(-A)^{-\ZEP}\dfrac{Z(s-\ztn+h)-Z(s-\ztn)}{h} [v(\ztn+h)-v(\ztn)] 
\]
is bounded in a space $L^p$ and converges   to zero. So, the following limit exists in an $L^p$-norm:
%\begin{subequations}
\begin{multline}\ZLA{eqDERIriccaCaso2Arisultato}
\lim _{h\to 0^+}\mbox{(\ref{eqDERIriccaCaso2A})} 
=-(-A)^{-\ZEP} Z'(\cdot-\ztn)v(\ztn)=-(-A)^{-\ZEP} Z'(\cdot-\ztn) \hat  v_\ztn\\
=-(-A)^{-\ZEP}\left [
(A+I)Z(\cdot-\ztn) \hat v_\ztn -\int_0^{\cdot-\ztn}e^{-(\cdot-\ztn-r)}Z(r) \hat v_\ztn\ZD r\right ]\,. 
\end{multline}

 We study the limit of~(\ref{eqDERIriccaCaso2B}). 
 We observe:~\emph{(\ref{eqDERIriccaCaso2B}) does not use the incremental quotient of $v$} and we can write
\begin{align*}
&\mbox{(\ref{eqDERIriccaCaso2B})}=\fbox{A$_h$}+\fbox{B$_h$}\\
\quad &\fbox{A$_h$}=\dfrac{1}{h}\int_{\ztn}^{\ztn+h} Z(s-r) e^{-(r-\ztn)}\int_0^\ztn e^{-\nu}v(\ztn+h-\nu)\ZD \nu\,\ZD r\\
\quad &\fbox{B$_h$}=-\int _{\ztn+h}^\cdot Z(\cdot -r)e^{-(r-\ztn)}\dfrac{1}{h}
\left [
 e^{h} \int_0^{\ztn+h} e^{-\nu}v(\ztn+h-\nu)\ZD \nu\right.\\
&\hskip 5.7cm   \left.-\int_0^{\ztn} e^{-\nu} v(\ztn+h-\nu)\ZD\nu
\right ]  \ZD r
\end{align*}

It is clear
\begin{equation}
\ZLA{LimAh}
\lim _{h\to 0^+}\fbox{A$_h$} =Z(s-\ztn)\int_0^{\ztn}e^{-\nu} v(\ztn-\nu)\ZD\nu
\end{equation}
because $v\ZCD$ is continuous.  
 
To elaborate $  \fbox{B$_h$} $ first we note
\begin{multline}\ZLA{eqPerPrimelabBh}
\dfrac{1}{h}
\left [
 e^{h} \int_0^{\ztn+h} e^{-\nu}v(\ztn+h-\nu)\ZD \nu -\int_0^{\ztn} e^{-\nu} v(\ztn+h-\nu)\ZD\nu
\right ]   \\
=\dfrac{1}{h} \int_{-h}^0
e^{-\nu}v(\ztn-\nu)\ZD\nu 
+ \int_0^{\ztn} e^{-\nu}\dfrac{1}{h}\left [v(\ztn-\nu)-v(\ztn-\nu+h)\right ]\ZD \nu\,.
 \end{multline}
 \emph{Note that the incremental quotient of $v$ in this expression  is inherited from the incremental quotient of the operator.}

 Thanks to the assumption $\Xi_{\ztn}\in\Dom\,\tilde{\mathcal{A}}_\ztn$ we have
\[% \begin{equation}\ZLA{conPerRippincrDebFIN}
\lim _{h\to 0^+}\mbox{(\ref{eqPerPrimelabBh})}= \hat  v_\ztn 
   -\int_0^{\ztn}e^{-\nu}(-\partial _\nu) v(\ztn-\nu)\ZD\nu \,.
\]% \end{equation}
 It follows:
 
\begin{multline}\ZLA{eq:LImiAhPiuBh}
\lim_{h\to 0^+}(-A)^{-\ZEP}
 \fbox{A$_{h}$}+
\lim_{h\to 0^+}(-A)^{-\ZEP}
 \fbox{B$_{h}$}  \\
=(-A)^{-\ZEP}\left [
Z(\cdot-\ztn)\int _0^{\ztn} e^{-\nu}v(\ztn-\nu)\ZD s\right.\\
\left.-\int _{\ztn}^\cdot
Z(\cdot-r)e^{-(r-\ztn)}  \hat  v_{\ztn}\ZD r
+\int_\ztn^\cdot Z(\cdot-r)e^{-(r-\ztn)}\int_0^{\ztn} e^{-\nu}(-\partial_\nu)v(\ztn-\nu)\ZD\nu\,\ZD r
\right ] \,.
\end{multline}
We note
\[
\int_\ztn^s Z(s-r)e^{-(r-\ztn)}  \hat  v_{\ztn}\ZD r=\int_0^{s-\ztn} Z(\nu)e^{-(s-\nu-\ztn)}  \hat  v_{\ztn}  \ZD \nu\,.
\]
This integral is equal, with opposite sign, to the corresponding integral in~(\ref{eqDERIriccaCaso2Arisultato}).

In conclusion, by summing~(\ref{eq:LImiAhPiuBh}) with~(\ref{eqDERIriccaCaso2Arisultato}) and by recalling $\xi_\ztn\ZCD=v(\ztn-\cdot)$, we get
 
\begin{multline}\ZLA{FineLIMITsubequatioPERdopobmbeta)}
\lim _{h\to 0^+}\mbox{(\ref{subequatioPERdopobmbeta})}=
 \int_\ztn^T\Bigl\ZL
 H_{\ztn} h_\ztn, \left [-Z(\cdot-\ztn)\left ((A+I) \hat v_\ztn-\int_0^{\ztn} e^{-\nu}  \xi_\ztn(\nu)\ZD\nu\right )\right.\Bigr.\\
\Bigl. \left. +\int_0^\cdot Z(\cdot-r)e^{-(r-\ztn)}\int_0^{\ztn} e^{-\nu}(-\partial_\nu)  \xi_\ztn(\nu)\ZD\nu\,\ZD r\right ]\Bigr\ZR\ZD s\\
= \int_\ztn^T\Bigl\ZL
 H_{\ztn} h_\ztn, \Gamma_\ztn\left ( - \tilde {\mathcal A}_{\ztn}\Xi_{\ztn}\right )\Bigr\ZR\ZD s\,.
\end{multline}
 
By invoking Lemma~\ref{LemmaINAppeDeriValue} 
we see that
\[
\lim _{h\to 0^+}(-A)^{\ZEP} H_{\ztn+h} h_{\ztn+h}=(-A)^{\ZEP} H_{\ztn } h_{\ztn }
\] 
boundedly and in $C([\ztn+h_0,T];H)$ (any $h_0>0$. So
have also
\begin{equation}\ZLA{FineLIMITrefsubequatioPERdopobmgamma)}
\lim _{h\to 0^+}\mbox{(\ref{subequatioPERdopobmgamma})}
=
 \int_\ztn^T\Bigl\ZL\Gamma_\ztn\left ( - \tilde {\mathcal A}_{\ztn}\Xi_{\ztn}\right ),
 H_{\ztn} h_\ztn 
 \Bigr\ZR\ZD s\,.
\end{equation}

Similar computations give
\[
\lim _{k\to 0^+}\mbox{(\ref{subequatioPERdopobmbeta}')}=
\lim _{h\to 0^+}\mbox{(\ref{subequatioPERdopobmbeta})}\qquad
\lim _{k\to 0^+}
 \mbox{(\ref{subequatioPERdopobmgamma}')}= \lim _{h\to 0^+}
  \mbox{(\ref{subequatioPERdopobmgamma})}\,.
\]
\paragraph{The contribution of the terms~(\ref{eqINCRequoDOLOriccati})
and of~(\ref{subequatioPERdopobmSOLOoperatore})} Let
\[\Sta_\ztn
=\left[\begin{array}{c}
 \hat v_\ztn\\
  \xi_\ztn\ZCD\\
\hat y_\ztn
\end{array}\right]\,.
\]
We sum~(\ref{FineLIMITsubequatioPERdopobmbeta)}),~(\ref{FineLIMITrefsubequatioPERdopobmgamma)})
  and~(\ref{subequatioPERdopobmSOLOoperatore}). We find
  \begin{multline}  \ZLA{LasommaDIPARTEderiRiccaOper}
  \mbox{(\ref{FineLIMITsubequatioPERdopobmbeta)})}+\mbox{(\ref{FineLIMITrefsubequatioPERdopobmgamma)})}+
\mbox{(\ref{subequatioPERdopobmSOLOoperatore})}=
  %%%%
  \ZL 
   H_\ztn h_\ztn  , 
 \Gamma_\ztn\mathcal{C}^{(1.2)}_\ztn(-\mathcal{A}_\ztn\Sta_\ztn
  \ZR _{L^2(\ztn,T;H)}\\
  +
   \ZL 
%%%% 
 \Gamma_\ztn\mathcal{C}^{(1.2)}_\ztn(-\mathcal{A}_\ztn\Sta_\ztn),  H_\ztn h_\ztn
  \ZR _{L^2(\ztn,T;H)}\,.
  \end{multline}
%\bibliographystyle{plain}
%\bibliography{bibliomemoria}

\enddocument